\documentclass[final,hidelinks,onefignum,onetabnum]{siamart251216}

\newsiamremark{remark}{Remark}
\newsiamremark{example}{Example}
\newsiamremark{property}{Property}
\crefname{example}{Example}{Examples}
\crefname{property}{Property}{Properties}
\AddToHook{env/example/begin}{\crefalias{theorem}{example}}
\AddToHook{env/property/begin}{\crefalias{theorem}{property}}
\AddToHook{env/remark/begin}{\crefalias{theorem}{remark}}

\newcounter{step}
\newcounter{substep}[step]
\newcounter{stepdepth}
\renewcommand{\thestep}{\arabic{step}}
\renewcommand{\thesubstep}{\thestep.\arabic{substep}}
\newenvironment{step}{%
	\par\addvspace{0.75em}%
	\ifnum\value{stepdepth}=0
		\refstepcounter{step}%
		\setcounter{substep}{0}%
		\def\currentsteplabel{\thestep}%
	\else
		\refstepcounter{substep}%
		\def\currentsteplabel{\thesubstep}%
	\fi
	\addtocounter{stepdepth}{1}%
	\noindent\textbf{Step~\currentsteplabel.}\enspace\ignorespaces
}{%
	\par\addtocounter{stepdepth}{-1}\addvspace{0.75em}%
}

\usepackage{amssymb}
\usepackage{amsfonts}
\usepackage{mathrsfs}
\usepackage{physics2}
\usephysicsmodule{ab, ab.legacy}
\usephysicsmodule{qtext.legacy}
\usephysicsmodule{op.legacy}
\usephysicsmodule{nabla.legacy}
\usepackage{derivative, fixdif}
\usepackage{newunicodechar}
\newunicodechar{∇}{\ensuremath{\nabla}}
\newunicodechar{·}{\ensuremath{\cdot}}
\usepackage{enumerate}

\usepackage{algpseudocode}

\crefname{step}{Step}{Steps}
\crefname{equation}{Eq.}{Eqs.}

\usepackage{indentfirst}
\usepackage{tabularx}

\usepackage[exponent-mode = scientific]{siunitx}

\usepackage{caption}
\DeclareCaptionFormat{siamtable}{%
  \centerline{#1}\par
  \setbox0\hbox{#3}%
  \ifdim\wd0>\hsize
    #3%
  \else
    \centerline{#3}%
  \fi}
\usepackage{subcaption}
\usepackage{float} 

\usepackage{booktabs} 
\usepackage{multirow}
\usepackage{cite} 

\def\hwenobase{HWENO}
\def\hwenodf{DFHWENO}

\usepackage{tikz}
\usetikzlibrary{arrows.meta, positioning, shapes.geometric, fit, backgrounds}

\usepackage{changes}
\definechangesauthor[color=teal]{chen}
\definechangesauthor[color=red]{zhao}

\def\half{\frac{1}{2}}
\def\dx{\Delta x}
\def\dy{\Delta y}
\def\dt{\Delta t}

\newcommand{\divh}[1]{\nabla_{h} \cdot #1}
\def\cleaneddivB{\widehat{\divh{\theB}}}
\def\divB{\div \theB}
\def\deg{\mathrm{deg}\,}
\usepackage{bm}
\def\vec#1{\bm{#1}}

\def\primitives{
	(\rho, \vec{u}^{\top}, \theB^{\top}, p)
}

\def\theU{\mathbf{U}}
\def\theV{\mathbf{V}}
\def\theW{\mathbf{W}}
\def\theF{\mathbf{F}}
\def\theG{\mathbf{G}}
\def\theB{\mathbf{B}}

\headers{Divergence-free HWENO scheme for ideal MHD}{P. Chen, S. Jin, J. Qiu, and Z. Zhao}

\title{A fifth-order divergence-free finite difference Hermite WENO scheme for ideal magnetohydrodynamics\thanks{
		\funding{This work was partially supported by the National Natural Science Foundation of China under Grant Number 12401541, the Natural Science Foundation of Xiamen, China under Grant Number 3502Z202472004, and the Tianyuan Fund for Mathematics of the National Natural Science Foundation of China under Grant Number 12426304.
	}
}}

\author{Peiwen Chen\thanks{School of Mathematical Sciences, Xiamen University, Xiamen, Fujian 361005, P.R. China
  (\email{chenpw@stu.xmu.edu.cn}).}
\and Shi Jin\thanks{School of Mathematical Sciences, Institute of Natural Sciences, MOE-LSC, Shanghai Jiao Tong University, Shanghai, 200240, P.R. China
  (\email{shijin-m@sjtu.edu.cn}).}
\and Jianxian Qiu\thanks{School of Mathematical Sciences and Fujian Provincial Key Laboratory of Mathematical Modeling and High-Performance Scientific Computing, Xiamen University, Xiamen, Fujian 361005, P.R. China
  (\email{jxqiu@xmu.edu.cn}).}
\and Zhuang Zhao\thanks{Corresponding author. School of Mathematical Sciences and Fujian Provincial Key Laboratory of Mathematical Modeling and High-Performance Scientific Computing, Xiamen University, Xiamen, Fujian 361005, P.R. China
  (\email{zzhao@xmu.edu.cn}).}}

\ifpdf
\hypersetup{
  pdftitle={A fifth-order divergence-free finite difference Hermite WENO scheme for ideal magnetohydrodynamics},
  pdfauthor={P. Chen, S. Jin, J. Qiu, and Z. Zhao}
}
\fi

\begin{document}

\maketitle

\begin{abstract}
In this paper, we present a fifth-order finite difference divergence-free Hermite weighted essentially non-oscillatory (HWENO) scheme for the ideal magnetohydrodynamics (MHD) equations.
In this framework, both the solution and its spatial partial derivatives are evolved in time and jointly employed in the spatial reconstruction procedure.
A major challenge in MHD simulations is preserving the divergence-free constraint of the magnetic field, which is generally violated by standard numerical methods designed solely for hyperbolic conservation laws.
To address this issue, we first solve the MHD equations within the HWENO framework for hyperbolic conservation laws, yielding a magnetic field divergence that remains zero up to high-order accuracy in smooth regions.
Subsequently, we apply a correction that evenly distributes the divergence error among the partial derivatives involved in the divergence-free constraint, thereby rendering the magnetic field discretely divergence-free at the new time level.
This approach offers several advantages.
First, the scheme retains the conservation property, as only the partial derivatives of the numerical solution are corrected, leaving the conserved variables unchanged.
Second, the correction applied to the partial derivatives of the magnetic field components introduces only a high-order perturbation, thereby preserving the overall accuracy.
Third, the divergence-free treatment significantly enhances robustness, as most benchmark test cases 
cannot be run stably without such a correction.
Fourth, the correction is a simple linear operation applied at each stage of the time integration, incurring negligible additional computational cost. Extensive numerical experiments demonstrate the accuracy, resolution, efficiency, effectiveness, and robustness of the proposed scheme.

\end{abstract}

\begin{keywords}
MHD, divergence-free, HWENO, finite difference, high-order accuracy
\end{keywords}

\begin{MSCcodes}
65M06, 
76W05, 
35L65 
\end{MSCcodes}

\section{Introduction}
This paper is concerned with the implementation of 
the fifth-order accurate divergence-free finite difference HWENO scheme
for the ideal compressible MHD equations.
The governing equations of ideal MHD read
\begin{equation}
	\label{eq:mhd}
	\pdv*{}{t}
	\begin{bmatrix}
	    		\rho \\
				\rho\vec{u}\\
				\theB\\
				\mathcal{E}
	\end{bmatrix}
	+\div \begin{bmatrix}
		\rho\vec{u}\\
		\rho\vec{u}\otimes\vec{u}+(p+\frac{1}{2}\norm{\theB}^{2})\mathbb{I}-\theB\otimes\theB\\
		\vec{u}\otimes\theB-\theB\otimes\vec{u}\\
		(\mathcal{E}+p+\frac{1}{2}\norm{\theB}^{2})\vec{u}-(\vec{u}\cdot\theB)\theB
	\end{bmatrix}
	= 0,
\end{equation}
where $\rho$ is the density, $\vec{u} = (u^{1}, u^{2}, u^{3})^{\top}$ is the velocity, $\theB = (B^{1}, B^{2}, B^{3})^{\top}$ is the magnetic field, $\mathcal{E} = \rho e + \frac{1}{2}\rho \norm{\vec{u}}^{2} + \frac{1}{2}\norm{\theB}^{2}$ is the total energy, and $ \rho e $ is the internal energy. The system is closed by the equation of state
for an ideal gas, $ p = (\gamma - 1)\rho e $, with the adiabatic index $ \gamma. $
Due to the absence of magnetic monopoles, the magnetic field must satisfy the divergence-free condition \begin{equation*}
    \div \theB = 0.
\end{equation*}
In numerical simulations of ideal MHD, preserving the divergence-free constraint is crucial,
since violations of this constraint can lead to nonphysical solutions and numerical instabilities \cite{Brackbill1980EffectNonzero}. 

The ideal MHD equations can be viewed as a hyperbolic system of conservation laws
together with the divergence-free constraint of the magnetic field.
Therefore, it is natural to consider applying high-order numerical methods for hyperbolic conservation laws directly.
Over the past few decades, various numerical methods have been developed for hyperbolic conservation laws, including the essentially non-oscillatory (ENO) scheme \cite{Harten1987UniformlyHighOrder}, the weighted ENO (WENO) schemes \cite{Liu1994WeightedEssentially, Jiang1996EfficientImplementation}, the discontinuous Galerkin (DG) method \cite{Cockburn1989TVBRungeKutta}, and the Hermite WENO (HWENO) schemes \cite{Qiu2004HermiteWENO, Qiu2005HermiteWENO, Liu2014FiniteDifference}.
HWENO schemes can be regarded as an extension of WENO schemes, but the first-order partial derivatives of the solution are additionally introduced and evolved in time together with the solution, and are used in the spatial reconstruction procedure to achieve high-order accuracy with more compact stencils than WENO schemes.
This compactness makes the HWENO schemes more computationally efficient with smaller numerical errors in smooth regions and enables higher resolution near discontinuities than WENO schemes of the same order \cite{Fan2023RobustFifth},
although the evolution of partial derivatives also introduces implementation complexity and incurs additional computational cost.

However, these numerical methods that have been successfully developed for  hyperbolic conservation laws cannot be directly applied to ideal MHD, since they
generally do not preserve the divergence-free condition.
For the continuous ideal MHD system, the divergence-free constraint is automatically preserved by analytical solutions if the initial magnetic field is divergence-free. 
However, numerical methods produce numerical errors in the magnetic field, and the divergence-free constraint is generally not preserved
exactly, or even approximately in a discrete sense.
This is in the same spirit as the well-balanced (WB) schemes for hyperbolic balance laws \cite{Xing2006HighOrder, Jin2001SteadystateCapturing, Jin2016WellBalancedStochastic}: the steady-state balance condition and the divergence-free constraint are both structural properties automatically preserved by the continuous equations, yet easily destroyed by naive discretizations, and both must be deliberately built into the numerical scheme to avoid spurious artifacts.
In this context, one of the major challenges in the numerical simulation of ideal MHD is to design effective strategies for handling the divergence-free constraint.
Over the past few decades, various techniques have been proposed to handle the divergence-free constraint in ideal MHD.
The projection method \cite{Brackbill1980EffectNonzero} enforces this constraint by projecting the magnetic field onto a divergence-free space through the solution of a Poisson equation for a scalar potential.
The eight-wave method \cite{Powell1997ApproximateRiemann, Powell1999SolutionAdaptiveUpwind, Powell1995UpwindScheme} augments the MHD system with the Godunov--Powell source term
to make the system symmetrizable.
Although this method does not enforce the divergence-free condition exactly, it allows divergence errors to be advected with the flow, thereby preventing their accumulation in the computational domain. However, the nonconservative nature of the eight-wave formulation may lead to certain drawbacks.
The constrained transport (CT) method \cite{Evans1988SimulationMagnetohydrodynamic, Gardiner2005UnsplitGodunov, Londrillo2004DivergencefreeCondition, Xu2016DivergenceFreeWENO} employs staggered grids to preserve the divergence-free condition at the discrete level, and unstaggered CT methods have also been developed in \cite{Christlieb2014FiniteDifference, Christlieb2015PositivityPreservingFinite, Christlieb2016HighorderPositivitypreserving, Christlieb2018HighOrderFinite, Rossmanith2006UnstaggeredHighresolution, Cakir2021HighOrder, Yu2020FreeStreamPreserving, Helzel2011UnstaggeredConstrained}.
The hyperbolic divergence cleaning method \cite{Dedner2002HyperbolicDivergence} introduces an additional scalar field together with an evolution equation to propagate and damp out divergence errors.
The locally divergence-free DG method \cite{Li2005LocallyDivergencefree} constructs basis functions for the magnetic field that are exactly divergence-free within each element.
Globally divergence-free DG methods have also been developed \cite{Li2011CentralDiscontinuous, Fu2018GloballyDivergenceFree}. The idea of directly reconstructing the magnetic field to satisfy the divergence-free constraint was introduced by Balsara and Spicer \cite{Balsara1999StaggeredMesh} and further developed by Balsara \cite{Balsara2004SecondOrderaccurateSchemes, Balsara2009DivergencefreeReconstruction}.
Recently, Ding and Wu \cite{Ding2024NewDiscretely} proposed a new discrete divergence-free projection approach, which projects the reconstructed point values of the magnetic field at cell interfaces onto a discretely divergence-free space without solving additional partial differential equations.

However, none of the aforementioned methods directly evolve the spatial derivatives of the magnetic field, while these derivatives directly appear in the expression of the divergence-free constraint. In contrast, the HWENO framework evolves these derivatives as auxiliary variables, making it naturally suitable for ideal MHD and its divergence-free constraint. 
For two-dimensional problems, the divergence of the magnetic field is given by the sum of the first-order spatial partial derivatives of the components of $\theB$:
$ \divB = \pdv*{}{x}B^{1} + \pdv*{}{y}B^{2}.$
These partial derivatives of the magnetic field $ \pdv*{}{x}B^{1} $ and $ \pdv*{}{y}B^{2} $
are conveniently available in the HWENO framework since the partial derivatives of all conserved quantities, including the magnetic field, are directly evolved by the HWENO schemes \cite{Qiu2004HermiteWENO, Qiu2005HermiteWENO, Liu2014FiniteDifference, Zhao2020ModifiedFifth}.
Motivated by this observation, we propose a novel discretization of the magnetic field divergence, based on which we design a divergence-free correction for the magnetic field derivatives
and develop a divergence-free HWENO (DFHWENO) scheme within the finite difference framework.

The main contributions and novelties of this work are summarized as follows.
\begin{enumerate}[I.]
\item \textbf{We propose a novel discretization of the magnetic field divergence using partial derivatives naturally incorporated in the HWENO framework.}
	Specifically, at each grid point $(x_i, y_j)$, the discrete divergence is computed directly from the Hermite auxiliary variables as $\divh{\theB} = \theV_{i,j}^{5} + \theW_{i,j}^{6}$, where $\theV_{i,j}^{5}$ and $\theW_{i,j}^{6}$ are the evolved first-order spatial partial derivatives of the magnetic field components $B^{1}$ and $B^{2}$ in the $x$- and $y$-directions, respectively.
	This provides an effective and high-order discretization of the divergence and plays a fundamental role in the divergence-free correction for the magnetic field derivatives.

\item \textbf{We propose a novel, concise and efficient divergence-free correction technique.} This divergence-free correction technique is different in principle from traditional divergence-free treatments, and it is effective in handling the divergence-free constraint of the magnetic field, as demonstrated by the numerical experiments.
The solution and its partial derivatives are first evolved by the HWENO scheme, which yields a magnetic field divergence that remains zero up to high-order accuracy in smooth regions.
We compute the divergence of the magnetic field using the proposed discrete divergence, and then apply a correction that evenly distributes the divergence error to the first-order spatial partial derivatives of the magnetic field to ensure that the corrected divergence is exactly zero.
Since the correction only distributes the discrete divergence among the spatial partial derivatives, only high-order accurate perturbations are introduced to them in smooth regions, ensuring the high-order accuracy of the overall scheme.
Moreover, the correction is simple to implement and computationally efficient, since it is a linear procedure and does not require solving any additional equations or performing any complex operations. Furthermore, our method preserves the conservation property, since the correction 
procedure is only applied to the auxiliary variables rather than the physical variables, namely, the magnetic field itself.

\item \textbf{We apply the high-order HWENO scheme to ideal MHD for the first time.} The HWENO scheme is a high-order accurate, high-resolution, and compact-stencil method that has been successfully applied to various hyperbolic conservation laws. To the best of our knowledge, its application to ideal MHD has not yet been explored. Its extension to ideal MHD is nontrivial due to the divergence-free constraint of the magnetic field. The methods proposed in this work provide effective solutions to this challenge, and the numerical experiments further demonstrate the broad applicability of the HWENO framework in computational fluid dynamics.
\end{enumerate}

The remainder of this paper is organized as follows. \Cref{sec:algorithm} presents the algorithm of the proposed fifth-order finite difference DFHWENO scheme for the ideal MHD equations. In \cref{subsec:semi-discrete}, we describe the semi-discrete scheme. The discrete divergence based on the Hermite auxiliary variables is introduced in \cref{subsec:divergence}. The divergence-free correction procedure for the magnetic field derivatives is presented in \cref{subsec:correction}. The third-order SSP-RK3 temporal discretization is discussed in \cref{subsec:time-discretization}. \Cref{sec:examples} presents a series of numerical examples to demonstrate the accuracy, robustness, and effectiveness of the proposed scheme. Finally, \cref{sec:conclusions} concludes the paper.

\section{The divergence-free HWENO scheme}
\label{sec:algorithm}
\subsection{The semi-discrete scheme}
\label{subsec:semi-discrete}
We consider the two-dimensional ideal MHD equations \cref{eq:mhd}, which can be written as hyperbolic conservation laws  as
\begin{equation}
	\left\{
	\begin{aligned}
		 & \theU_{t} + \theF(\theU)_{x} + \theG(\theU)_{y} = 0, \quad          \\
    & \theU(\vec{x}, 0) = \theU_{0}(\vec{x}), 
	\end{aligned}
	\right.
	\label{eq:hyperbolic-conservation-laws}
\end{equation}
with the divergence-free constraint $\divB = \pdv*{B^{1}}{x}+\pdv*{B^{2}}{y}=0$. Here, the conserved variable $ \theU $ reads
\begin{equation*}
	\theU =
	\begin{bmatrix}
	 	\rho,
		\rho u^{1},
		\rho u^{2},
		\rho u^{3},
		B^{1},
		B^{2},
		B^{3},
		\mathcal{E}
	\end{bmatrix}^{\top},
\end{equation*}
and the fluxes $ \theF(\theU) $ and $ \theG(\theU) $ are given by
\begin{equation*}
	 \theF(\theU) = \begin{bmatrix}
		 \rho u^{1}\\
		 \rho u^{1} u^{1} + p^{\mathrm{tot}} - B^{1}B^{1} \\
		 \rho u^{1}u^{2} - B^{1}B^{2}\\
		 \rho u^{1}u^{3} - B^{1}B^{3}\\
		 0\\
		 u^{1}B^{2} - B^{1}u^{2}\\
		 u^{1}B^{3} - B^{1}u^{3}\\
		 (\mathcal{E}+p^{\mathrm{tot}})u^{1} - (\vec{u}\cdot\theB)B^{1}
	 \end{bmatrix},
	 \theG(\theU) = \begin{bmatrix}
	     		 \rho u^{2}\\
		 \rho u^{1}u^{2} - B^{1}B^{2}\\
		 \rho u^{2}u^{2} + p^{\mathrm{tot}} - B^{2}B^{2}   \\
		 \rho u^{2}u^{3} - B^{2}B^{3}\\
		 u^{2}B^{1} -B^{2} u^{1}\\
		 0\\
		 u^{2}B^{3} -B^{2} u^{3}\\
		 (\mathcal{E}+p^{\mathrm{tot}})u^{2} - (\vec{u}\cdot\theB)B^{2}
	 \end{bmatrix},
\end{equation*}
where $ p^{\mathrm{tot}} = p + \frac{1}{2}\norm{\theB}^{2} $ is the total pressure.
Note that the partial derivatives of the magnetic field $\pdv*{B^1}{x}$ and $\pdv*{B^2}{y}$ appear in the divergence-free constraint $\divB = \pdv*{B^{1}}{x}+\pdv*{B^{2}}{y}=0$. In order to directly evolve them, we differentiate \cref{eq:hyperbolic-conservation-laws} with respect to $ x $ and $ y $ as in the HWENO scheme \cite{Liu2014FiniteDifference}, to obtain  
\begin{equation}
	\left\{
	\begin{aligned}
		 & {\theU}_{t} + \theF(\theU)_{x} + \theG(\theU)_{y} = 0, \quad             \\
		 & {\theV}_{t} + \theF^{v}(\theU,\theV)_{x} + \theG^{v}(\theU,\theV)_{y} = 0, \quad \\
		 & {\theW}_{t} + \theF^{w}(\theU,\theW)_{x} + \theG^{w}(\theU,\theW)_{y} = 0. \quad \\
	\end{aligned}
	\right.
	\label{eq:2D-uvw-cartesian}
\end{equation}
where $ \theV = \pdv*{\theU}{x}  $, $ \theW = \pdv*{\theU}{y} $, $ \theF^{v}(\theU,\theV) = \theF'(\theU) \theV $, $ \theG^{v}(\theU,\theV) = \theG'(\theU) \theV $, $ \theF^{w}(\theU,\theW) = \theF'(\theU) \theW $, and $ \theG^{w}(\theU,\theW) = \theG'(\theU) \theW $.
Consider a uniform Cartesian grid 
\begin{equation*}
  \{(x_{i}, y_{j}) :i = 1, 2, \ldots, N_x ,  j = 1, 2, \ldots, N_y  \} ,   
\end{equation*}
with mesh size $ \dx = x_{i+1} - x_{i} $, $ \dy = y_{j+1} - y_{j} $, and half grid points $ x_{i+\half} = x_{i} + \dx $, $ y_{j+\half} = y_{j} + \dy $.
Then we obtain the  semi-discrete finite difference scheme for \cref{eq:2D-uvw-cartesian} as  
\begin{equation}
    \left\{
    \begin{aligned}
    	 & \odv*{\theU_{i,j}}{t} = -\frac{1}{\dx}(\hat{\theF}_{i+\half, j} - \hat{\theF}_{i-\half, j}) - \frac{1}{\dy}(\hat{\theG}_{i, j+\half} - \hat{\theG}_{i, j - \half}), \quad  \\
    	 & \odv*{\theV_{i,j}}{t} = -\frac{1}{\dx}(\hat{\theF}^{v}_{i+\half, j} - \hat{\theF}^{v}_{i-\half, j}) - \frac{1}{\dy}(\hat{\theG}^{v}_{i, j+\half} - \hat{\theG}^{v}_{i, j - \half}), \quad  \\
    	 & \odv*{\theW_{i,j}}{t} = -\frac{1}{\dx}(\hat{\theF}^{w}_{i+\half, j} - \hat{\theF}^{w}_{i-\half, j}) - \frac{1}{\dy}(\hat{\theG}^{w}_{i, j+\half} - \hat{\theG}^{w}_{i, j - \half}). \quad  \\
    \end{aligned}
    \right.
	\label{eq:semi-discrete}
\end{equation}
Here $ \hat{\theF}_{i+\half, j}, \hat{\theG}_{i, j+\half}, \hat{\theF}^{v}_{i+\half, j}, \hat{\theG}^{v}_{i, j+\half}, 
\hat{\theF}^{w}_{i+\half, j}$, and $\hat{\theG}^{w}_{i, j+\half} $ are the numerical fluxes, which can be
computed by the HWENO reconstruction procedure.
The details of the HWENO reconstruction  procedure can be found in \cite{Zhao2020ModifiedFifth, Liu2014FiniteDifference, Zhao2023WellbalancedFifthorder},
and we provide a brief review in Appendix~\ref{sec:hweno-detail}.
After the spatial discretization,
the strong stability preserving (SSP) Runge--Kutta method \cite{Gottlieb2001StrongStabilitypreserving} is used for the temporal discretization of the semi-discrete scheme \cref{eq:semi-discrete}.
It is detailed in \cref{subsec:time-discretization}.

The semi-discrete scheme \cref{eq:semi-discrete}, along with its temporal and spatial discretizations, follows a standard finite difference HWENO framework for hyperbolic conservation laws, as described in \cite{Zhao2020ModifiedFifth, Liu2014FiniteDifference, Zhao2023WellbalancedFifthorder}. However, directly applying this general HWENO  framework to  the MHD equations may lead to a nonzero divergence of the magnetic field, i.e., $\divB = \pdv*{B^{1}}{x}+\pdv*{B^{2}}{y}\neq 0$, and eventually 
to nonphysical solutions. We therefore next introduce a divergence-free correction technique based on a discrete divergence computed from $\theV^5$ (i.e., $\pdv*{B^{1}}{x}$) and $\theW^6$ (i.e., $\pdv*{B^{2}}{y}$),  both of which  evolve over time in the HWENO scheme.

\subsection{The divergence of the magnetic field}
\label{subsec:divergence}
We observe that, in the semi-discrete scheme \cref{eq:semi-discrete}, the auxiliary variables $\theV_{i,j}$ and $\theW_{i,j}$ represent the spatial derivatives of the conserved variable $\theU_{i,j}$ in the $x$- and $y$-directions, respectively. For the ideal MHD system
\begin{equation*}
	\theU =
	\begin{bmatrix}
	 	\rho,
		\rho u^{1},
		\rho u^{2},
		\rho u^{3},
		B^{1},
		B^{2},
		B^{3},
		\mathcal{E}
	\end{bmatrix}^{\top},
\end{equation*}
the fifth component of $\theV_{i,j}$ and the sixth component of $\theW_{i,j}$ correspond to the derivatives of the magnetic field components $B^{1}$ and $B^{2}$ in the $x$- and $y$-directions, respectively, namely
\begin{equation*}
		\theV_{i,j}^{5} \approx \pdv*{B^{1}(t, x_{i}, y_{j})}{x}, 
		\quad 
		\theW_{i,j}^{6} \approx \pdv*{B^{2}(t, x_{i}, y_{j})}{y}
\end{equation*}
Therefore, the magnetic field divergence can be naturally approximated by
\begin{equation*}
\divh{\theB} = \theV_{i,j}^{5} + \theW_{i,j}^{6}.
\end{equation*}
We call $\divh{\theB}$ the discrete divergence, since it is computed directly from the Hermite auxiliary variables $\theV$ and $\theW$, rather than from the physical magnetic field variables themselves.

\begin{property}
    The discrete divergence $ \divh{\theB} $ is a $ (k-1) $-th order approximation to $ \div \theB $ in smooth regions for the $ k $-th order HWENO scheme, namely
	\begin{equation*}
				\divh{\theB} = \div \theB + \mathcal{O}(h^{k-1}) = \mathcal{O}(h^{k-1}),
	\end{equation*}
	where $ h = \max\{\dx, \dy\} $.
\begin{proof}
    Since the auxiliary variables $ \theV $ and $ \theW $ are evolved by the HWENO scheme, they are high-order approximations to the spatial derivatives of $ \theU $ in smooth regions. In particular, we have
	\begin{equation*}
		\theV_{i,j}^{5} = \pdv*{B^{1}(t, x_{i}, y_{j})}{x} + \mathcal{O}(h^{k-1}),
		\quad 
		\theW_{i,j}^{6} = \pdv*{B^{2}(t, x_{i}, y_{j})}{y} + \mathcal{O}(h^{k-1}).
	\end{equation*}
Therefore, by summing the two components, we have
\begin{equation*}
	\divh{\theB} = \theV_{i,j}^{5} + \theW_{i,j}^{6} = \pdv*{B^{1}(t, x_{i}, y_{j})}{x} + \pdv*{B^{2}(t, x_{i}, y_{j})}{y} + \mathcal{O}(h^{k-1}) = \div \theB + \mathcal{O}(h^{k-1}).
\end{equation*}
Since the magnetic field is divergence-free, we have $\div \theB = 0$, and therefore $\divh{\theB} = \mathcal{O}(h^{k-1})$ in smooth regions.
\end{proof}
\end{property}

This discrete divergence has several advantages.
First, its expression is extremely simple and easy to implement: it only requires the summation of two already available auxiliary variables, and therefore introduces almost no additional computational cost. 
Second, it provides a high-order and compact approximation to $\div \theB$ in smooth regions, because the auxiliary variables $\theV$ and $\theW$ are themselves high-order compact approximations to the spatial derivatives of $\theU$ ensured by the HWENO scheme. Accordingly,  we next introduce the divergence-free correction technique for the magnetic field derivatives based on this discrete divergence.

\subsection{The divergence-free correction of magnetic field derivatives}
\label{subsec:correction}

The magnetic field divergence is not necessarily zero in numerical simulations,
which may lead to non-physical solutions and numerical instabilities.
In our scheme, if the auxiliary variables $ \theV_{i,j} $ and $ \theW_{i,j} $ are evolved 
by the general HWENO scheme without any treatment for the magnetic field divergence, 
the discrete divergence $ \divh{\theB} $  also becomes nonzero. 
This may further result in non-physical solutions and numerical instabilities, 
as demonstrated by the numerical experiments in \cref{sec:examples}.

To address this issue, we propose a divergence-free correction for the 
ideal MHD system within the finite difference HWENO framework. The proposed method 
evenly distributes the nonzero discrete divergence $ \divh{\theB} $ to the Hermite auxiliary variables 
$ \theV_{i,j} $ and $ \theW_{i,j} $, generating the corrected auxiliary variables $ \hat{\theV}_{i,j} $ and $ \hat{\theW}_{i,j} $.
After this correction, the corrected auxiliary variables $ \hat{\theV}_{i,j} $ and $ \hat{\theW}_{i,j} $ 
yield an exactly zero discrete divergence, i.e.,
$ \cleaneddivB = \hat{\theV}^{5}_{i,j} + \hat{\theW}_{i,j}^{6} = 0 $.

We now describe the divergence-free correction procedure in detail.
Given $ \theV_{i,j} $ and $ \theW_{i,j} $ evolved by the HWENO scheme, with $ \divh{\theB} = \theV_{i,j}^{5} + \theW_{i,j}^{6}=\mathcal{O}(h^{k-1}) \neq 0 $, 
the corrected auxiliary variables $ \hat{\theV}_{i,j} $ and $ \hat{\theW}_{i,j} $ are obtained by
\begin{equation*}
		\hat{\theV}_{i,j}^{5} = \theV_{i,j}^{5} - \frac{1}{2} \divh{\theB}, \quad \hat{\theW}_{i,j}^{6} = \theW_{i,j}^{6} - \frac{1}{2} \divh{\theB}.
\end{equation*}
This divergence-free correction has several advantages.
First, since the discrete divergence $\divh{\theB} = \mathcal{O}(h^{k-1})$, the divergence-free correction is a high-order perturbation to the partial derivatives of the magnetic field, i.e., 
\begin{equation*}
    \hat{\theV}_{i,j}^{5} =\theV_{i,j}^{5}+\mathcal{O}(h^{k-1}), \quad
    \hat{\theW}_{i,j}^{6} =\theW_{i,j}^{6}+\mathcal{O}(h^{k-1}).
\end{equation*}
Second, the corrected auxiliary variables $ \hat{\theV}_{i,j} $ and $ \hat{\theW}_{i,j} $ are 
divergence-free, since
\begin{equation*}
	\cleaneddivB = 	\hat{\theV}_{i,j}^{5} + \hat{\theW}_{i,j}^{6} = (\theV_{i,j}^{5} - \frac{1}{2} \divh{\theB}) + (\theW_{i,j}^{6} - \frac{1}{2} \divh{\theB}) = \theV_{i,j}^{5} + \theW_{i,j}^{6} - \divh{\theB} = 0.
\end{equation*}
Third, this correction incurs negligible additional computational cost, since it is a linear and pointwise operation acting on the magnetic field derivatives.
The corrected auxiliary variables $ \hat{\theV}^5_{i,j} $ and $ \hat{\theW}^6_{i,j} $, which are essentially the original values augmented by a high-order perturbation, are then used for spatial reconstruction in the HWENO scheme at the new time level.  And the numerical fluxes in \cref{eq:semi-discrete},
especially the fluxes for the physical variables, 
i.e., $ \hat{\theF}_{i+\half, j} $ and $ \hat{\theG}_{i, j+\half} $, are reconstructed from 
divergence-free data. Therefore the reconstructed numerical fluxes
are more consistent with the physical constraint,
which helps enhance the stability of the numerical solution, 
suppress nonphysical oscillations 
and avoid nonphysical structures.
For clarity, we summarize the divergence-free correction procedure in \cref{alg:DFC}.
\begin{algorithm}[htbp]
	\caption{Divergence-free correction (DFC)}
	\label{alg:DFC}
	\begin{algorithmic}[1]
		\Require non-divergence-free data $ \theV $, $ \theW $
		\Ensure divergence-free data $ \hat{\theV} $, $ \hat{\theW} $
		\For{$i = 1, \ldots, N_x$}
			\For{$j = 1, \ldots, N_y$}
				\State $\divh{\theB} \gets \theV_{i,j}^{5} + \theW_{i,j}^{6}$
				\State $ \hat{\theV}_{i,j}^{5} \gets \theV_{i,j}^{5} - \frac{1}{2} \divh{\theB} $
				\State $ \hat{\theW}_{i,j}^{6} \gets \theW_{i,j}^{6} - \frac{1}{2} \divh{\theB} $
			\EndFor
		\EndFor
	\end{algorithmic}
\end{algorithm}
This approach shares a similar goal with WB schemes for shallow water equations. WB schemes preserve the steady-state balance condition, which can be achieved either by acting a priori on the coefficient factor of the source term \cite{Jin2001SteadystateCapturing, Jin2016WellBalancedStochastic} or by modifying the numerical flux \cite{Xing2006HighOrder}. In contrast, the proposed method for MHD simply adds a single linear divergence-free correction step to the magnetic field derivatives after the evolution stage of the standard HWENO scheme, with the aim of preserving the divergence-free constraint.

\subsection{The temporal discretization}
\label{subsec:time-discretization}
The fully discrete scheme is the third-order strong stability
preserving Runge--Kutta  (SSP-RK3) method,
given in \cref{eq:ssp-rk3}, where the divergence-free correction 
is applied at each Runge--Kutta stage.
Starting from $(\theU^n,\theV^n,\theW^n)$, the SSP-RK3 update is
\begin{equation}
\label{eq:ssp-rk3}
\left\{
\begin{aligned}
	&\left\{
	\begin{alignedat}{2}
		&(\theU^{(1)},\theV,\theW)
		&&=(\theU^n,\tilde{\theV}^n,\tilde{\theW}^n)+\dt\,\mathcal{L}(\theU^n,\theV^n,\theW^n),\\
		&(\theV^{(1)},\theW^{(1)})
		&&=\mathrm{DFC}(\theV,\theW),
	\end{alignedat}
	\right.\\
	&\left\{
	\begin{alignedat}{2}
		&(\theU^{(2)},\theV,\theW)
		&&=\frac{3}{4}(\theU^n,\tilde{\theV}^n,\tilde{\theW}^n)\\
		&
		&&+\frac{1}{4}\left[
		(\theU^{(1)},\tilde{\theV}^{(1)},\tilde{\theW}^{(1)})
		+\dt\,\mathcal{L}(\theU^{(1)},\theV^{(1)},\theW^{(1)})
		\right],\\
		&(\theV^{(2)},\theW^{(2)})
		&&=\mathrm{DFC}(\theV,\theW),
	\end{alignedat}
	\right.\\
	&\left\{
	\begin{alignedat}{2}
		&(\theU^{n+1},\theV,\theW)
		&&=\frac{1}{3}(\theU^n,\tilde{\theV}^n,\tilde{\theW}^n)\\
		&
		&&+\frac{2}{3}\left[
		(\theU^{(2)},\tilde{\theV}^{(2)},\tilde{\theW}^{(2)})
		+\dt\,\mathcal{L}(\theU^{(2)},\theV^{(2)},\theW^{(2)})
		\right],\\
		&(\theV^{n+1},\theW^{n+1})
		&&=\mathrm{DFC}(\theV,\theW).
	\end{alignedat}
	\right.
\end{aligned}
\right.
\end{equation}
Here, for each stage $ \# \in \{n, (1), (2)\} $, $ \mathcal{L}(\theU^{\#}, \theV^{\#}, \theW^{\#}) $ denotes the right-hand side of the semi-discrete scheme \cref{eq:semi-discrete},
$ (\tilde{\theV}^{\#}, \tilde{\theW}^{\#})  = \mathrm{LIMITER} (\theU^{\#}, \theV^{\#}, \theW^{\#}) $,
and $\mathrm{DFC}$ denotes the divergence-free correction in
\cref{alg:DFC}. 
$ \mathrm{LIMITER} $ denotes a HWENO limiter for derivatives, which makes the HWENO scheme robust; one can refer to \cite{Zhao2020ModifiedFifth} and Appendix~\ref{sec:hweno-detail} for more details.
In summary, the proposed divergence-free HWENO scheme achieves high-order accuracy by evolving both the conserved variables and their spatial derivatives,  while the divergence-free correction is applied at every stage of the SSP-RK3 time integration to maintain a divergence-free magnetic field.

\section{Numerical examples}
\label{sec:examples}
In this section, we apply the  base HWENO scheme and the divergence-free HWENO scheme, termed \hwenobase{} and \hwenodf{}, respectively, to a wide range of benchmark examples for the ideal MHD equations to demonstrate their accuracy, robustness, shock-capturing capability and effectiveness. 
The CFL number is set to $0.6$ for all simulations.
The total CPU time is reported in \cref{tab:elapsed-time}, from which it can be seen that the additional cost introduced by the linear divergence-free correction technique in the DFHWENO scheme is almost negligible compared to the HWENO scheme. The reported CPU times were obtained from serial computations performed on a single physical core of an Intel Xeon Gold 6130 processor running at 2.10 GHz. The numerical code was written in Fortran 90 and compiled using the Intel Fortran Compiler (ifort), version 19.0.1.144, 
with the \texttt{-O3} optimization option enabled.

\crefname{example}{Ex.}{Exs.}
\begin{table}[htbp]
\centering
\caption{CPU time comparison between the DFHWENO and HWENO schemes.}
\label{tab:elapsed-time}
\begin{tabularx}{\linewidth}{@{}Xrrr@{}}
\toprule
\multirow{2}{*}{Example}
& \multicolumn{2}{c}{CPU time (s)}
& \multirow{2}{*}{Additional cost (\%)} \\
\cmidrule(lr){2-3}
& DFHWENO & HWENO & \\
\midrule
\cref{ex:orszag-tang} Orszag--Tang vortex      & 967.70   & failed      & --- \\
\cref{ex:rotor} Rotor                 & 2681.22  & 2557.04  & 4.86 \\
\cref{ex:extreme-blast} Extreme blast & 332.55   & 339.53   & $-2.06$ \\
\cref{ex:cloud-shock} Cloud--shock interaction      & 16030.62 & 16512.07 & $ -2.92 $ \\
\cref{ex:kelvin-helmholtz} Kelvin--Helmholtz instability & 16042.90 & 15380.07 & 4.31 \\
\cref{ex:field-loop} Field loop advection        & 4599.04  & 4345.67  & 5.83 \\
\cref{ex:brio-wu} 2D rotated shock tube  & 48.59 & 47.24 & 2.84 \\
\bottomrule
\end{tabularx}
\end{table}

\begin{example}[The smooth Alfv\'en wave]
\label{ex:alfven-wave}
To examine the high-order accuracy of the proposed schemes, we simulate the 
nonlinear circularly polarized Alfv\'en wave \cite{Toth2000B0Constraint}. 
The computational domain is 
$[0, 1/\cos \alpha]\times[0, 1/\sin \alpha]$, 
where $\alpha$ is the angle between the wave propagation direction and 
the $x$-axis, taken as $ \alpha = \pi/4 $. 
The initial density, pressure, and adiabatic index are given by
$
    \rho = 1, p = 1, \gamma = 5/3.
	$
The initial vector primitive variables in the directions parallel and 
perpendicular to the wave propagation direction are prescribed as
\begin{equation*}
	u_{\parallel} = 0,
	\quad u_{\perp} = 0.1\sin(2\pi\xi),
	\quad u^{3} = 0.1\cos(2\pi\xi),\\
	\quad B_{\parallel} = 1,
	\quad B_{\perp} = u_{\perp},
	\quad B^{3} = u^{3}.
\end{equation*}
Here, $\xi = x\cos \alpha + y\sin \alpha$ is the coordinate along the 
wave propagation direction. 
Using the rotation between the parallel-perpendicular coordinates and the 
Cartesian coordinates
\begin{equation*}
	\begin{aligned}
		&u^{1} = u_{\parallel}\cos \alpha - u_{\perp}\sin \alpha, \quad u^{2} = u_{\parallel}\sin \alpha + u_{\perp}\cos \alpha,\\
		&B^{1} = B_{\parallel}\cos \alpha - B_{\perp}\sin \alpha, \quad B^{2} = B_{\parallel}\sin \alpha + B_{\perp}\cos \alpha,
	\end{aligned}
\end{equation*}
we obtain the initial primitive variables in the $x$- and $y$-directions.
The Alfv\'en wave propagates along the direction 
$\vec{n} = (\cos \alpha, \sin \alpha)$ with speed 
$v_{A} = B_{\parallel}/\sqrt{\rho}$. 
The exact solution at time $t$ is obtained by replacing $\xi$ with 
$\xi - v_{A}t$ in the initial condition. 
We simulate this problem until $t = 2$ and demonstrate the $L^{1}$ errors of 
several representative primitive variables in \cref{tab:alfven}. 
The results show that the proposed schemes achieve fifth-order accuracy with very close numerical errors, which indicates that the divergence-free correction of magnetic field derivatives in the DFHWENO scheme does not compromise the accuracy of the HWENO scheme. To further investigate this, we present the divergence error of the magnetic field and its associated derivatives in \cref{tab:alfven-div-accuracy-hweno}. 
 From the table, we can first observe that the divergence error of the \hwenobase{} scheme is a high-order perturbation, while that of the \hwenodf{} scheme is near machine precision, confirming that the \hwenodf{} scheme is divergence-free. We can also see that the divergence error is much smaller than the numerical error of the solution, which explains why the numerical errors of \hwenobase{} and \hwenodf{} are very close.

\begin{table}[htbp]
\centering
\setlength{\tabcolsep}{3pt}
\caption{\cref{ex:alfven-wave}. The smooth Alfv\'en wave, $L^1$ errors and convergence orders of  \hwenodf{}  and \hwenobase{} schemes at $t = 2$. $ N_{x} = N_{y} $.}
\label{tab:alfven}
\begin{tabularx}{\linewidth}{@{}Xllllllll@{}}
\toprule
$N_x$ & \multicolumn{2}{c}{$u_{1}$} & \multicolumn{2}{c}{$u_{2}$} & \multicolumn{2}{c}{$B_{1}$} & \multicolumn{2}{c}{$B_{2}$} \\
\cmidrule(lr){2-3}
\cmidrule(lr){4-5}
\cmidrule(lr){6-7}
\cmidrule(lr){8-9}
& $L^1$ error & Order & $L^1$ error & Order & $L^1$ error & Order & $L^1$ error & Order \\
\midrule
\multicolumn{9}{c}{\hwenodf{}} \\
  40 & 3.6324e-06 & -- & 3.6324e-06 & -- & 3.6291e-06 & -- & 3.6291e-06 & -- \\
  80 & 1.1087e-07 & 5.0340 & 1.1087e-07 & 5.0340 & 1.1053e-07 & 5.0371 & 1.1053e-07 & 5.0371 \\
  120 & 1.4309e-08 & 5.0497 & 1.4309e-08 & 5.0497 & 1.4265e-08 & 5.0497 & 1.4265e-08 & 5.0497 \\
  160 & 3.3206e-09 & 5.0776 & 3.3206e-09 & 5.0776 & 3.3078e-09 & 5.0803 & 3.3078e-09 & 5.0803 \\
  200 & 1.0616e-09 & 5.1102 & 1.0616e-09 & 5.1102 & 1.0567e-09 & 5.1140 & 1.0567e-09 & 5.1140 \\
  240 & 4.1555e-10 & 5.1446 & 4.1555e-10 & 5.1446 & 4.1341e-10 & 5.1472 & 4.1341e-10 & 5.1472 \\
\midrule
\multicolumn{9}{c}{\hwenobase{}} \\
  40 & 3.6391e-06 & -- & 3.6391e-06 & -- & 3.6219e-06 & -- & 3.6219e-06 & -- \\
  80 & 1.1095e-07 & 5.0356 & 1.1095e-07 & 5.0356 & 1.1044e-07 & 5.0354 & 1.1044e-07 & 5.0354 \\
  120 & 1.4318e-08 & 5.0499 & 1.4318e-08 & 5.0499 & 1.4255e-08 & 5.0495 & 1.4255e-08 & 5.0495 \\
  160 & 3.3229e-09 & 5.0775 & 3.3229e-09 & 5.0775 & 3.3054e-09 & 5.0804 & 3.3054e-09 & 5.0804 \\
  200 & 1.0624e-09 & 5.1101 & 1.0624e-09 & 5.1101 & 1.0559e-09 & 5.1141 & 1.0559e-09 & 5.1141 \\
  240 & 4.1587e-10 & 5.1444 & 4.1587e-10 & 5.1444 & 4.1307e-10 & 5.1475 & 4.1307e-10 & 5.1475 \\
\bottomrule
\end{tabularx}
\end{table}

\begin{table}[htbp]
\centering
\caption{ \cref{ex:alfven-wave}. The smooth Alfv\'en wave, $L^1$ errors and convergence orders of magnetic field derivatives and divergence for \hwenodf{} and \hwenobase{} at $t = 2$. $ N_x=N_y. $}
\label{tab:alfven-div-accuracy-hweno}
\begin{tabularx}{\linewidth}{@{}Xllllll@{}}
\toprule
$N_x$ & \multicolumn{2}{c}{$\pdv*{B^{1}}{x}$} & \multicolumn{2}{c}{$\pdv*{B^{2}}{y}$} & \multicolumn{2}{c}{$\div\theB$} \\
\cmidrule(lr){2-3}
\cmidrule(lr){4-5}
\cmidrule(lr){6-7}
& $L^1$ error & Order & $L^1$ error & Order & $L^1$ error & Order \\
\midrule
\multicolumn{7}{c}{\hwenodf{}} \\
  40 & 4.4934e-05 & -- & 4.4934e-05 & -- & 6.6050e-18 & -- \\
  80 & 1.6105e-06 & 4.8022 & 1.6105e-06 & 4.8022 & 7.1102e-18 & -0.1063 \\
  120 & 2.3019e-07 & 4.7979 & 2.3019e-07 & 4.7979 & 7.9939e-18 & -0.2889 \\
  160 & 5.8649e-08 & 4.7530 & 5.8649e-08 & 4.7530 & 8.5824e-18 & -0.2469 \\
  200 & 2.0302e-08 & 4.7541 & 2.0302e-08 & 4.7541 & 8.8787e-18 & -0.1521 \\
  240 & 8.4542e-09 & 4.8050 & 8.4542e-09 & 4.8050 & 7.4534e-18 & 0.9597 \\
\midrule
\multicolumn{7}{c}{\hwenobase{}} \\
  40 & 4.4906e-05 & -- & 4.4906e-05 & -- & 1.0973e-07 & -- \\
  80 & 1.6101e-06 & 4.8017 & 1.6101e-06 & 4.8017 & 1.2642e-09 & 6.4395 \\
  120 & 2.3014e-07 & 4.7978 & 2.3014e-07 & 4.7978 & 1.1154e-10 & 5.9878 \\
  160 & 5.8638e-08 & 4.7529 & 5.8638e-08 & 4.7529 & 1.5800e-11 & 6.7934 \\
  200 & 2.0299e-08 & 4.7540 & 2.0299e-08 & 4.7540 & 2.7613e-12 & 7.8171 \\
  240 & 8.4528e-09 & 4.8050 & 8.4528e-09 & 4.8050 & 1.0251e-12 & 5.4347 \\
\bottomrule
\end{tabularx}
\end{table}
    
\end{example}

\begin{example}[Orszag-Tang vortex]
	\label{ex:orszag-tang}
	The Orszag-Tang vortex \cite{Orszag1979SmallscaleStructure} is a widely used benchmark example for testing the performance of numerical schemes for ideal MHD equations, since complex interaction between several shocks develops as the vortex system evolves. The computational domain is $ [0, 2\pi]\times[0, 2\pi] $ with periodic boundary conditions. The initial condition is given by
	\begin{equation*}
	    \primitives = (\gamma^{2}, -\sin(y), \sin(x), 0, -\sin(y), \sin(2x), 0, \gamma),
	\end{equation*}
with $ \gamma = 5/3. $
The numerical solution is simulated until $ t = 4.0 $ with a uniform mesh of $ 200\times 200 $ cells.
As depicted in \cref{fig:orszag-tang}, the proposed \hwenodf{} scheme can capture the complex vortex structure without noticeable spurious oscillations. At $ t = 3.0 $, the \hwenodf{} scheme produces better results than the \hwenobase{} scheme. Consistent with \cite{Ding2024NewDiscretely}, the HWENO scheme without any divergence-free treatment fails before $t=4.0$. Specifically, the \hwenobase{} scheme produces nonphysical solutions and breaks down near $t\approx 3.06$, while the \hwenodf{} scheme remains stable up to $t=4.0$. This agreement confirms that the proposed divergence-free correction of magnetic field derivatives is both crucial and effective.
\begin{figure}[htbp]
	\centering
	\begin{subfigure}{.23\textwidth}
		\centering
		\includegraphics[width=\linewidth]{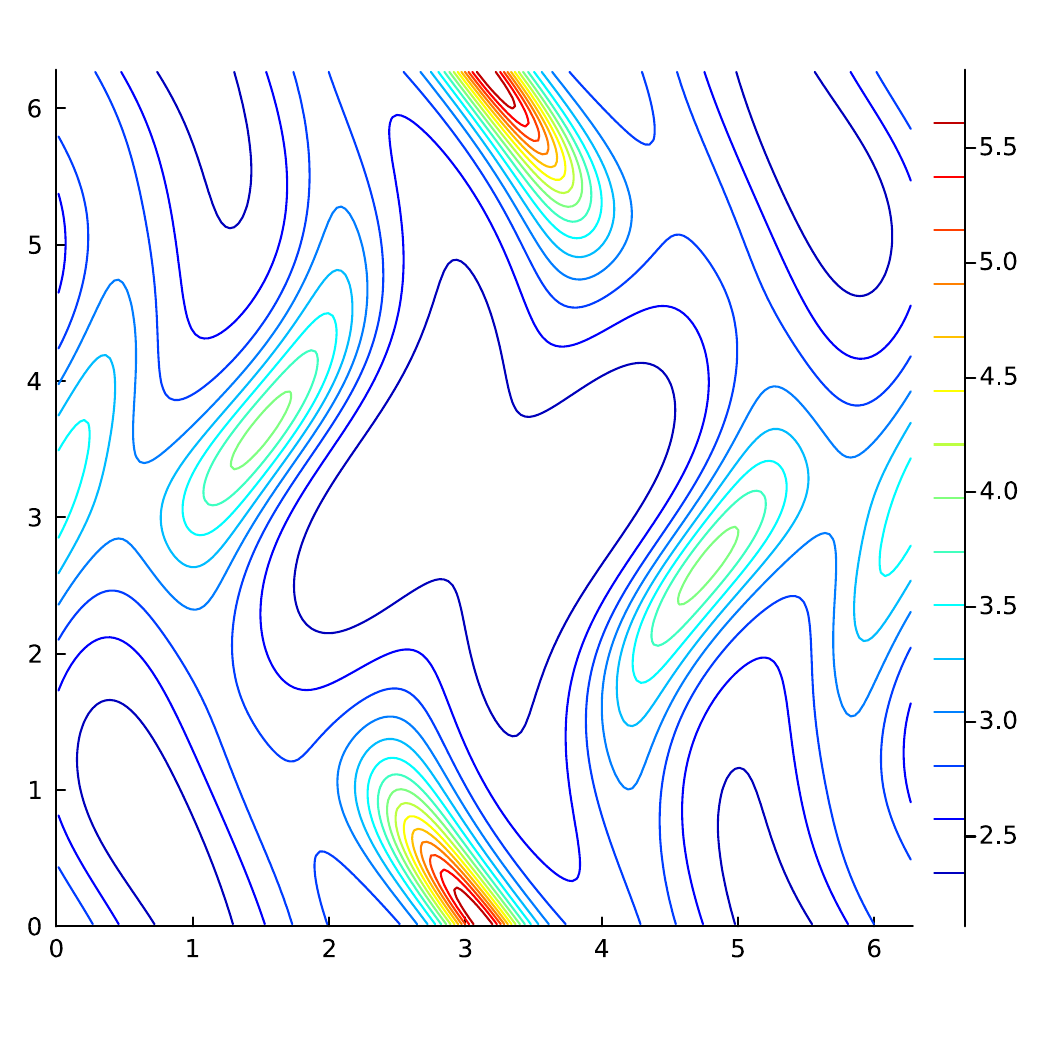}
		\caption{\hwenodf{}, $ t = 0.5 $}
	\end{subfigure}
	\hfill
	\begin{subfigure}{.23\textwidth}
		\centering
		\includegraphics[width=\linewidth]{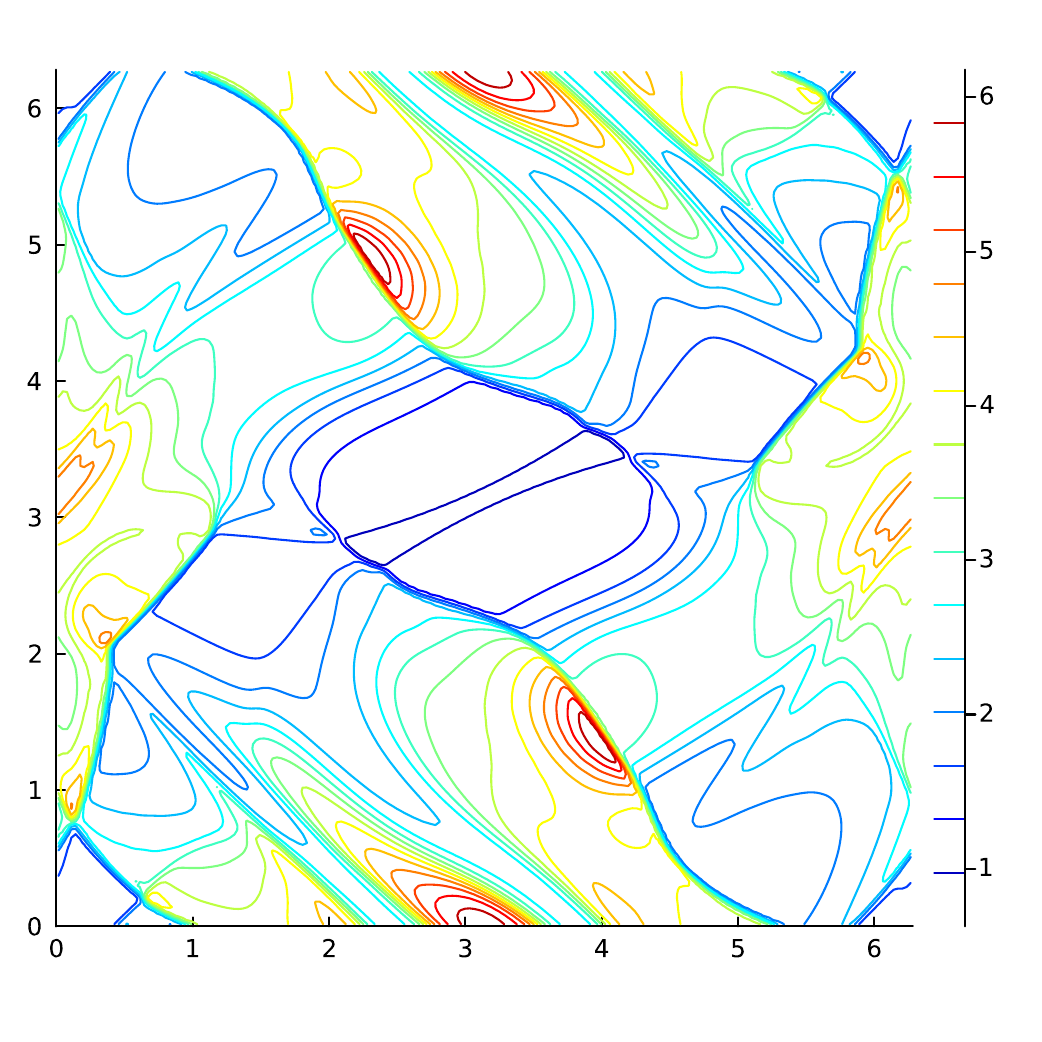}
		\caption{\hwenodf{}, $ t = 2.0 $}
	\end{subfigure}
	\hfill
	\begin{subfigure}{.23\textwidth}
		\centering
		\includegraphics[width=\linewidth]{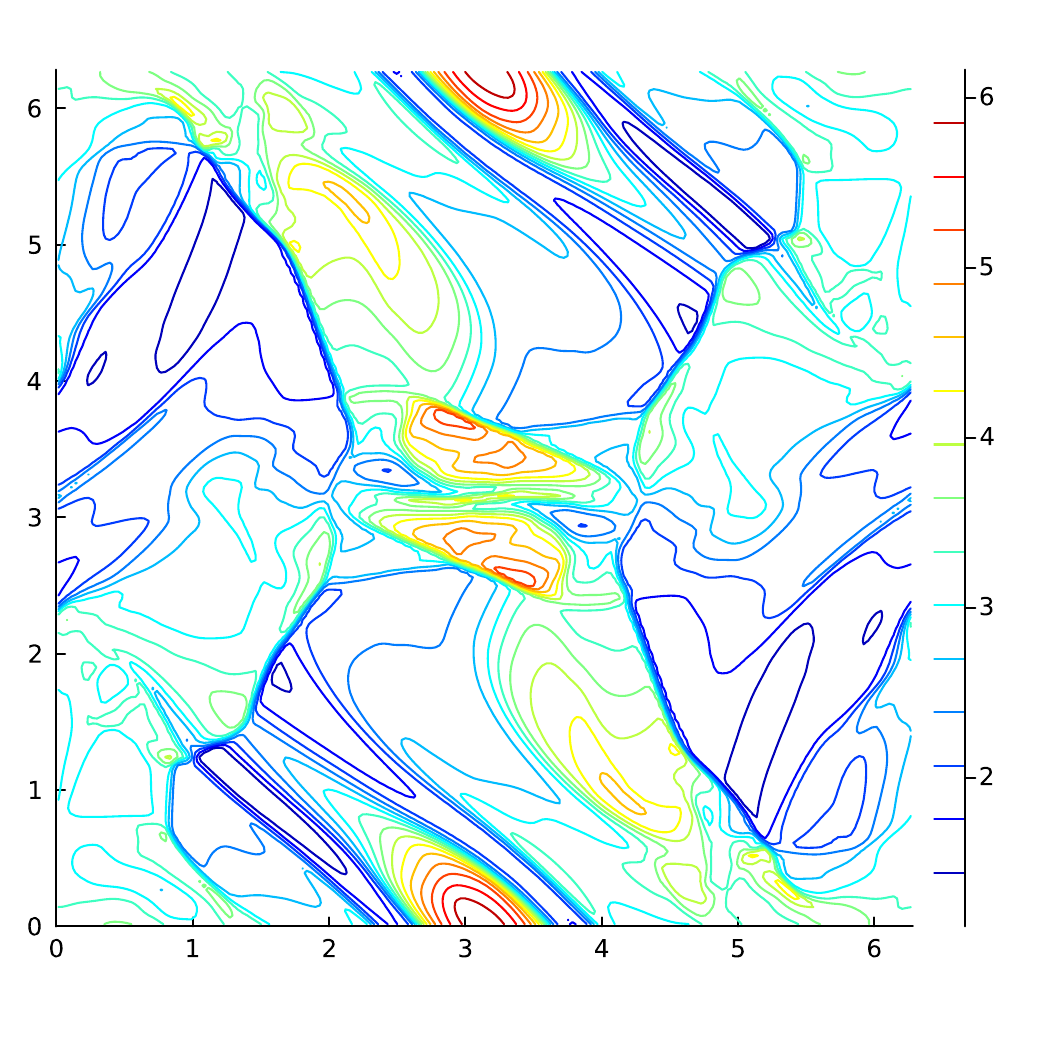}
		\caption{\hwenodf{}, $ t = 3.0 $}
	\end{subfigure}
	\hfill
	\begin{subfigure}{.23\textwidth}
		\centering
		\includegraphics[width=\linewidth]{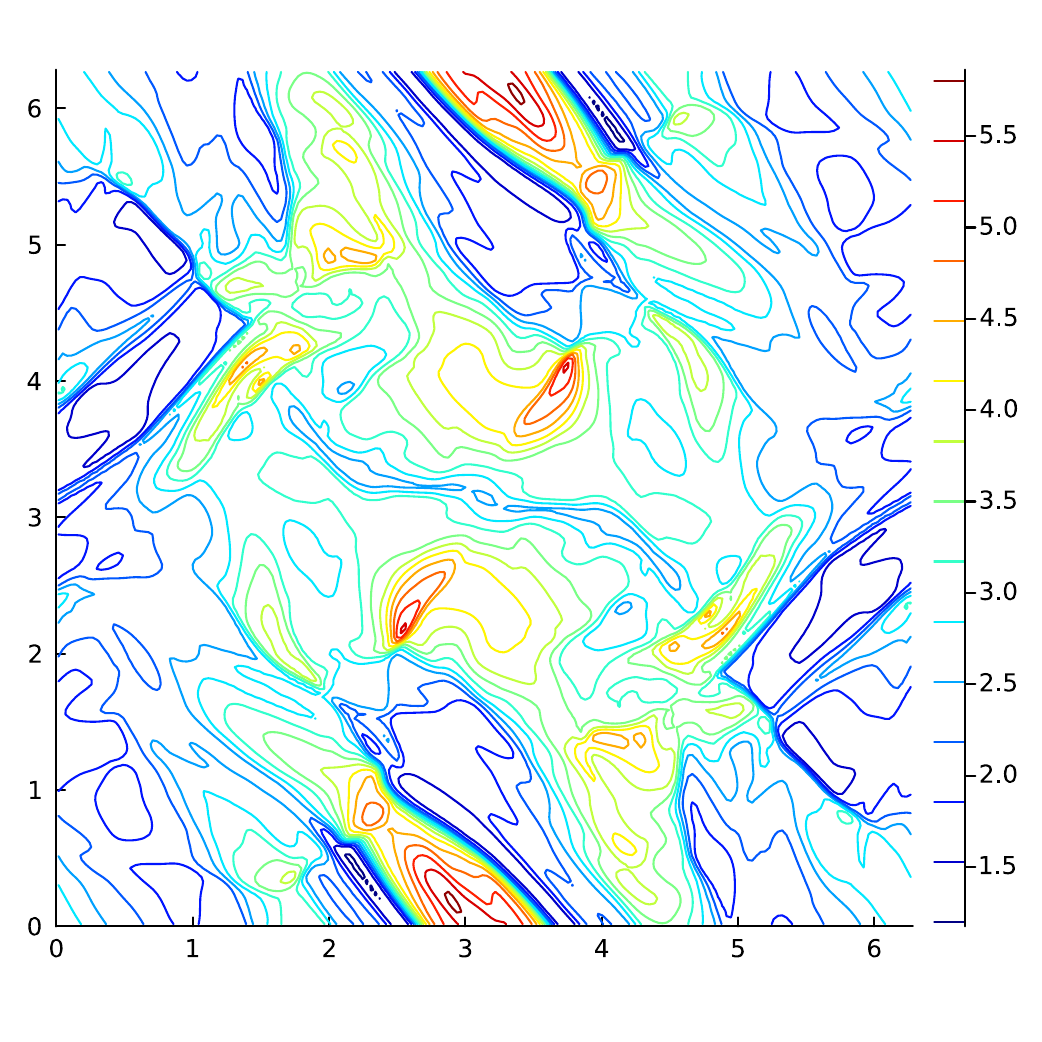}
		\caption{\hwenodf{}, $ t = 4.0 $}
	\end{subfigure}

	\begin{subfigure}{.23\textwidth}
		\centering
		\includegraphics[width=\linewidth]{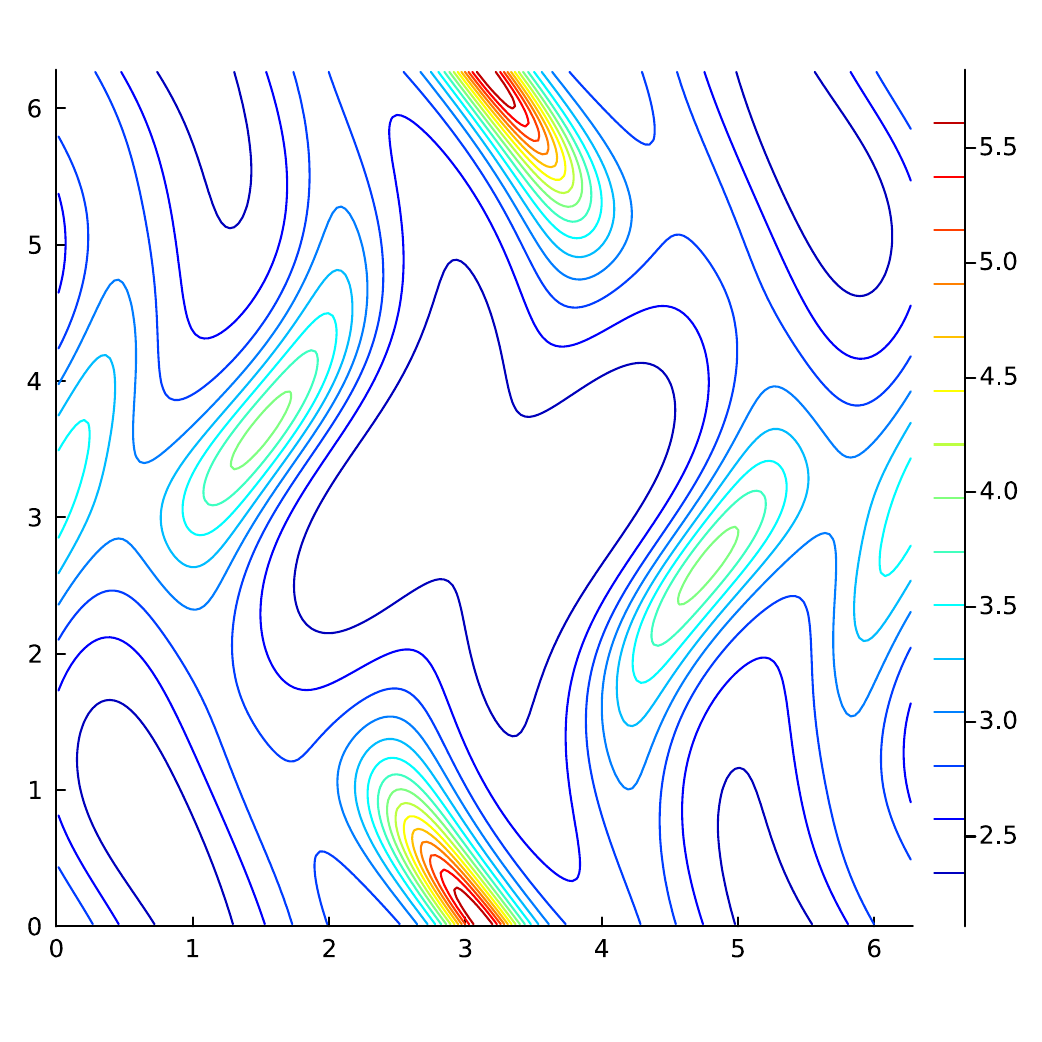}
		\caption{\hwenobase{}, $ t = 0.5 $}
	\end{subfigure}
	\hfill
	\begin{subfigure}{.23\textwidth}
		\centering
		\includegraphics[width=\linewidth]{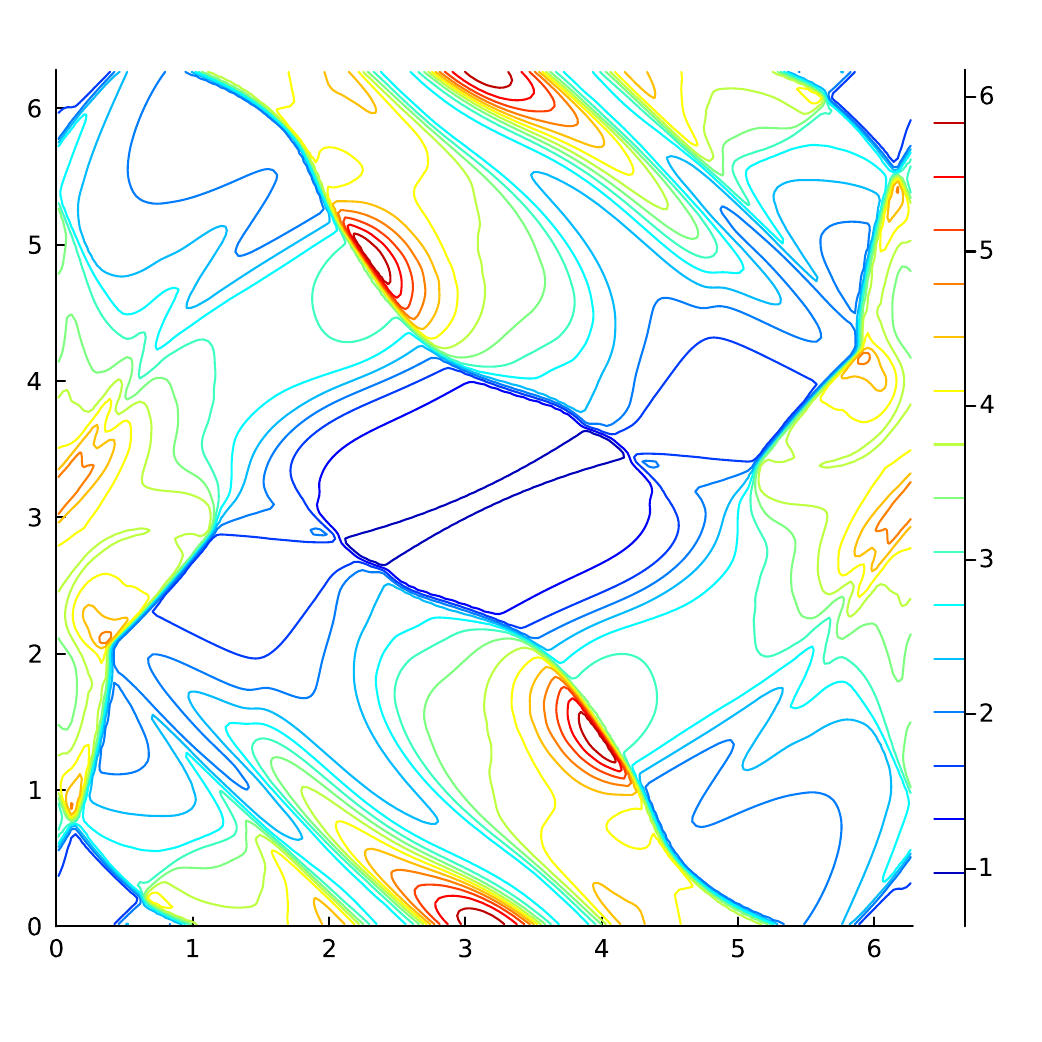}
		\caption{\hwenobase{}, $ t = 2.0 $}
	\end{subfigure}
	\hfill
	\begin{subfigure}{.23\textwidth}
		\centering
		\includegraphics[width=\linewidth]{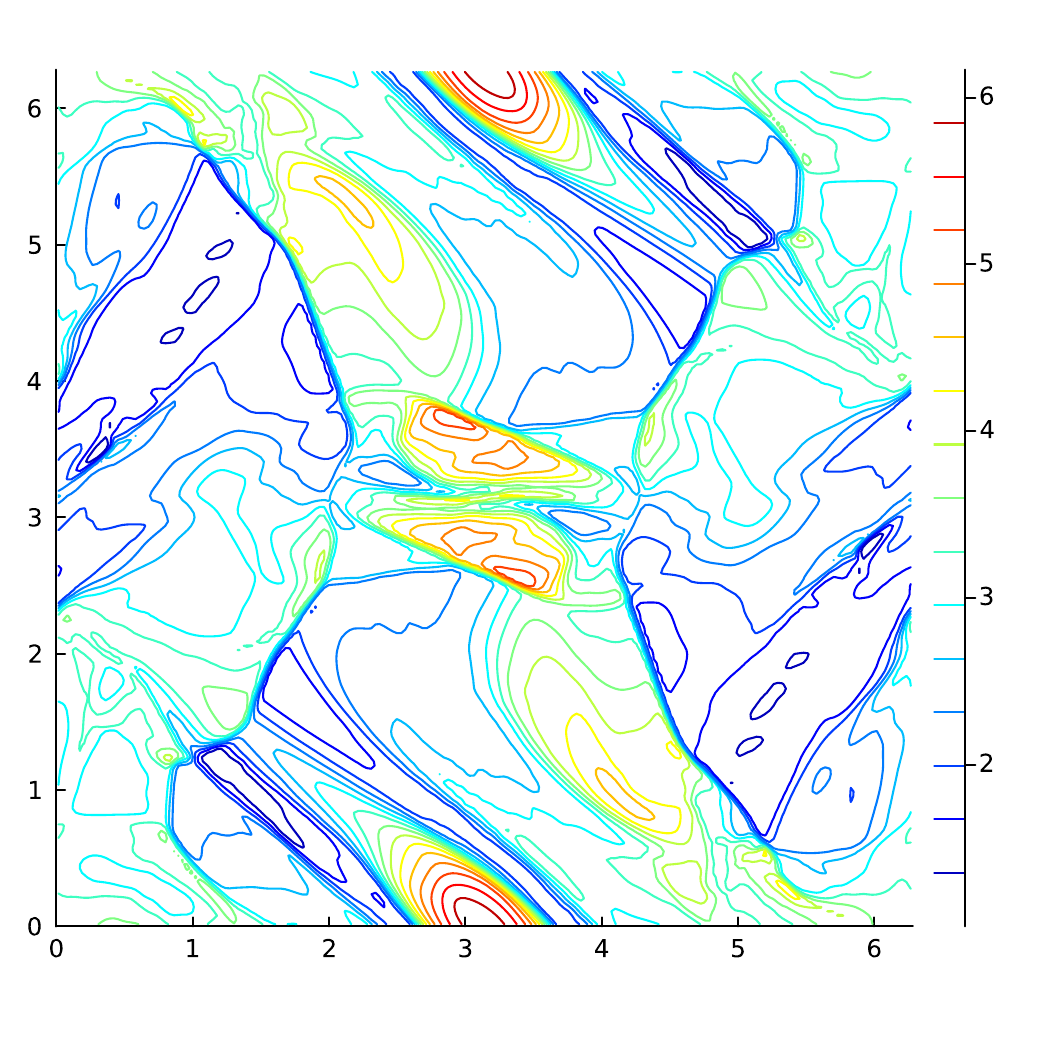}
		\caption{\hwenobase{}, $ t = 3.0 $}
	\end{subfigure}
	\hfill
	\begin{subfigure}{.23\textwidth}
		\centering
		\includegraphics[width=\linewidth]{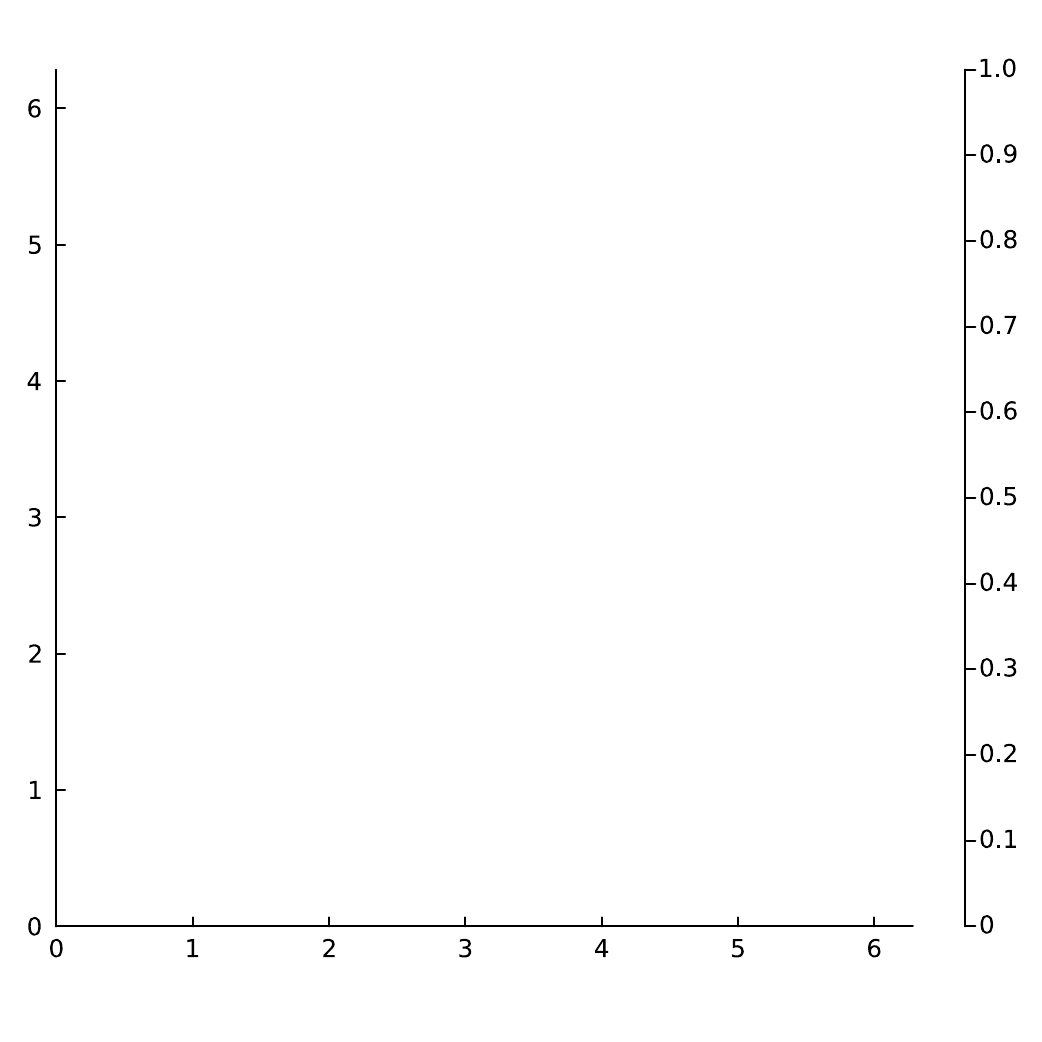}
		\caption{\hwenobase{}, failed.}
	\end{subfigure}
	\caption{\cref{ex:orszag-tang}. Orszag--Tang vortex, simulated by \hwenodf{} (top row) and \hwenobase{} (bottom row) with $200\times 200$ meshes. Density contours at $t=0.5,2.0,3.0,4.0$. 15 equally spaced contour lines are plotted.}
	\label{fig:orszag-tang}
\end{figure}
\end{example}
\begin{example}[Rotor]
	\label{ex:rotor}
	The rotor problem \cite{Toth2000B0Constraint} is computed on the domain $[0,1]\times[0,1]$ with outflow boundary conditions.  A dense rotating disk centered at $(0.5,0.5)$ is embedded in a uniform magnetized medium. 
The initial condition is given by
\begin{equation*}
\primitives =
\left\{
\begin{aligned}
&
(
10,
-\frac{y-0.5}{r_1},
\frac{x-0.5}{r_1},
0,
\frac{2.5}{\sqrt{4\pi}},
0,
0,
0.5
),
&& r<r_1,
\\
&
(
1+9\phi,
-\phi\frac{y-0.5}{r},
\phi\frac{x-0.5}{r},
0,
\frac{2.5}{\sqrt{4\pi}},
0,
0,
0.5
),
&& r_1\leqslant r\leqslant r_2,
\\
&
(
1,
0,
0,
0,
\frac{2.5}{\sqrt{4\pi}},
0,
0,
0.5
),
&& r>r_2,
\end{aligned}
\right.
\end{equation*}
where
$ r=\sqrt{(x-0.5)^2+(y-0.5)^2},
\phi=\frac{r_2-r}{r_2-r_1},
r_1=0.1,
r_2=0.115 $,
and $ \gamma =  5/3.$
We simulate this problem until $ t = 0.295 $ with a uniform mesh of $ 400\times 400 $ cells.
In \cref{fig:rotor}, the presented results show that the \hwenodf{} scheme simulates this problem well and the \hwenobase{} scheme can introduce noticeable oscillations.
In \cref{fig:rotor-mach-zoomed}, the zoomed view of the central area of contour plots of Mach number shows that significant distortions of the velocity field are observed in the \hwenobase{} scheme, while the \hwenodf{} scheme captures the velocity field without visible distortion. 
This example clearly demonstrates 
the importance of the divergence-free constraint for simulating MHD problems,
and the effectiveness of the proposed divergence-free correction for handling this constraint. 
\begin{figure}[htbp]
	\centering
	\begin{subfigure}{.3\textwidth}
		\centering
		\includegraphics[width=\linewidth]{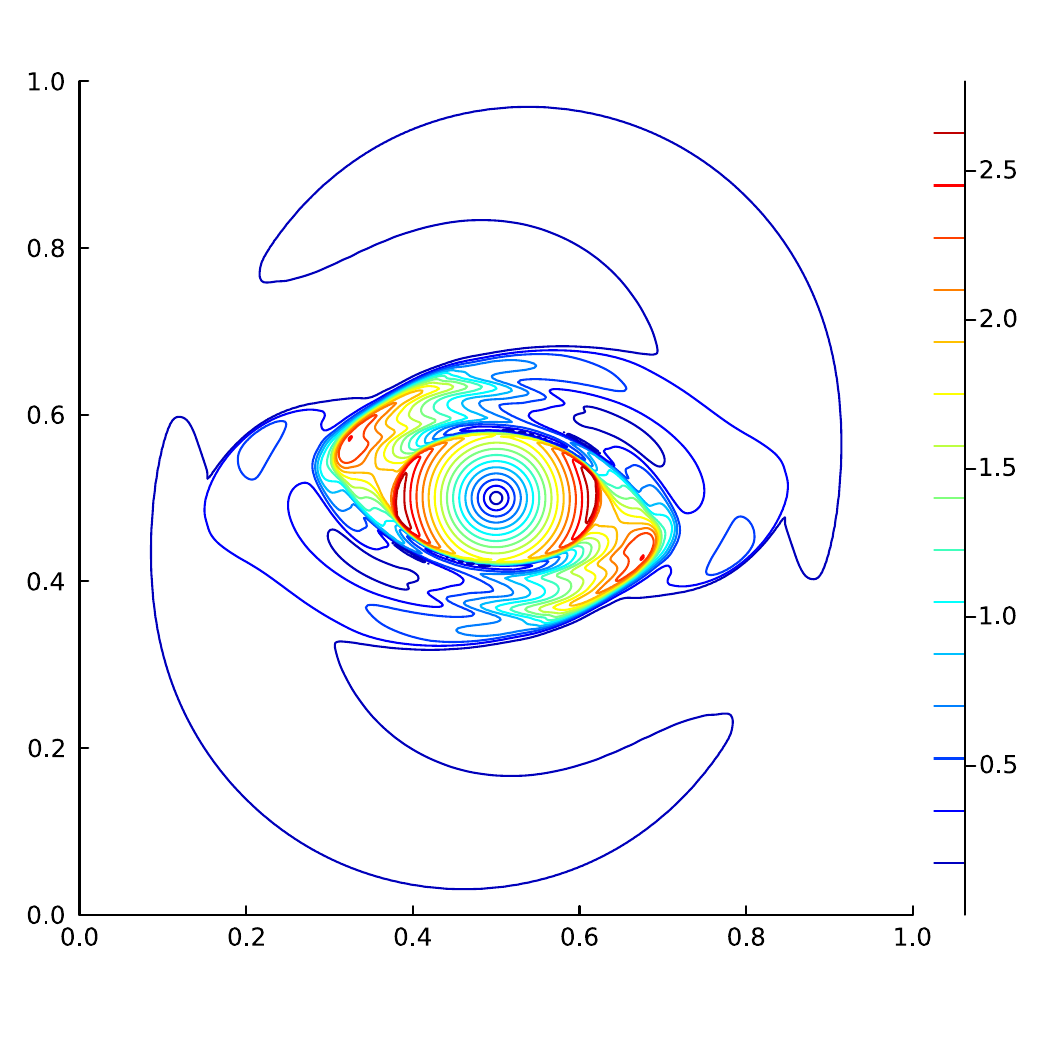}
		\caption{\hwenodf{}, $ \abs{\vec{u}}/c $}
	\end{subfigure}
	\hfill
	\begin{subfigure}{.3\textwidth}
		\centering
		\includegraphics[width=\linewidth]{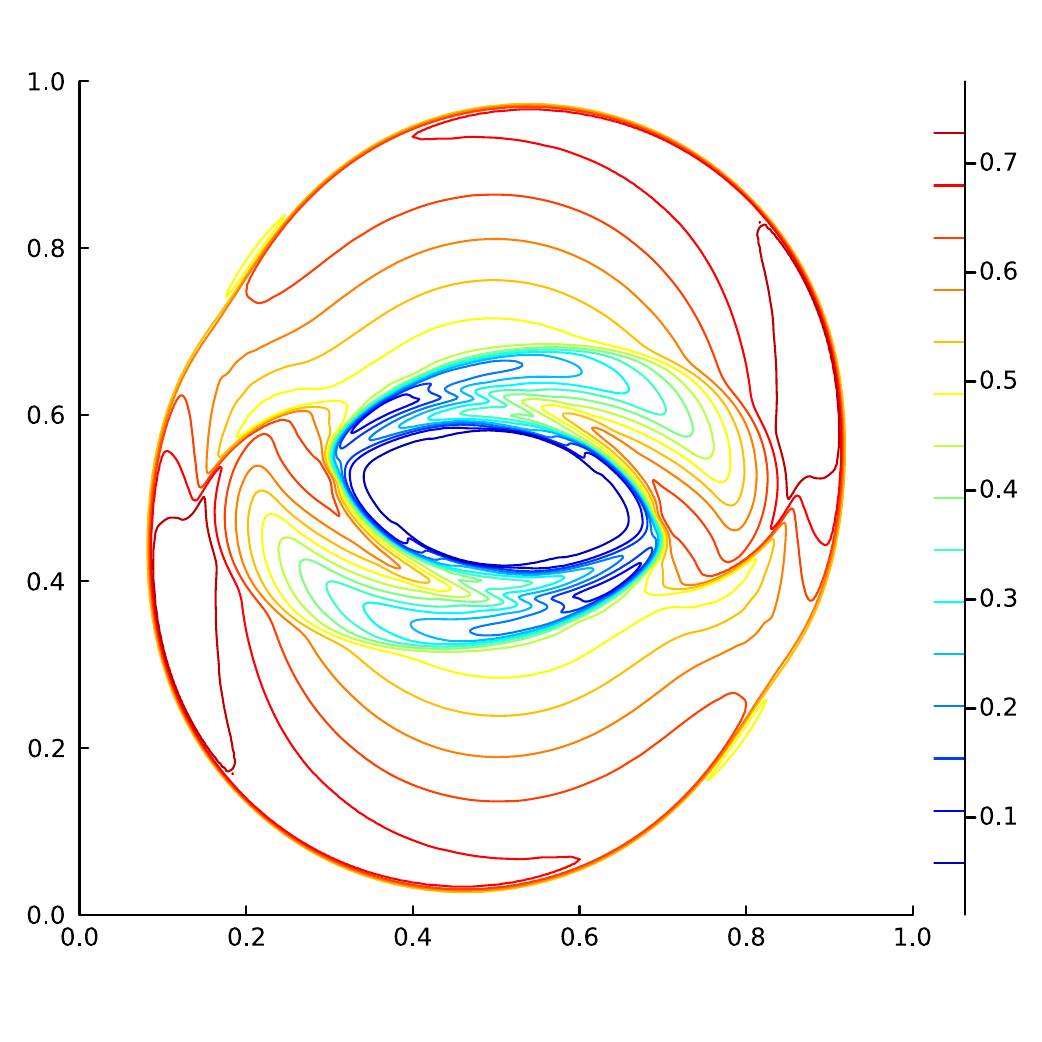}
		\caption{\hwenodf{}, $ p $}
	\end{subfigure}
	\hfill
	\begin{subfigure}{.3\textwidth}
		\centering
		\includegraphics[width=\linewidth]{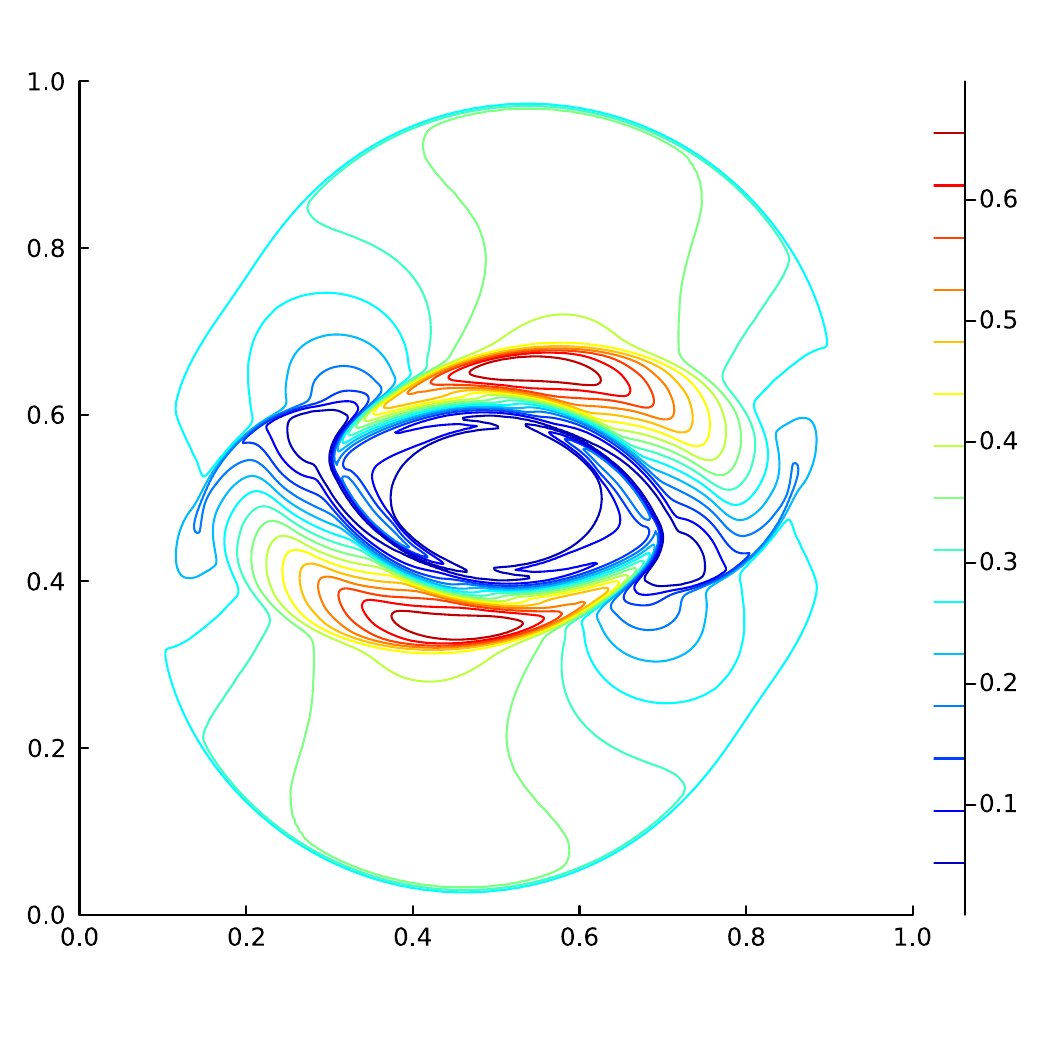}
		\caption{\hwenodf{}, $ \frac{1}{2}\norm{\theB}^{2} $}
	\end{subfigure}

	\begin{subfigure}{.3\textwidth}
		\centering
		\includegraphics[width=\linewidth]{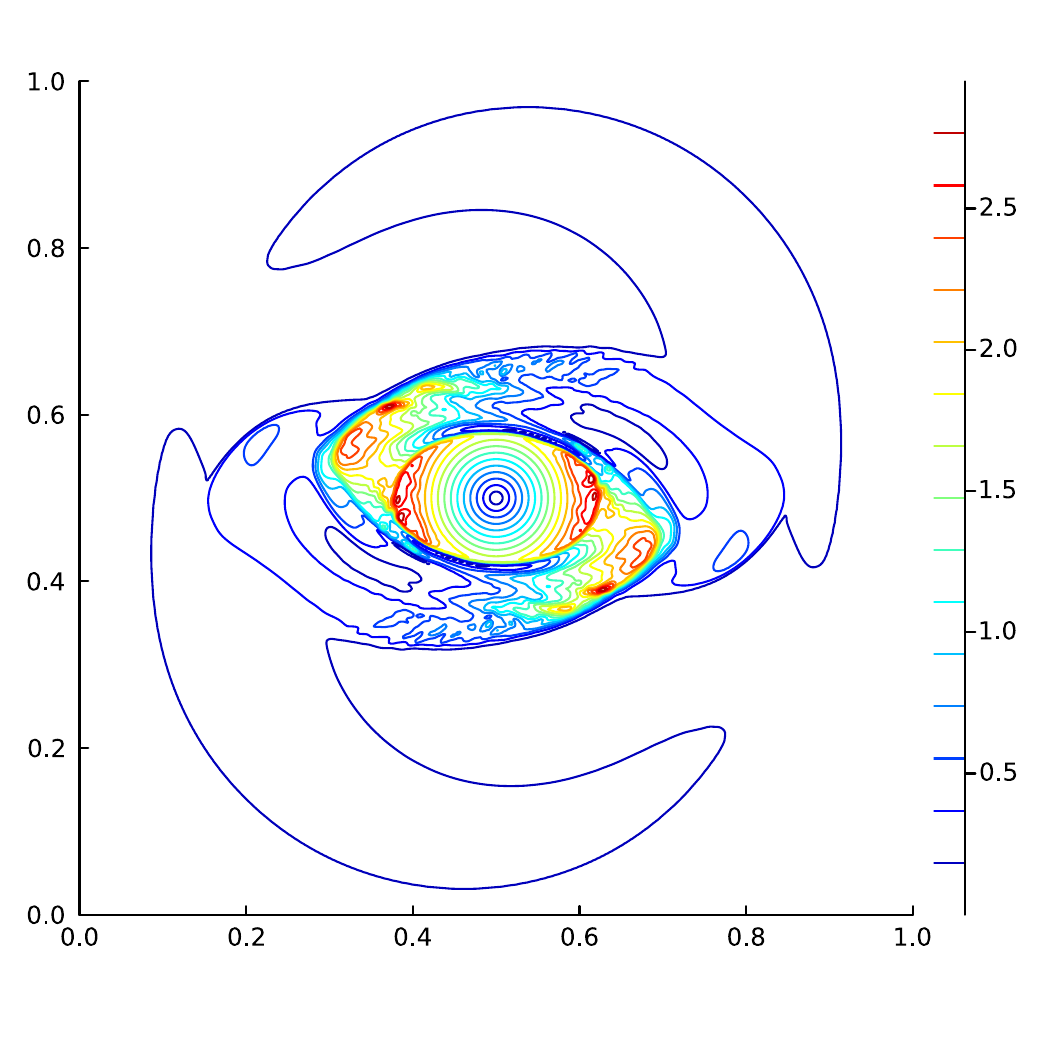}
		\caption{\hwenobase{}, $ \abs{\vec{u}}/c $}
	\end{subfigure}
	\hfill
	\begin{subfigure}{.3\textwidth}
		\centering
		\includegraphics[width=\linewidth]{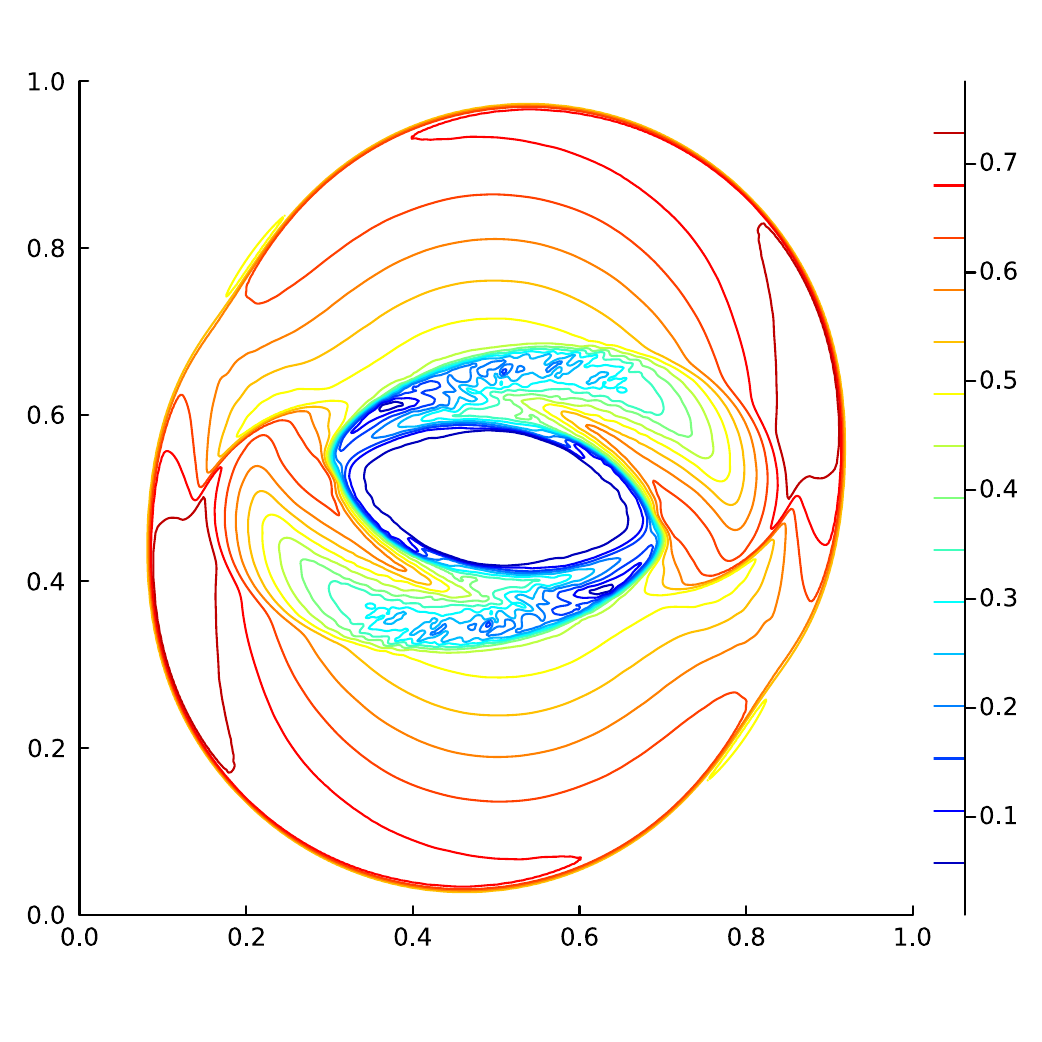}
		\caption{\hwenobase{}, $ p $}
	\end{subfigure}
	\hfill
	\begin{subfigure}{.3\textwidth}
		\centering
		\includegraphics[width=\linewidth]{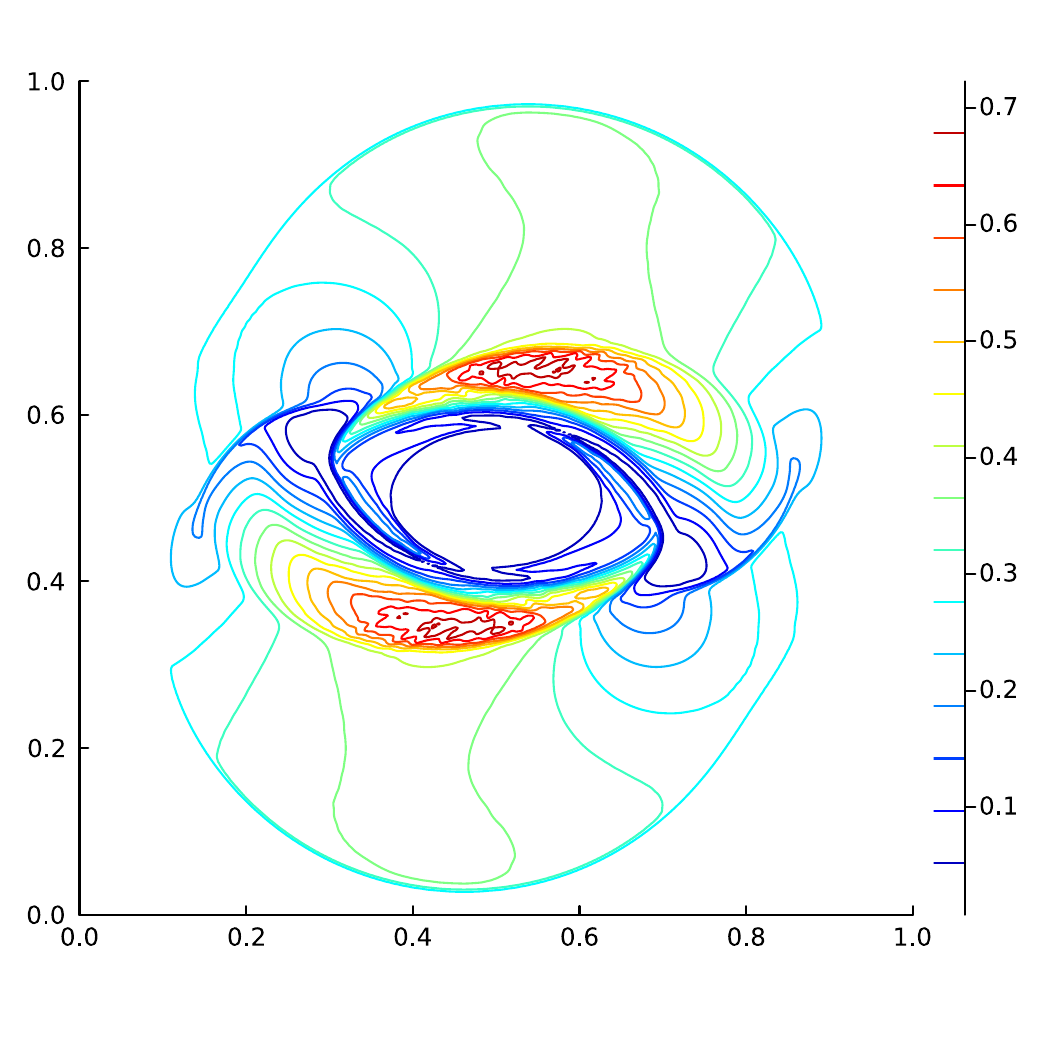}
		\caption{\hwenobase{}, $ \frac{1}{2}\norm{\theB}^{2} $}
	\end{subfigure}
	\caption{\cref{ex:rotor}. Rotor, simulated by \hwenodf{} (top row) and \hwenobase{} (bottom row) with $400\times 400$ meshes. Contour plots of Mach number, pressure, and magnetic pressure at $t=0.295$. 15 equally spaced contour lines are plotted.}
	\label{fig:rotor}
\end{figure}

\begin{figure}[htbp]
	\centering
	\begin{subfigure}{.4\textwidth}
    \centering
		\includegraphics[width=\linewidth]{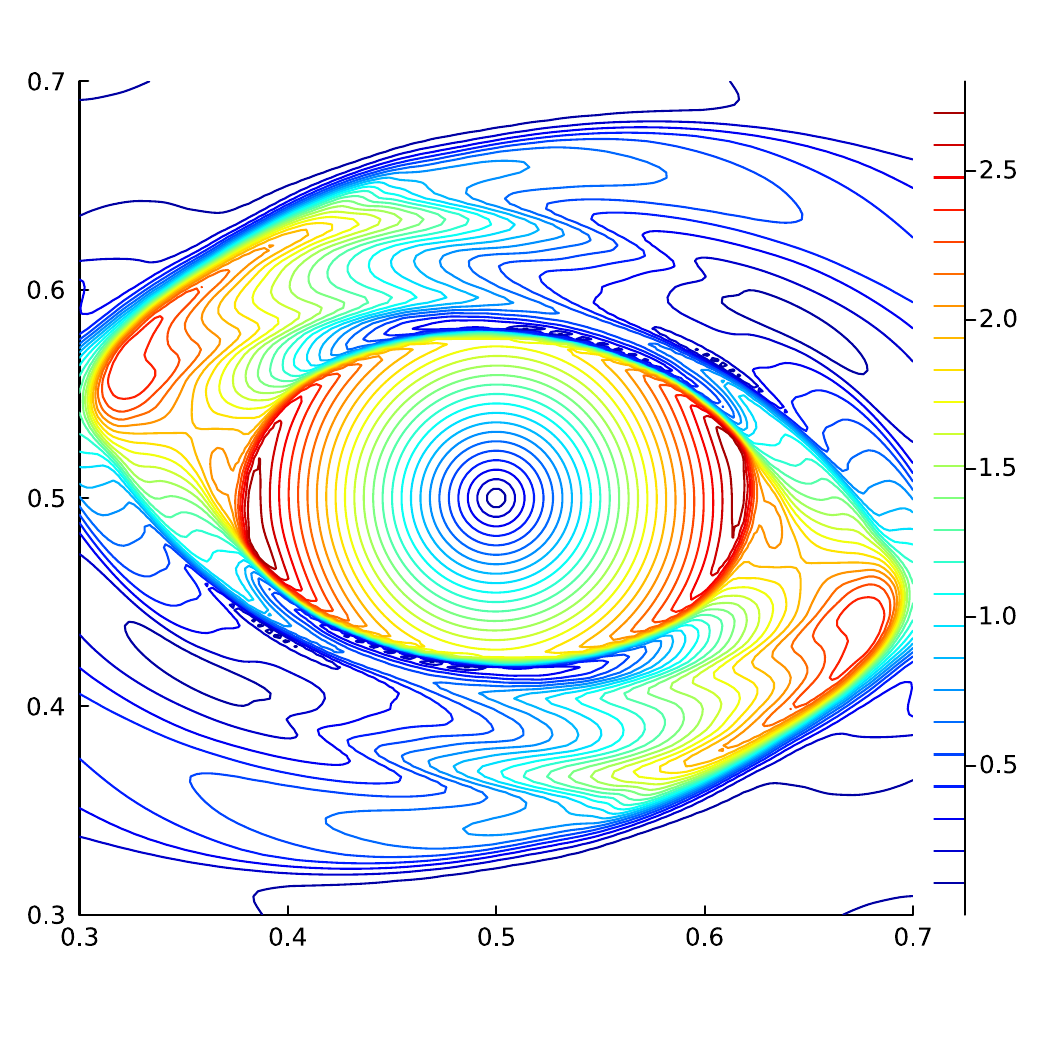}
		\caption{\hwenodf{}}
	\end{subfigure}
	\hfil
	\begin{subfigure}{.4\textwidth}
    \centering
		\includegraphics[width=\linewidth]{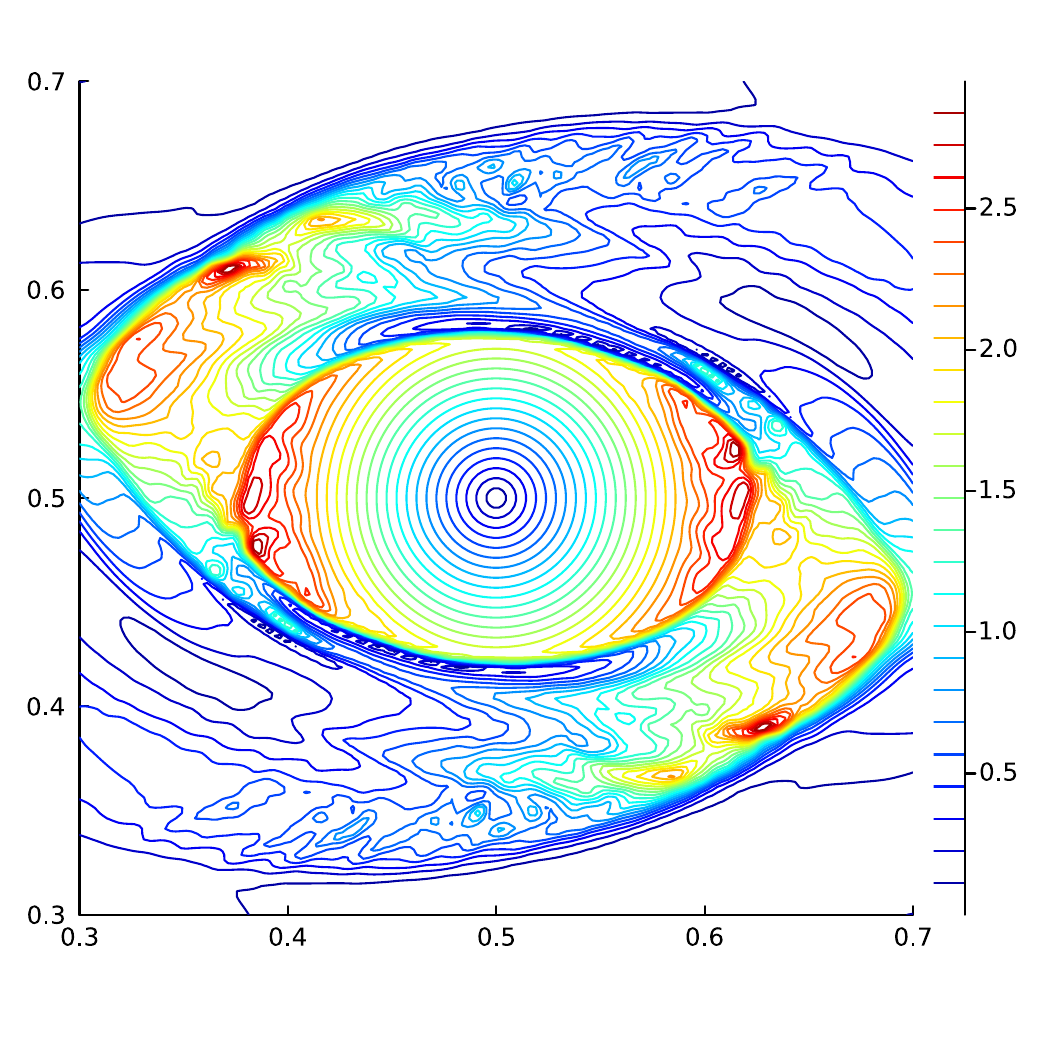}
		\caption{\hwenobase{}}
	\end{subfigure}
	\caption{\cref{ex:rotor}. Rotor, simulated by \hwenodf{} (left) and \hwenobase{} (right) with $400\times 400$ meshes. Zoomed view of Mach number contours in the central area at $t=0.295$. 30 equally spaced contour lines are plotted.}
	\label{fig:rotor-mach-zoomed}
\end{figure}
\end{example}

\begin{example}[Extreme blast]
	\label{ex:extreme-blast}
	The extreme MHD blast problem \cite{Balsara1999StaggeredMesh} is a benchmark test, featuring a strong central pressure in a uniform magnetized medium. The computational domain is $[-0.5,0.5]\times[-0.5,0.5]$ with periodic boundary conditions. The initial condition is given by
\begin{equation*}
\primitives =
\left\{
\begin{aligned}
&
(
1,0,0,0,
\dfrac{100}{\sqrt{4\pi}},0,0,
1000
),
\quad && r \le r_1,\\[1ex]
&
(
1,0,0,0,
\dfrac{100}{\sqrt{4\pi}},0,0,
0.1
),
\quad && r > r_1,
\end{aligned}
\right.
\end{equation*}
where
$ r = \sqrt{x^2 + y^2}, r_1 = 0.1, $ and $ \gamma = 1.4 $.
We simulate this problem until $t = 0.01$ with a uniform mesh of $ 200\times200 $.
The plasma beta of the fluid pulse is very small ($ \beta = p/\norm{\theB}^2 \approx 2.5\times10^{-4} $), which means that negative pressure may occur during the simulation.
Since this problem is very stringent \cite{Balsara1999StaggeredMesh} and no positivity-preserving limiter is applied, negative pressure may arise during the simulation \cite{Li2011CentralDiscontinuous, Fu2018GloballyDivergenceFree, Li2008HighOrder, Ding2024NewDiscretely}.
 If the negative pressure is ignored, the simulation can continue and produce the results shown in \cref{fig:extreme-blast}. 
The \hwenodf{} scheme exhibits superior numerical performance, while the \hwenobase{} scheme suffers from noticeable oscillations as depicted in \cref{fig:extreme-blast}.
Our results compare well with those reported in \cite{Balsara1999StaggeredMesh, Ding2024NewDiscretely, Liu2025EntropyStable}.
\begin{figure}[htbp]
	\centering
	\begin{subfigure}{.23\textwidth}
		\centering
		\includegraphics[width=\linewidth]{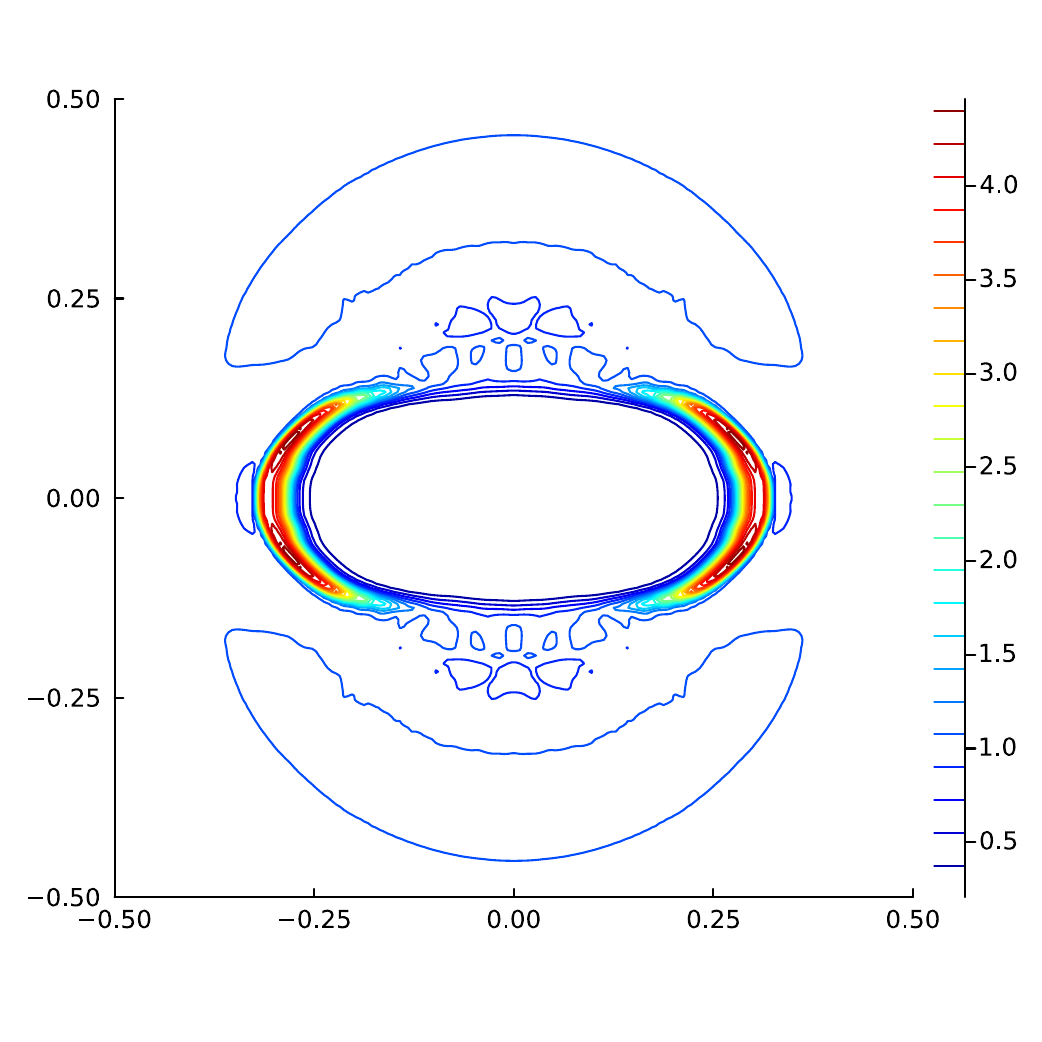}
		\caption{\hwenodf{}, $ \rho $}
	\end{subfigure}
	\hfill
	\begin{subfigure}{.23\textwidth}
		\centering
		\includegraphics[width=\linewidth]{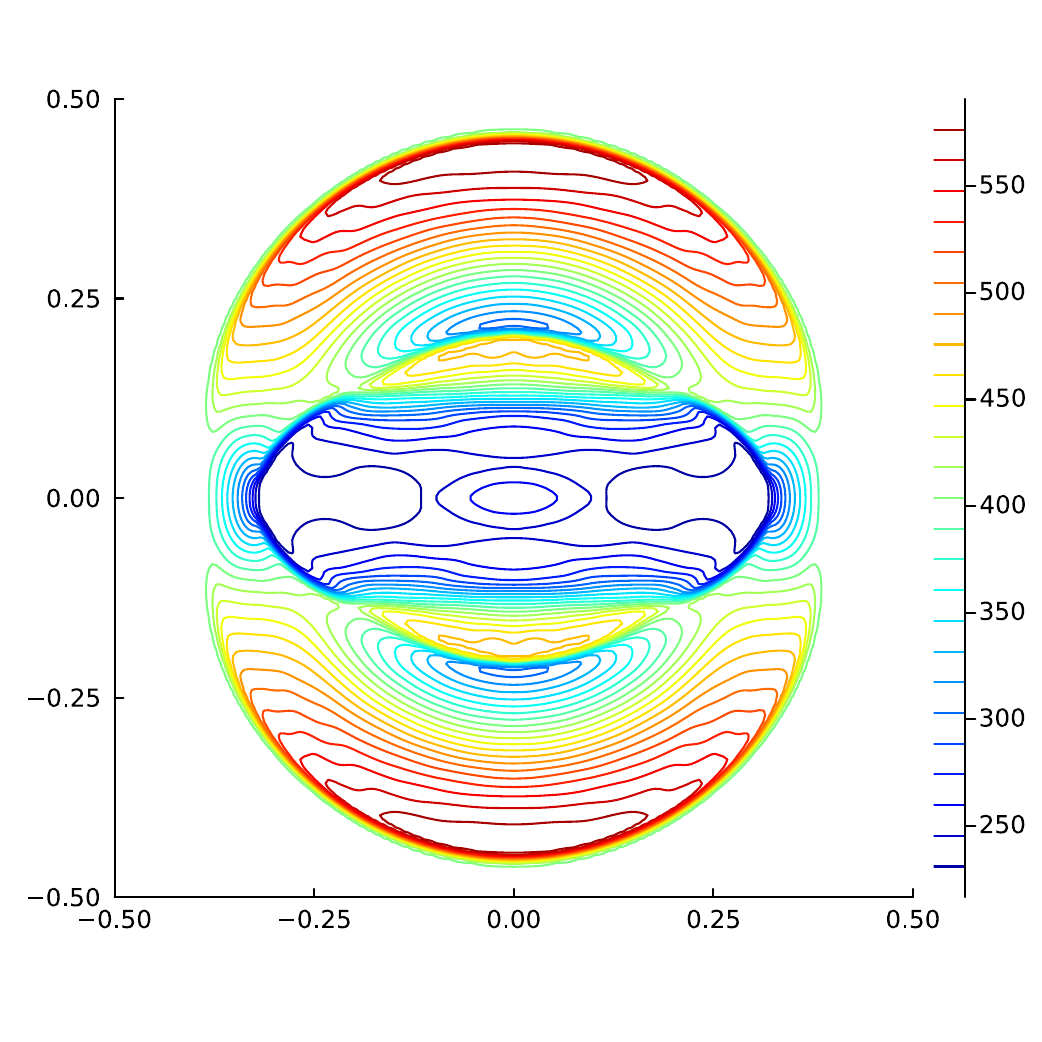}
		\caption{\hwenodf{}, $ \frac{1}{2}\norm{\theB}^{2} $}
	\end{subfigure}
	\hfill
	\begin{subfigure}{.23\textwidth}
		\centering
		\includegraphics[width=\linewidth]{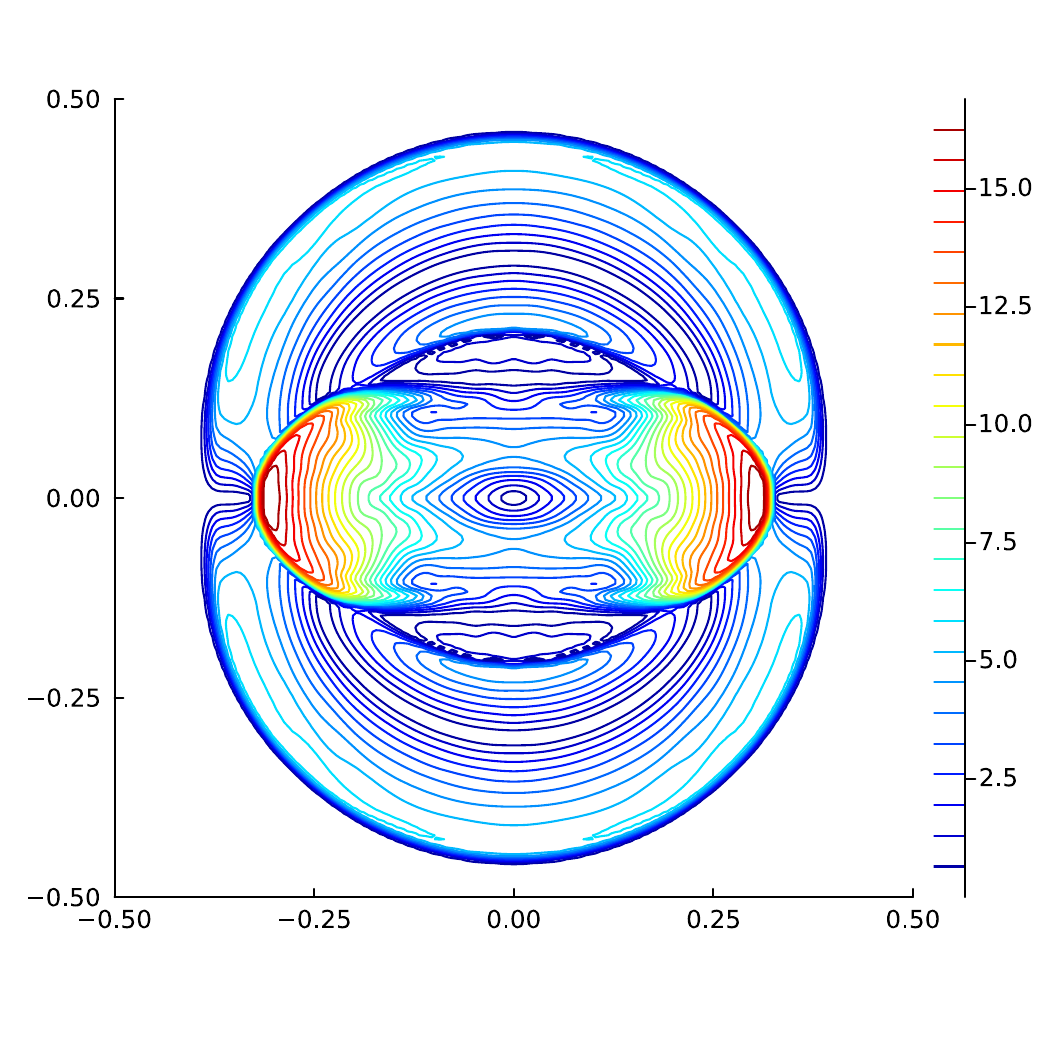}
		\caption{\hwenodf{}, $ \abs{\vec{u}} $}
	\end{subfigure}
	\hfill
	\begin{subfigure}{.23\textwidth}
		\centering
		\includegraphics[width=\linewidth]{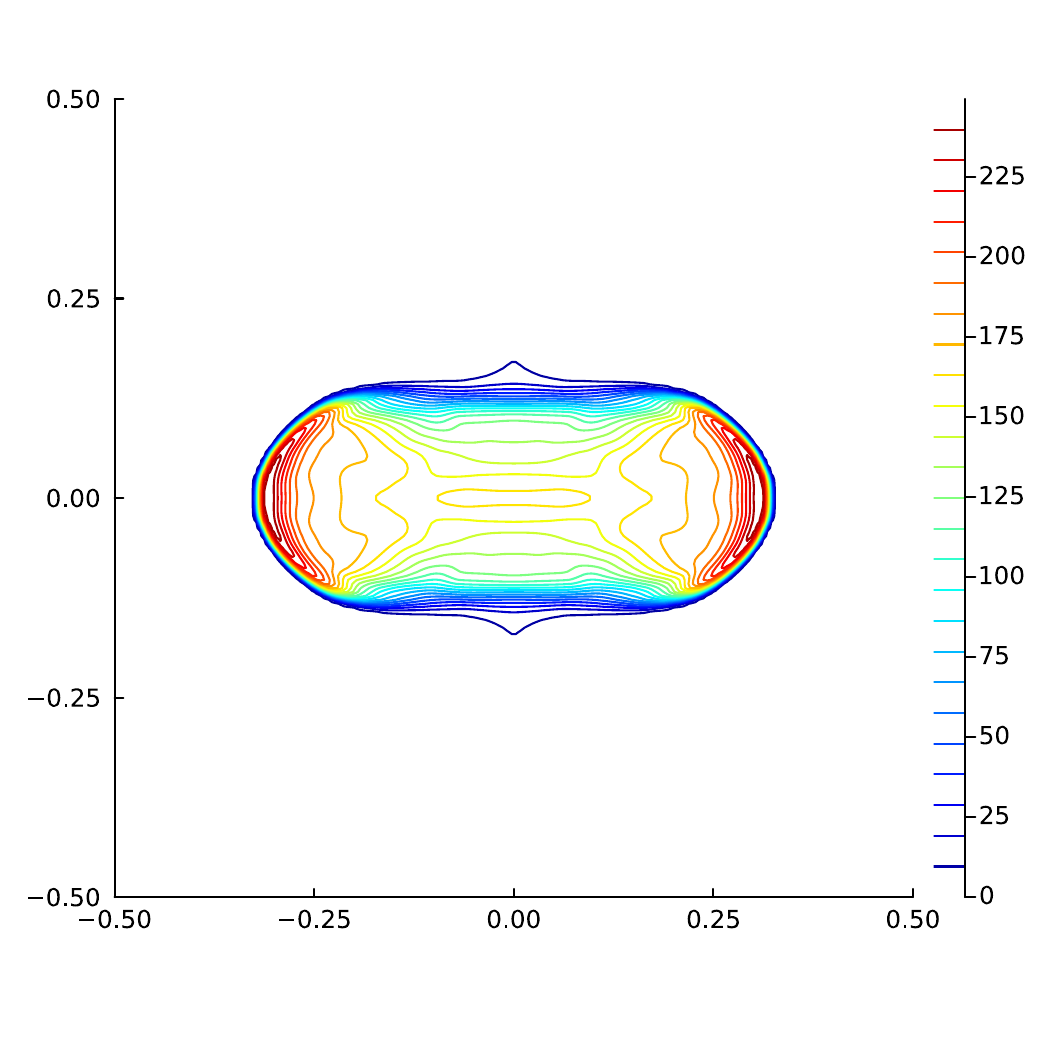}
		\caption{\hwenodf{}, $ p $}
	\end{subfigure}

	\begin{subfigure}{.23\textwidth}
		\centering
		\includegraphics[width=\linewidth]{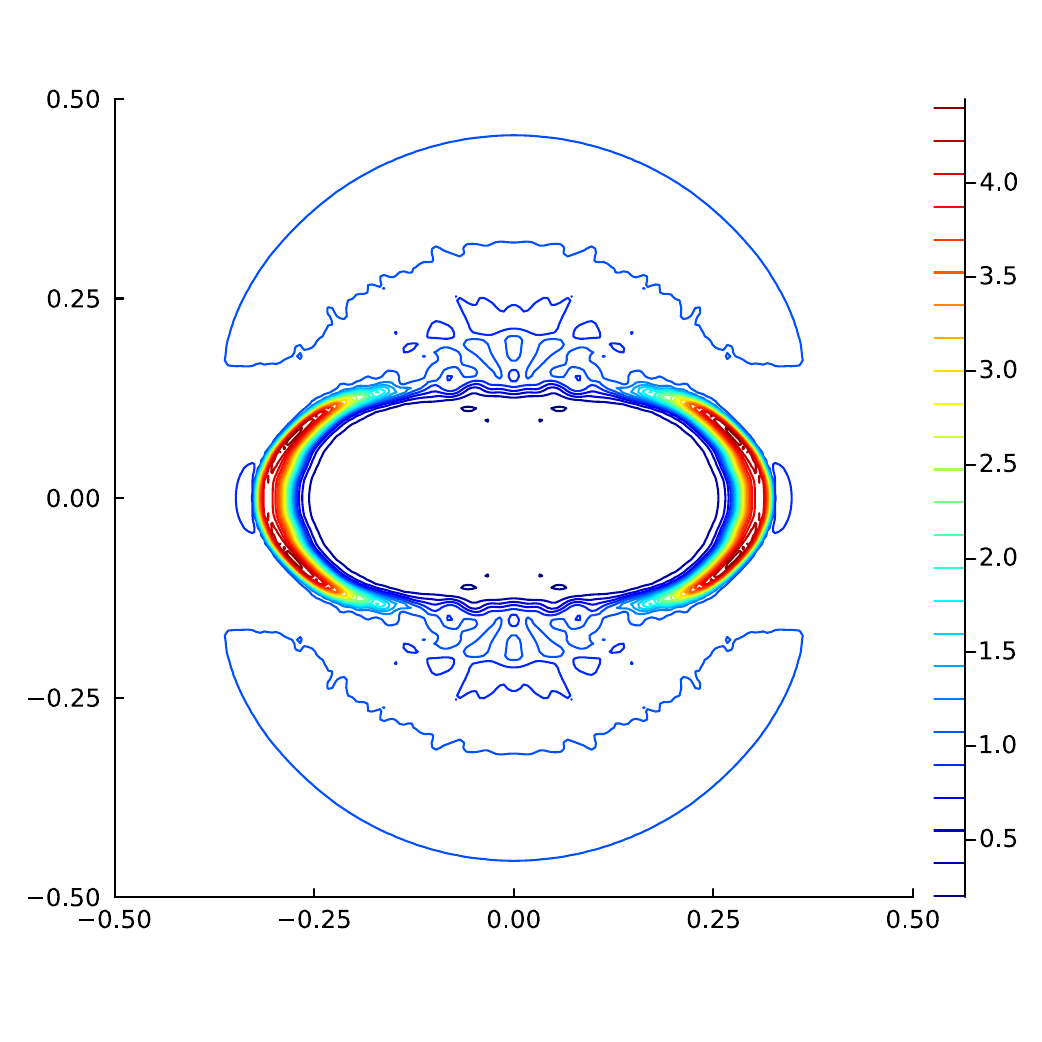}
		\caption{\hwenobase{}, $ \rho $}
	\end{subfigure}
	\hfill
	\begin{subfigure}{.23\textwidth}
		\centering
		\includegraphics[width=\linewidth]{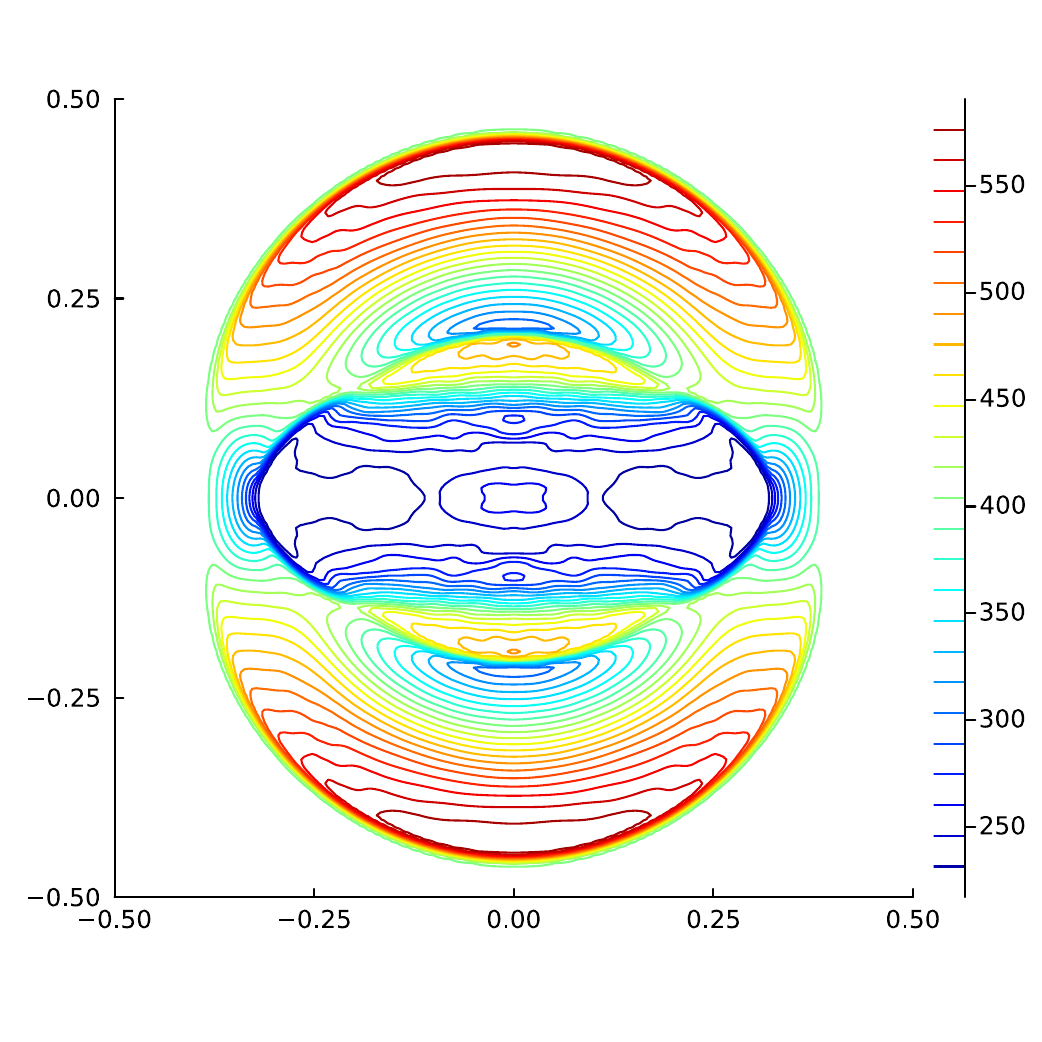}
		\caption{\hwenobase{}, $ \frac{1}{2}\norm{\theB}^{2} $}
	\end{subfigure}
	\hfill
	\begin{subfigure}{.23\textwidth}
		\centering
		\includegraphics[width=\linewidth]{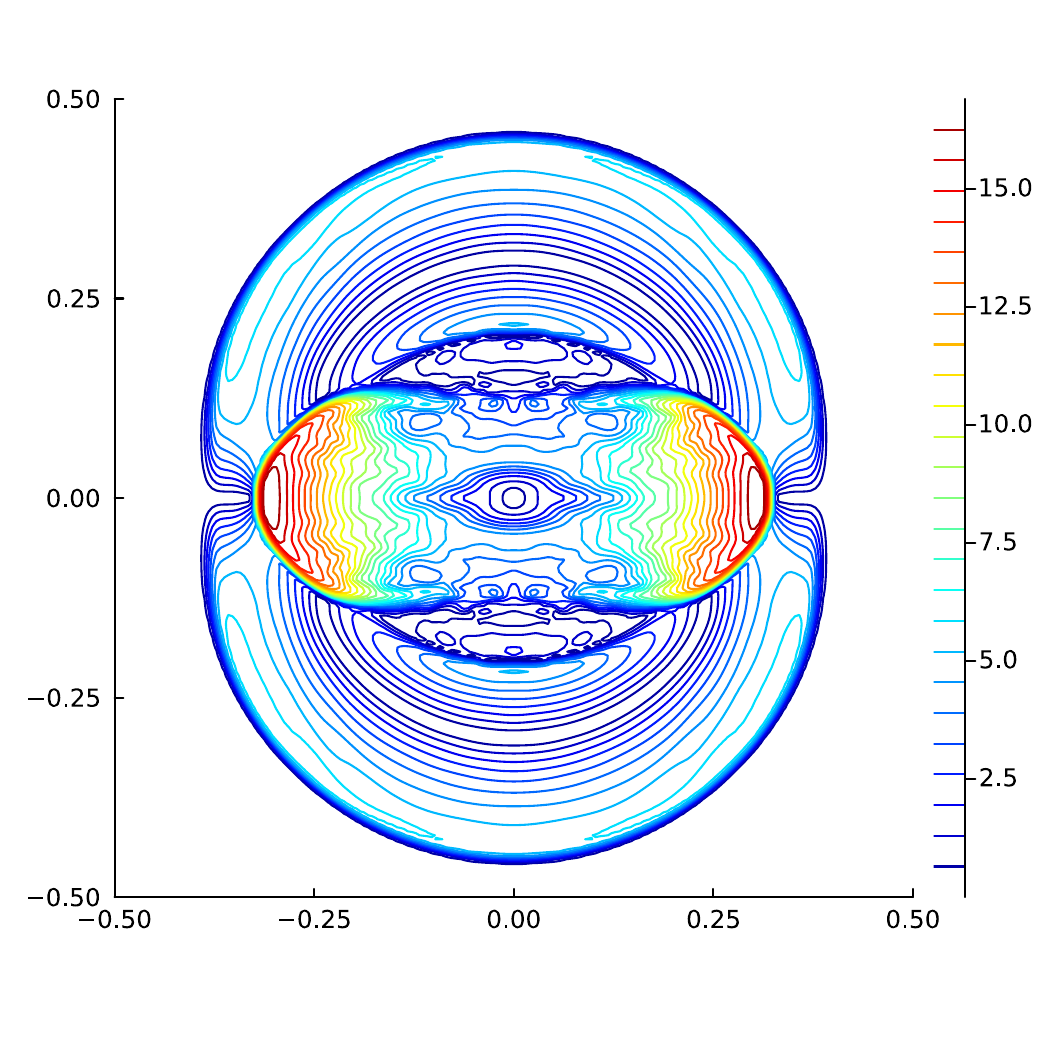}
		\caption{\hwenobase{}, $ \abs{\vec{u}} $}
	\end{subfigure}
	\hfill
	\begin{subfigure}{.23\textwidth}
		\centering
		\includegraphics[width=\linewidth]{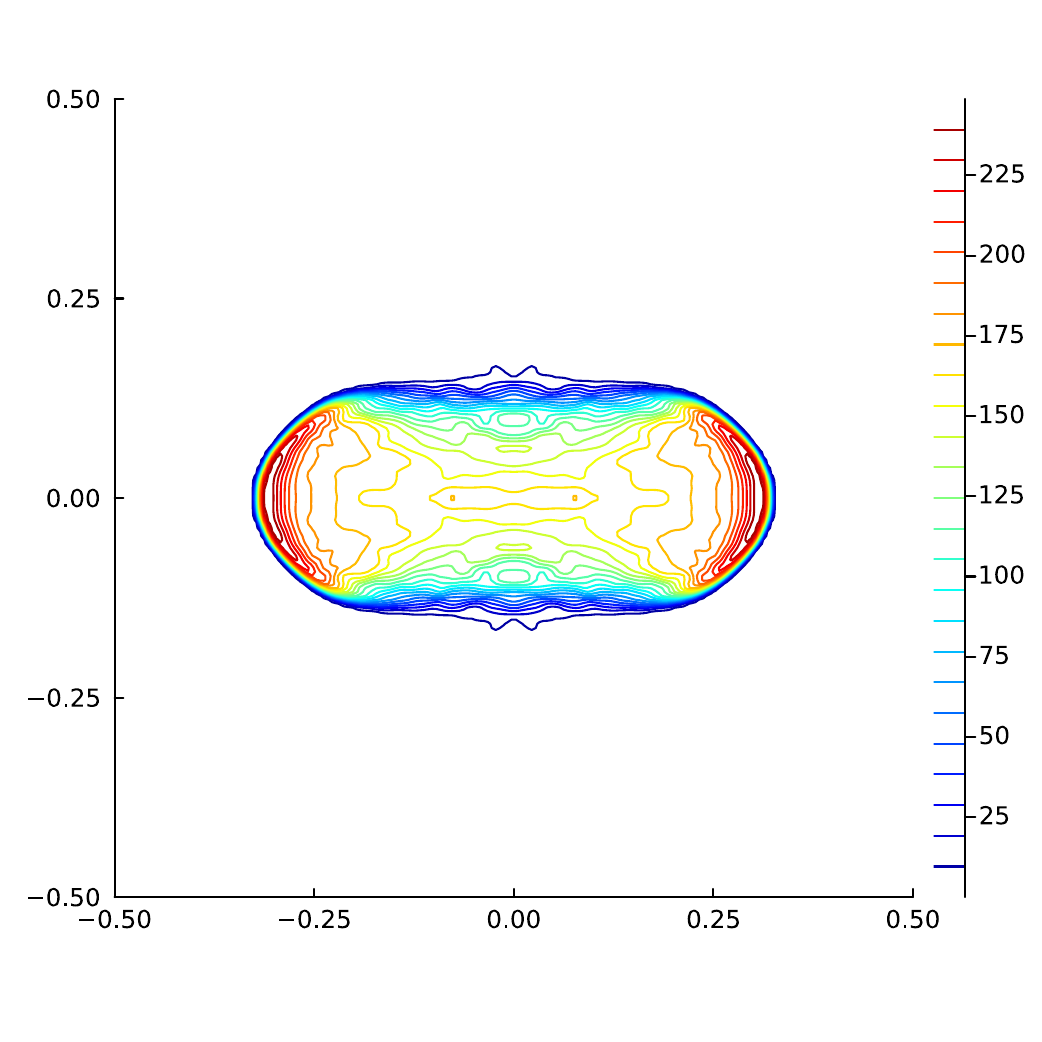}
		\caption{\hwenobase{}, $ p $}
	\end{subfigure}
	\caption{\cref{ex:extreme-blast}. Extreme blast, simulated by \hwenodf{} (top row) and \hwenobase{} (bottom row) with $200\times 200$ meshes. Contour plots of density, magnetic pressure, velocity magnitude, and thermal pressure at $t=0.01$. 25 equally spaced contour lines are plotted.}
	\label{fig:extreme-blast}
\end{figure}
\end{example}
\begin{example}[Cloud shock interaction]
	\label{ex:cloud-shock}
	The cloud shock interaction \cite{Dai1998SimpleFinitea, Toth2000B0Constraint} problem considers an MHD shock propagating toward a stationary dense cloud. The computational domain is $[0,1]\times[0,1]$ with $\gamma=5/3$. The initial primitive variables
are given by
\begin{equation*}
\primitives=
\left\{
\begin{array}{@{}l@{\hspace{2pt}}l@{}}
(3.86859,11.2536,0,0,0,2.1826182,-2.1826182,167.345), & x{\leqslant}0.05,
\\
(1,0,0,0,0,0.56418958,0.56418958,1), & x{>}0.05.
\end{array}
\right.
\end{equation*}
A circular cloud of density $\rho=10$ and radius $r=0.15$ is centered at $(0.25,0.5)$, namely
\begin{equation*}
\rho(0,x,y)=10,
\quad
(x-0.25)^2+(y-0.5)^2\le 0.15^2.
\end{equation*}
The boundary conditions are set as inflow at $x=0$ and outflow on the other three sides. The solution is evolved until $t=0.06$ with $600\times 600$ meshes.
In \cref{figs:cloud-shock-contour}, the \hwenodf{} scheme resolves the interaction between the shock and the dense cloud.
As shown in \cref{figs:cloud-shock-contour}, slight oscillations can be observed in the \hwenobase{} scheme near shock interfaces, while the \hwenodf{} scheme captures the shock structure without visible oscillations.
The numerical results of \hwenodf{} agree well with those in the literature \cite{Christlieb2014FiniteDifference, Duan2025ActiveFlux, Liu2025GloballyDivergencefreea}.
This example further verifies that the proposed divergence-free correction is effective in suppressing spurious oscillations and improving the robustness of the scheme.

\begin{figure}[htbp]
	\centering
	\begin{subfigure}{.23\textwidth}
		\centering
		\includegraphics[width=\linewidth]{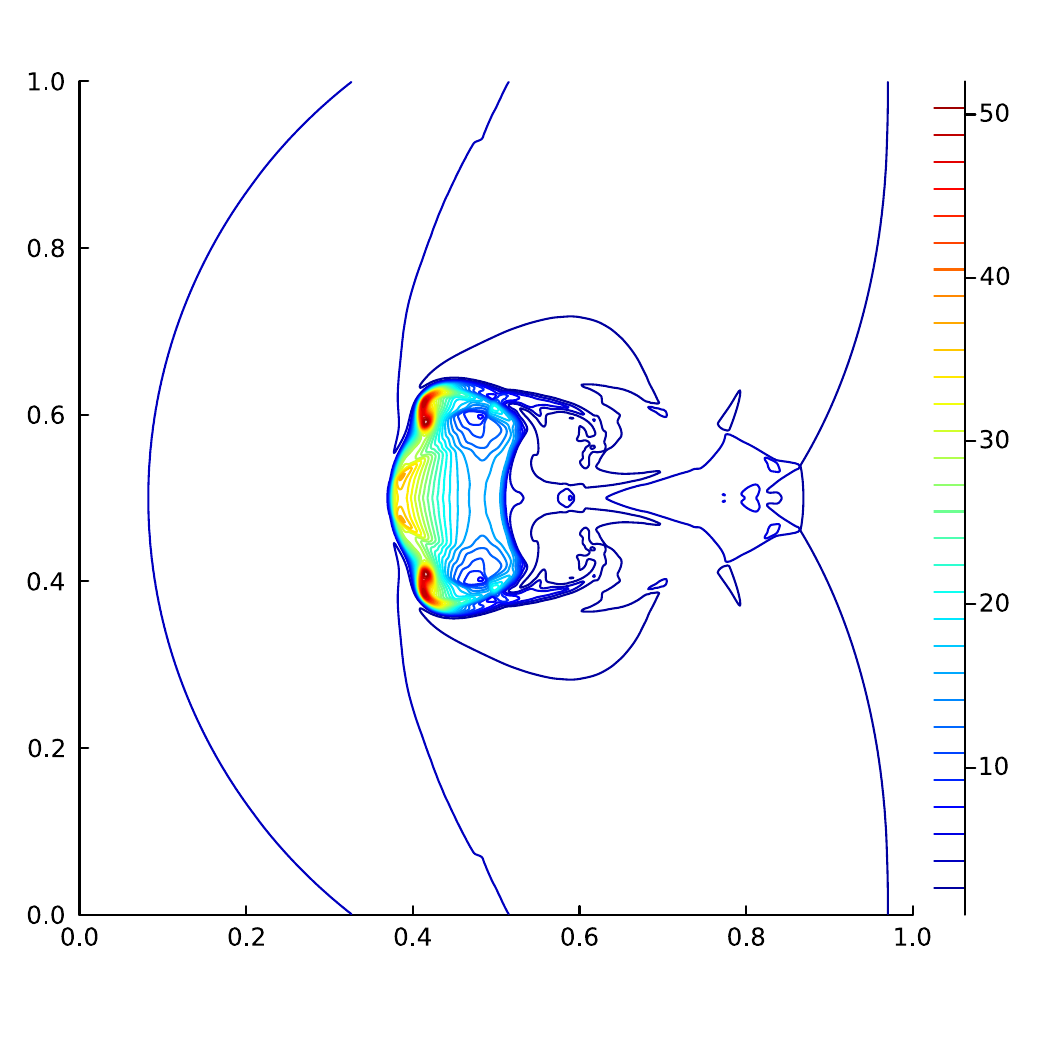}
		\caption{\hwenodf{}, $ \rho $}
	\end{subfigure}
	\hfill
	\begin{subfigure}{.23\textwidth}
		\centering
		\includegraphics[width=\linewidth]{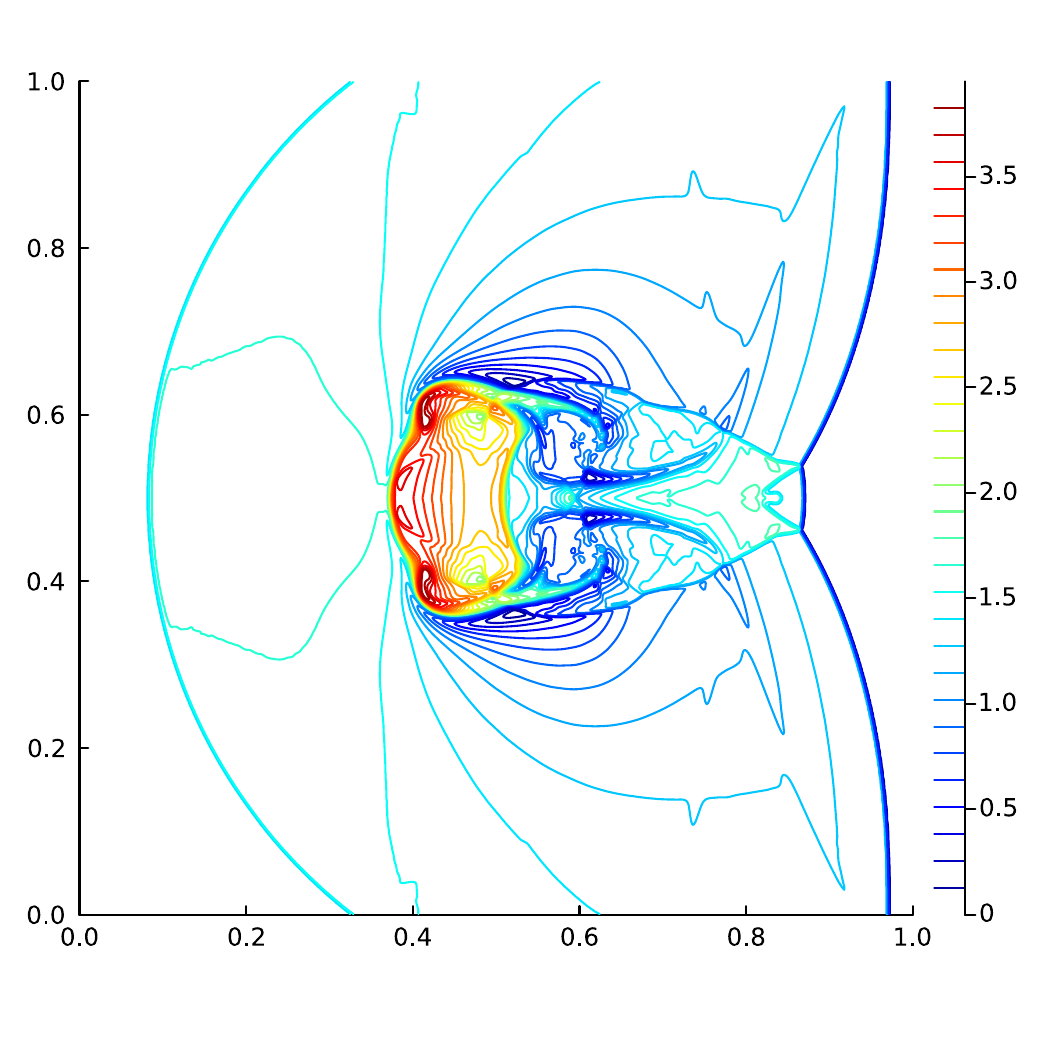}
		\caption{\hwenodf{}, $ \ln \rho $}
	\end{subfigure}
	\hfill
	\begin{subfigure}{.23\textwidth}
		\centering
		\includegraphics[width=\linewidth]{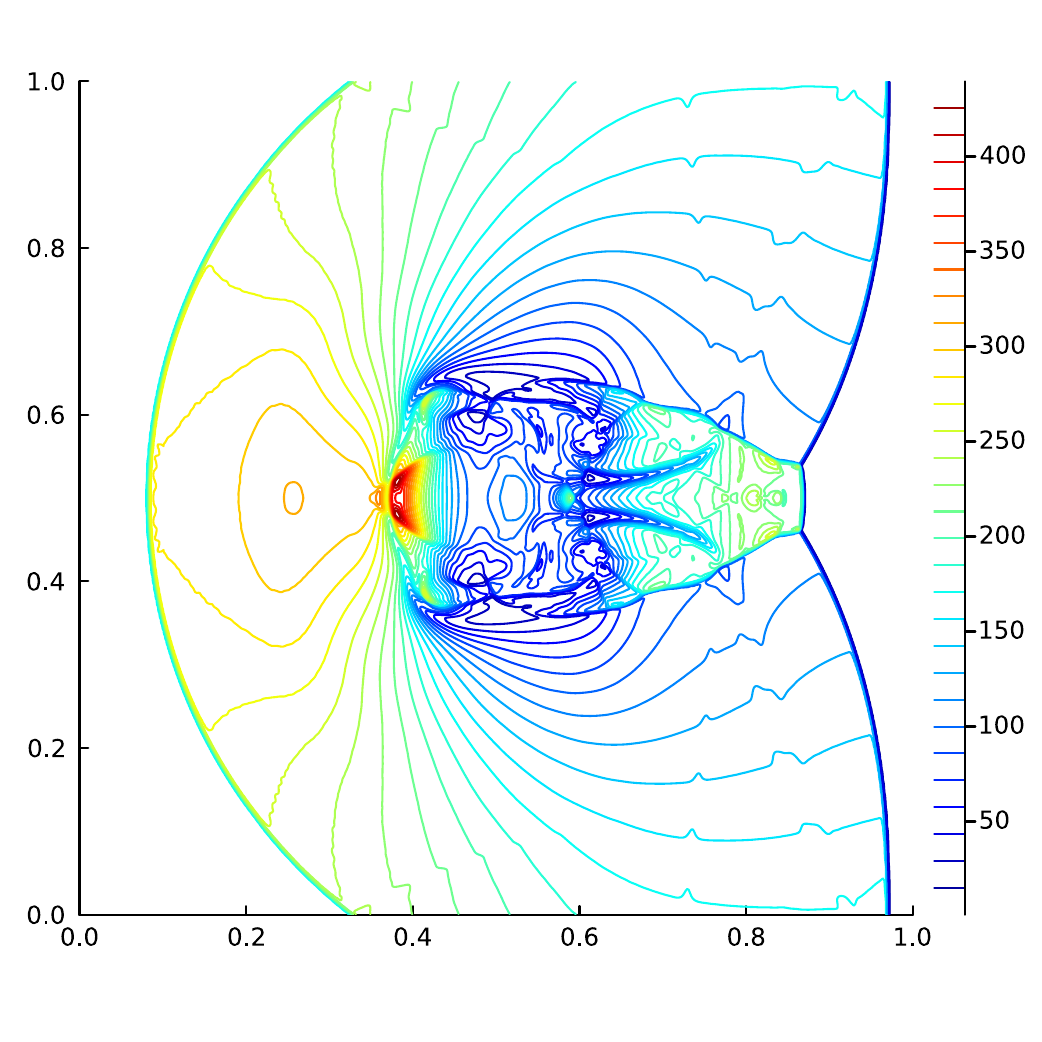}
		\caption{\hwenodf{}, $ p $}
	\end{subfigure}
	\hfill
	\begin{subfigure}{.23\textwidth}
		\centering
		\includegraphics[width=\linewidth]{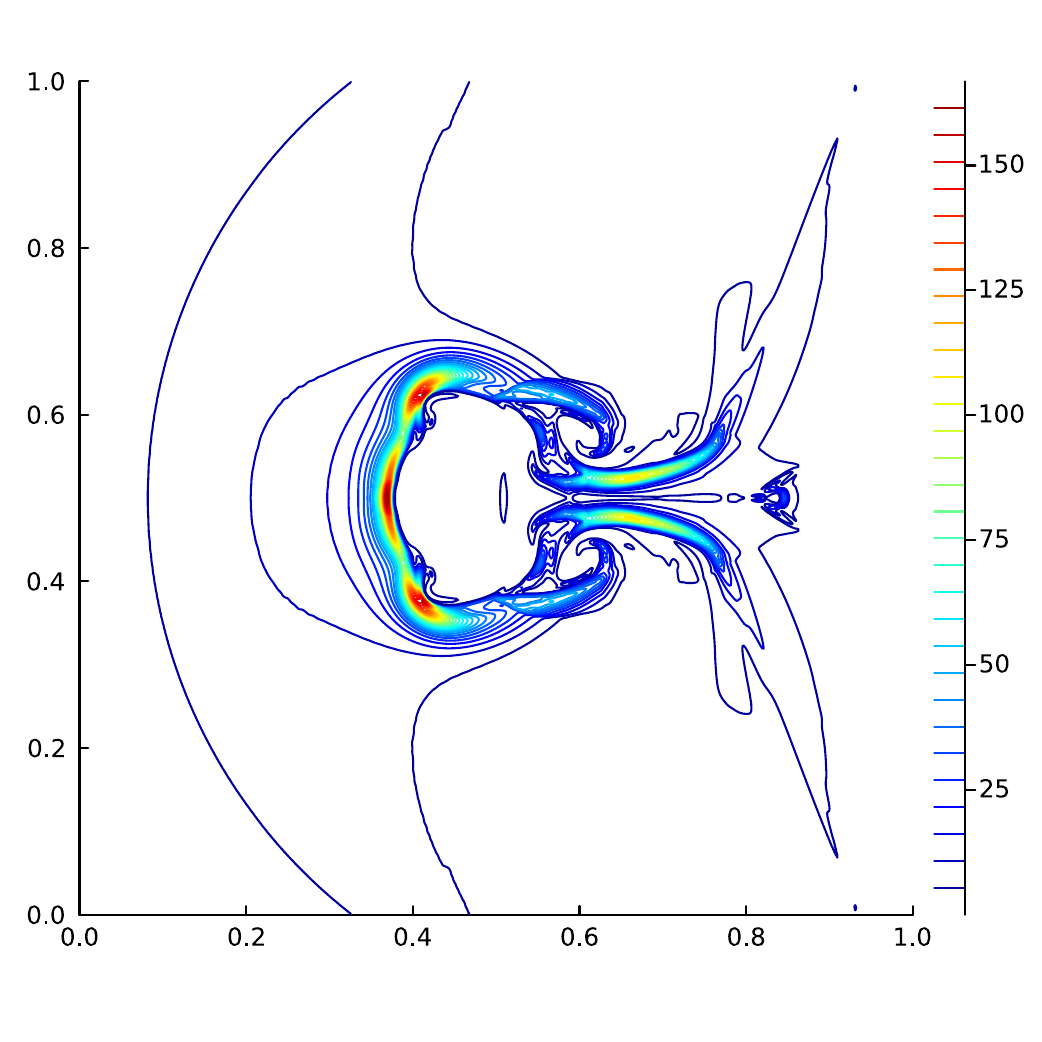}
		\caption{\hwenodf{}, $ \frac{1}{2} {\norm{\theB}}^{2} $}
	\end{subfigure}

	\begin{subfigure}{.23\textwidth}
		\centering
		\includegraphics[width=\linewidth]{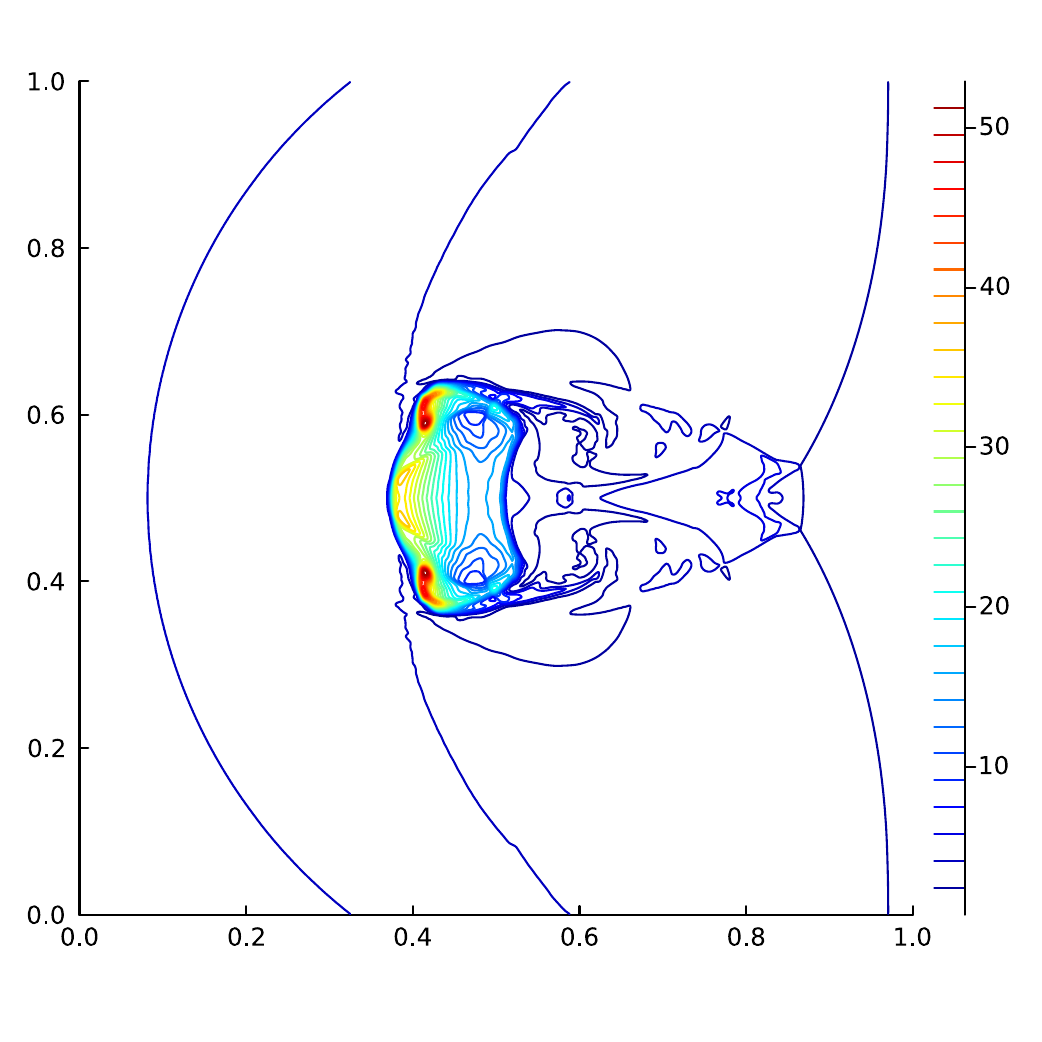}
		\caption{\hwenobase{}, $ \rho $}
	\end{subfigure}
	\hfill
	\begin{subfigure}{.23\textwidth}
		\centering
		\includegraphics[width=\linewidth]{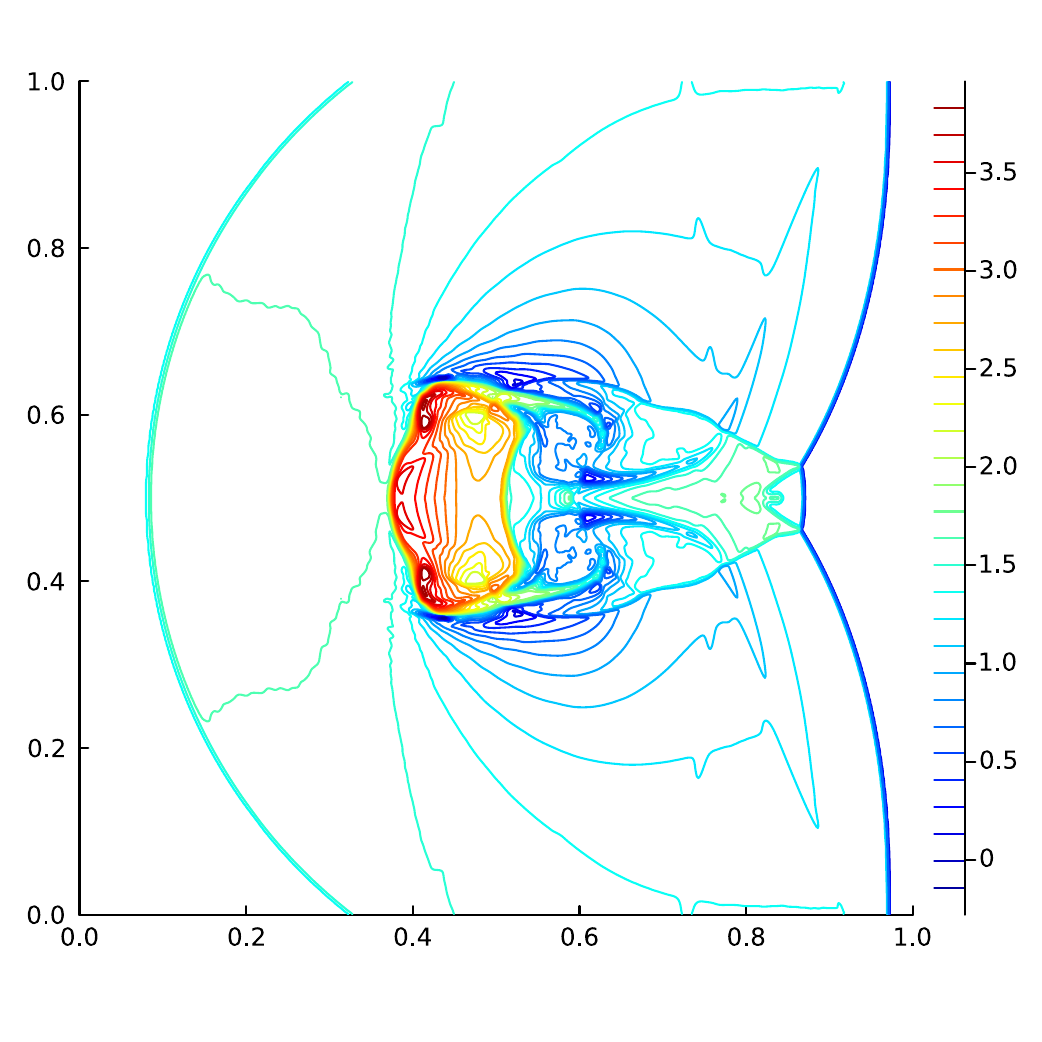}
		\caption{\hwenobase{}, $ \ln \rho $}
	\end{subfigure}
	\hfill
	\begin{subfigure}{.23\textwidth}
		\centering
		\includegraphics[width=\linewidth]{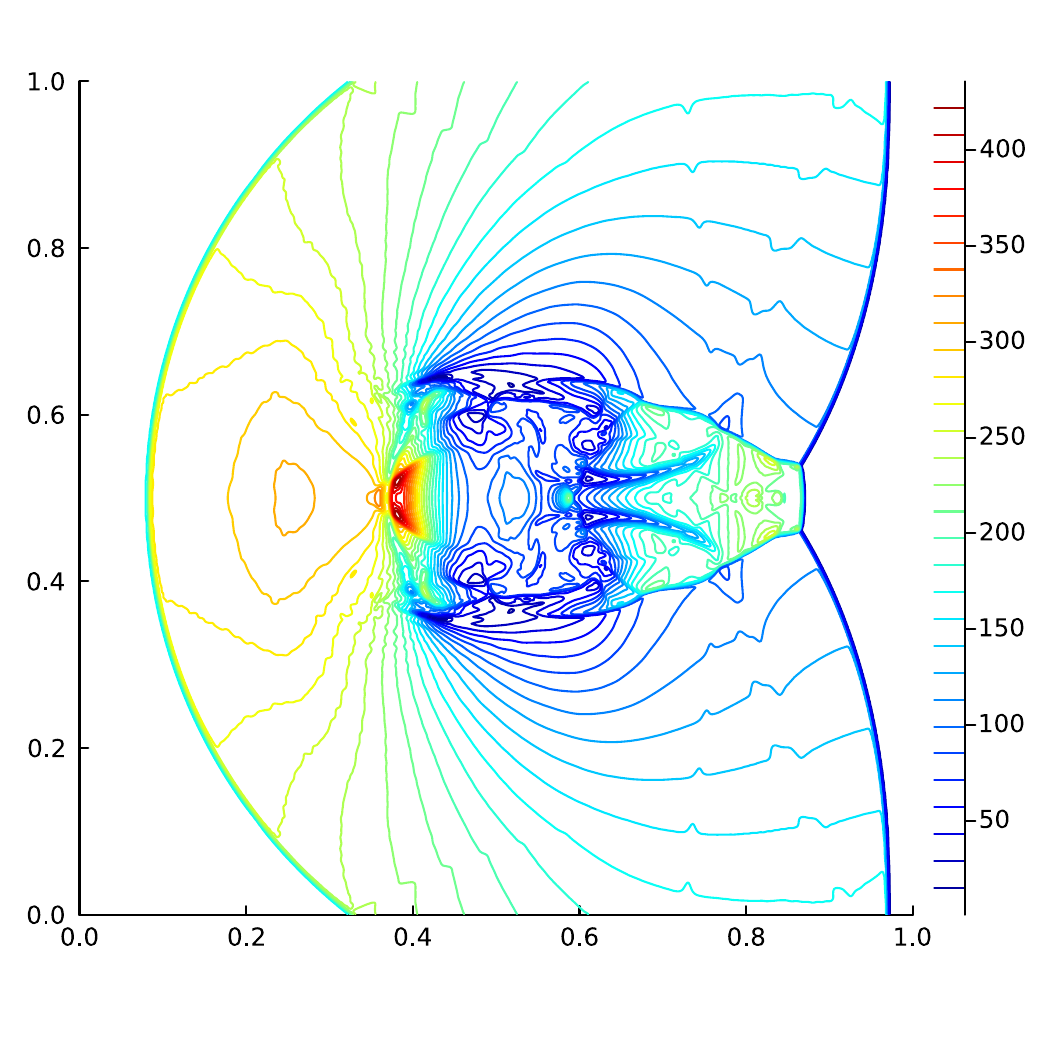}
		\caption{\hwenobase{}, $ p $}
	\end{subfigure}
	\hfill
	\begin{subfigure}{.23\textwidth}
		\centering
		\includegraphics[width=\linewidth]{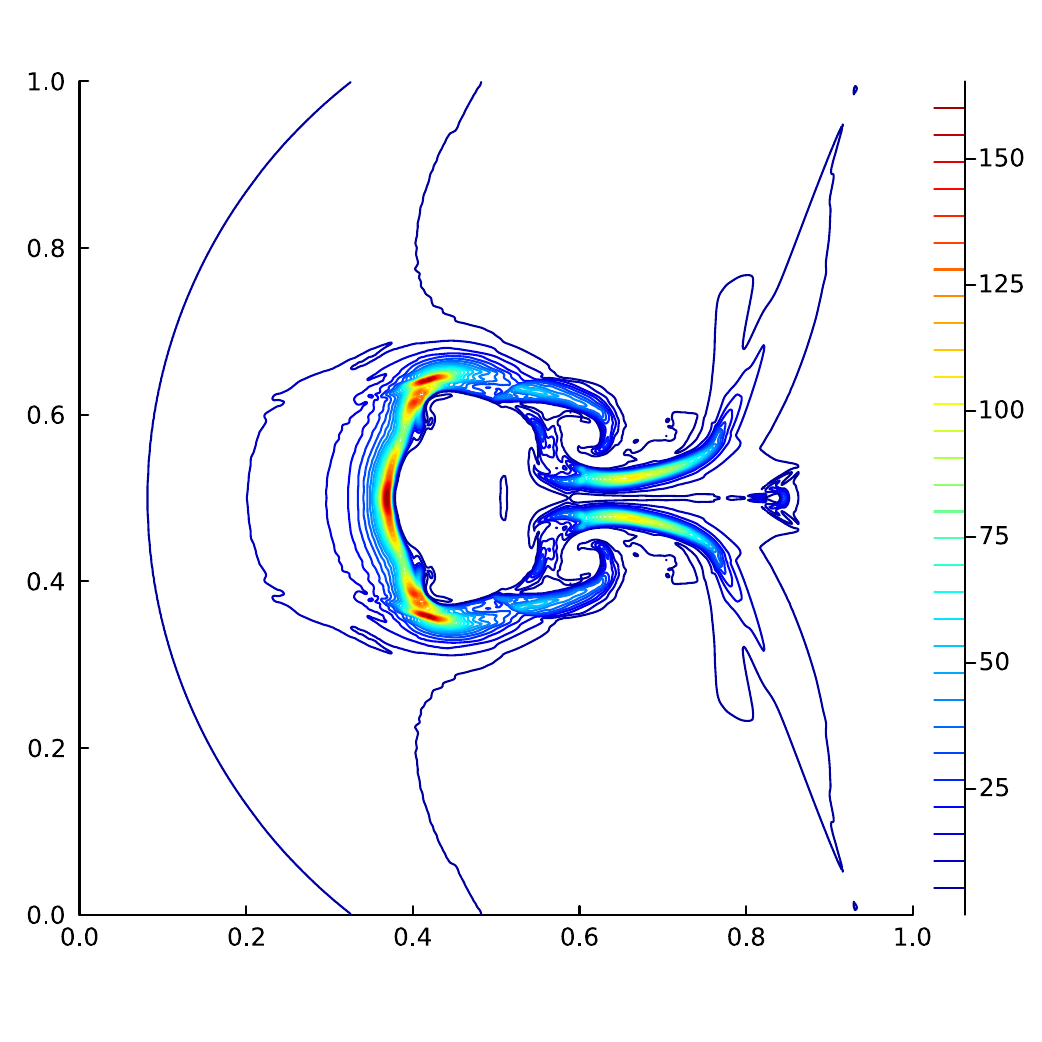}
		\caption{\hwenobase{}, $ \frac{1}{2} {\norm{\theB}}^{2} $}
	\end{subfigure}
	\caption{\cref{ex:cloud-shock}. Cloud--shock interaction, simulated by \hwenodf{} (top row) and \hwenobase{} (bottom row) with $600\times 600$ meshes. Contour lines of density, log density, thermal pressure, and magnetic pressure at $t=0.06$.}
	\label{figs:cloud-shock-contour}
\end{figure}
\end{example}

\begin{example}[Kelvin-Helmholtz instability]
	\label{ex:kelvin-helmholtz}
	Consider the Kelvin-Helmholtz instability problem in Yee and Sj\"ogreen \cite{Yee2005DivergenceFree}.
The computational domain is $[0,1]\times[-1,1]$ with periodic boundary conditions.
The initial primitive variables are given by
\begin{equation*}
\begin{aligned}
&\rho = 1, \quad
u^{1} = 5(
\tanh\!\bigl(20(y+0.5)\bigr)
-\tanh\!\bigl(20(y-0.5)\bigr)-1
), \\
&u^{2} = 0.25\sin(2\pi x)
(
e^{-100(y+0.5)^2}
-e^{-100(y-0.5)^2}
), \quad u^{3} = 0,\\
&\theB^\top = (1,0,0), \quad p = 50.
\end{aligned}
\end{equation*}
The adiabatic index is $ \gamma = 1.4.  $
The solution is evolved until $t=0.5$ with $ 256\times 512 $ meshes.
In \cref{figs:kelvin-helmholtz}, the \hwenodf{} scheme resolves the complex features of 
the magnetic field frozen into the fluid and amplified by the vortical flow and the 
results are consistent with those in the literature \cite{Yee2005DivergenceFree, Liu2025EntropyStable}. In contrast, the \hwenobase{} scheme exhibits noticeable oscillations and nonphysical distortions in the solution, which can be attributed to the violation of the divergence-free constraint.
This example further demonstrates the effectiveness of the proposed divergence-free correction in suppressing spurious oscillations and improving the robustness of the scheme.
\begin{figure}[htbp]
	\centering
	\begin{subfigure}{.23\textwidth}
		\centering
		\includegraphics[width=\linewidth]{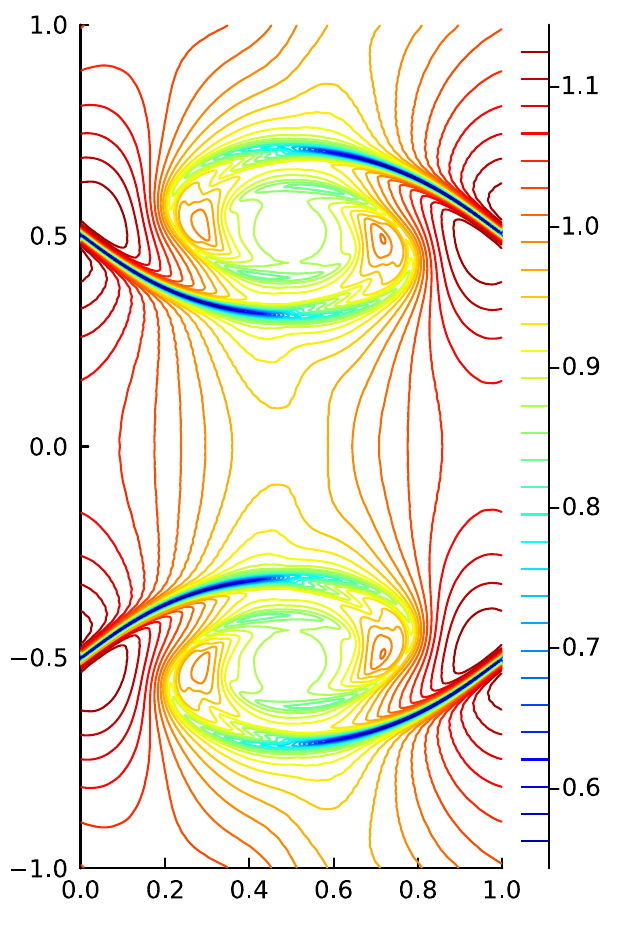}
		\caption{\hwenodf{}, $ \rho $}
	\end{subfigure}
	\hfill
	\begin{subfigure}{.23\textwidth}
		\centering
		\includegraphics[width=\linewidth]{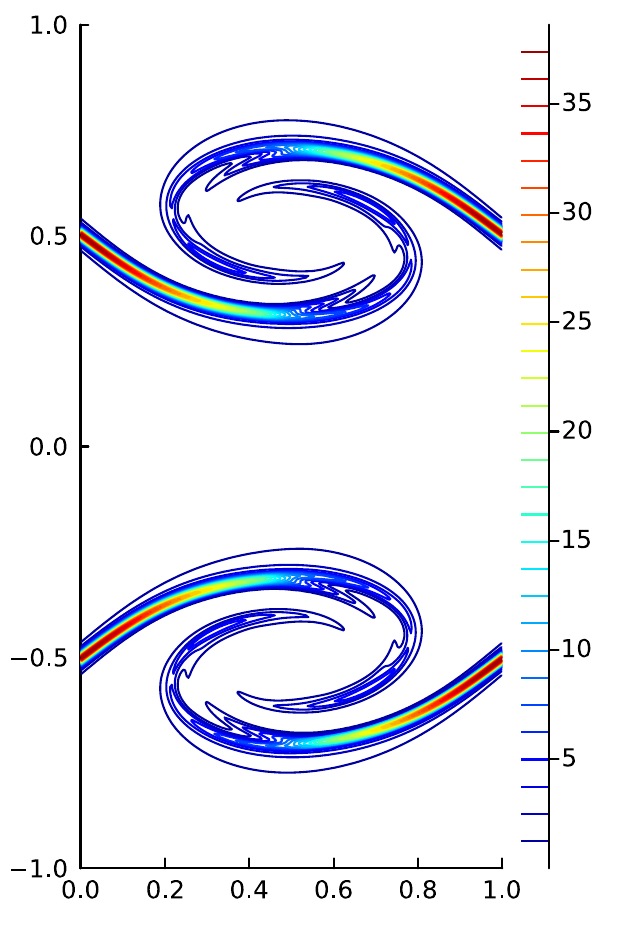}
		\caption{\hwenodf{}, $ \frac{1}{2} \norm{\theB}^{2} $}
	\end{subfigure}
	\hfill
	\begin{subfigure}{.23\textwidth}
		\centering
		\includegraphics[width=\linewidth]{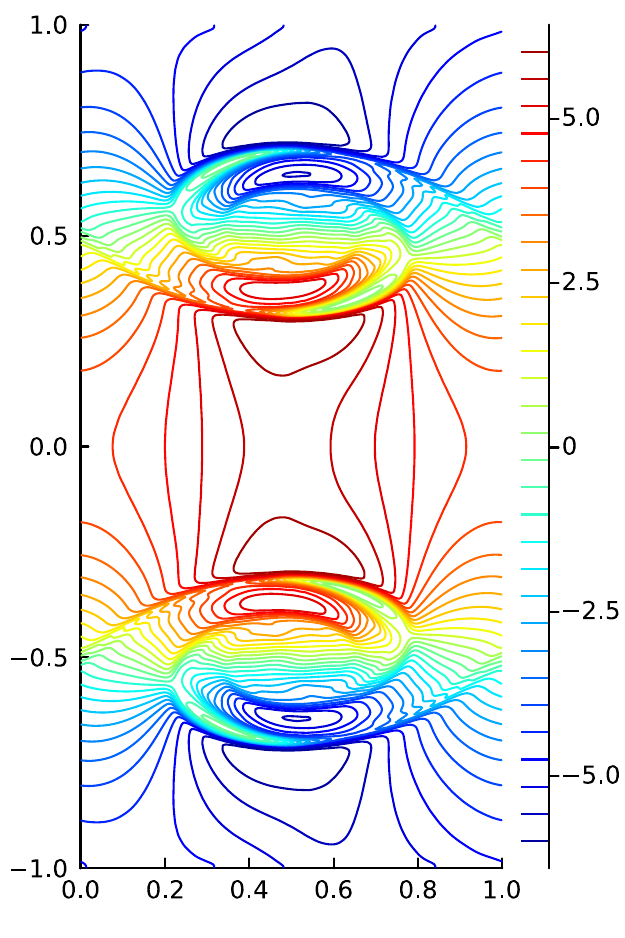}
		\caption{\hwenodf{}, $ u^{1} $}
	\end{subfigure}
	\hfill
	\begin{subfigure}{.23\textwidth}
		\centering
		\includegraphics[width=\linewidth]{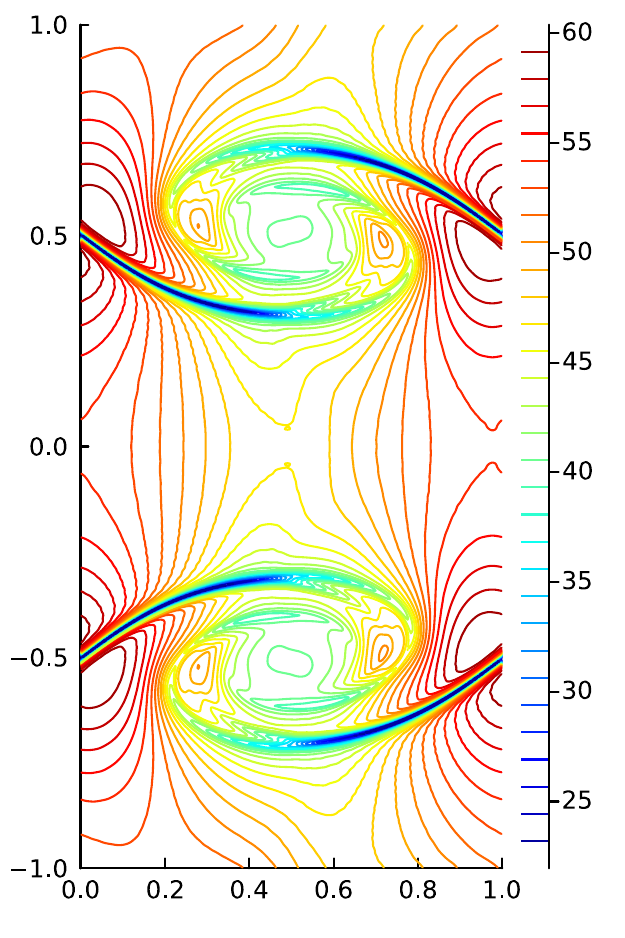}
		\caption{\hwenodf{}, $ p $}
	\end{subfigure}

	\begin{subfigure}{.23\textwidth}
		\centering
		\includegraphics[width=\linewidth]{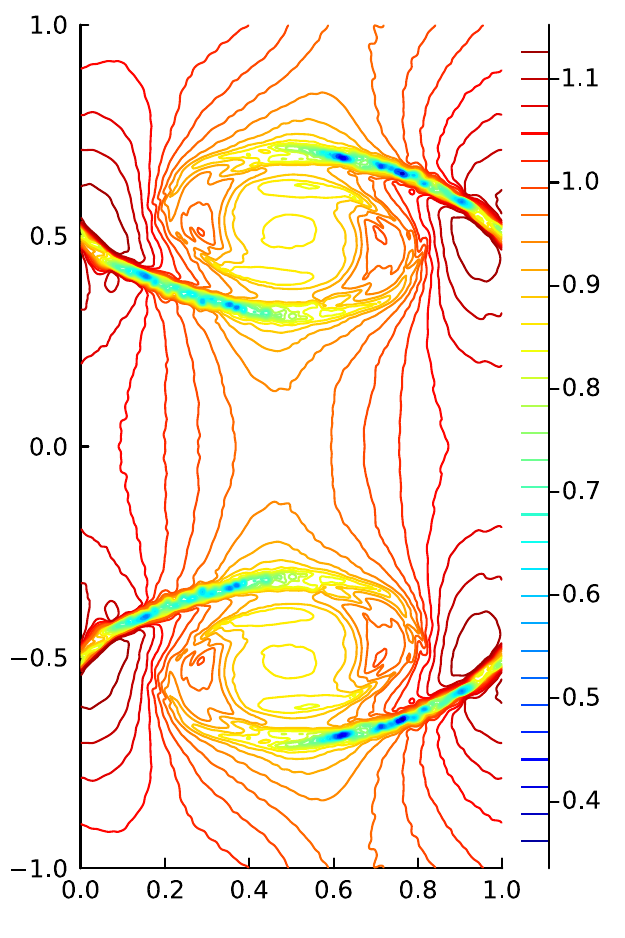}
		\caption{\hwenobase{}, $ \rho $}
	\end{subfigure}
	\hfill
	\begin{subfigure}{.23\textwidth}
		\centering
		\includegraphics[width=\linewidth]{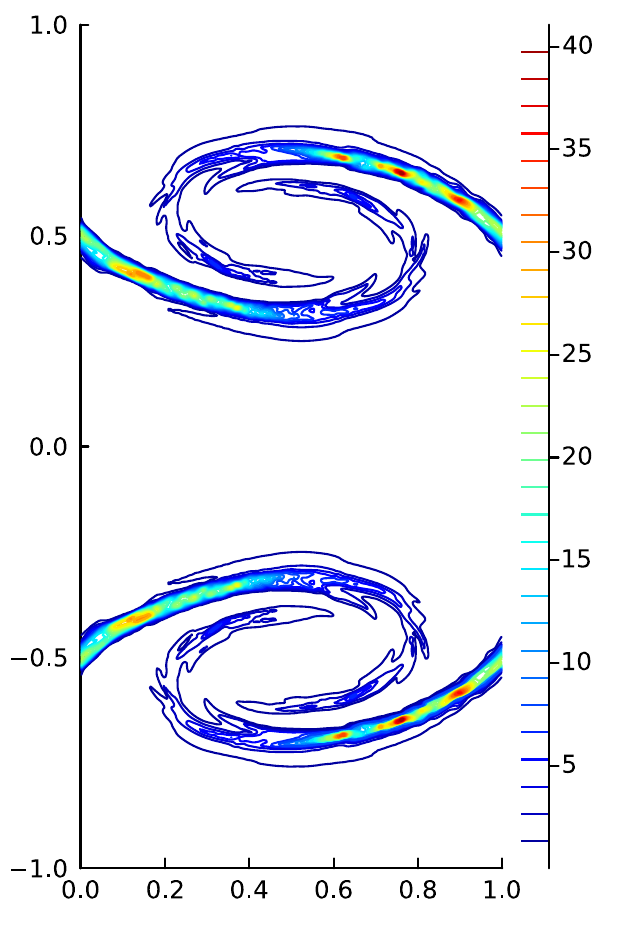}
		\caption{\hwenobase{}, $ \frac{1}{2} \norm{\theB}^{2} $}
	\end{subfigure}
	\hfill
	\begin{subfigure}{.23\textwidth}
		\centering
		\includegraphics[width=\linewidth]{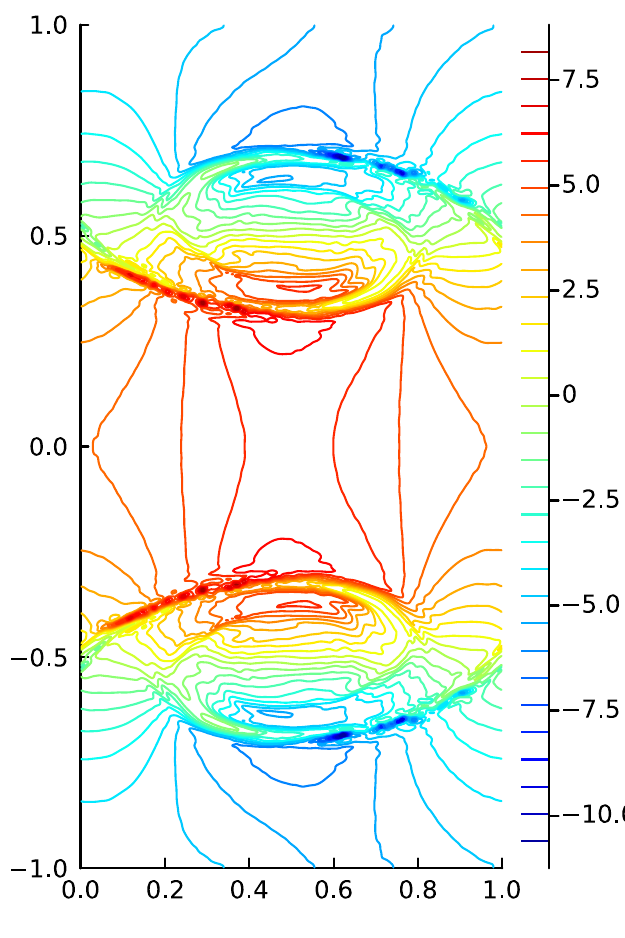}
		\caption{\hwenobase{}, $ u^{1} $}
	\end{subfigure}
	\hfill
	\begin{subfigure}{.23\textwidth}
		\centering
		\includegraphics[width=\linewidth]{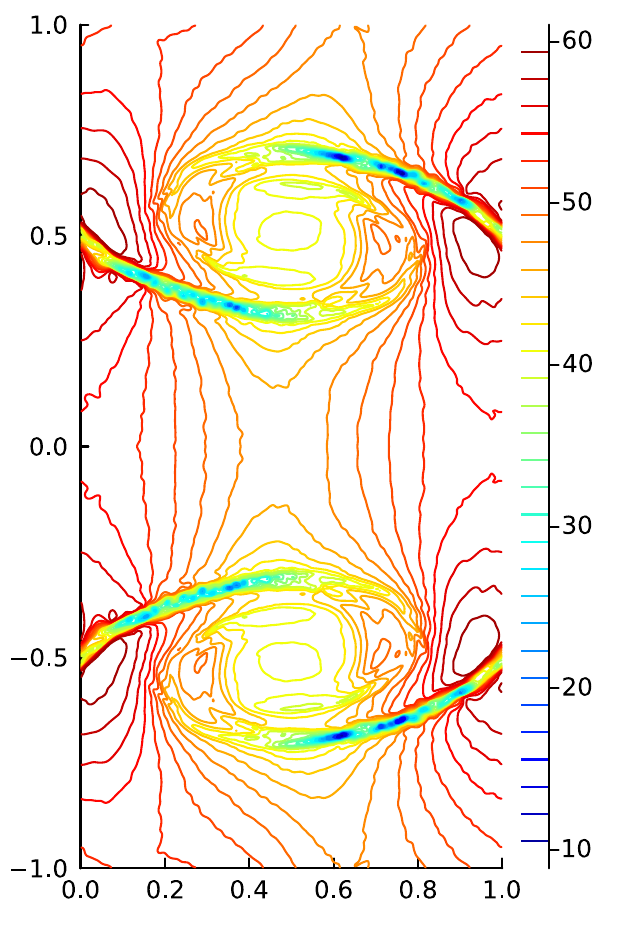}
		\caption{\hwenobase{}, $ p $}
	\end{subfigure}
	\caption{\cref{ex:kelvin-helmholtz}. Kelvin--Helmholtz instability, simulated by \hwenodf{} (top row) and \hwenobase{} (bottom row) with $256\times 512$ meshes. Contour plots of density, magnetic pressure, $x$-velocity, and thermal pressure at $t=0.5$. 30 equally spaced contour lines are plotted.}
	\label{figs:kelvin-helmholtz}
\end{figure}
\end{example}
\begin{example}[The field loop advection]
	\label{ex:field-loop}
	The field loop advection \cite{Gardiner2005UnsplitGodunov} problem is considered. The computational domain is 
$[-1,1]\times[-0.5,0.5]$ with periodic boundary conditions.
The initial condition is given by
\begin{equation*}
(\rho,u^{1},u^{2},u^{3},B^{3},p)=(1,2,1,1,0,1),\quad
(B^{1},B^{2})=\curl A_z,
\end{equation*}
where
\begin{equation*}
A_z=
\left\{
\begin{aligned}
&A_0(r_0-r),
\quad && r\le r_0,
\\
&0,
\quad && r>r_0,
\end{aligned}
\right.
\end{equation*}
with
$ A_0=10^{-3},
r_0=0.3,
r=\sqrt{x^2+y^2}. $
This initial condition describes a circular region centered at the origin with
a magnetic field loop generated by the magnetic potential $A_z$, while
the magnetic field is zero outside this region. 
The initial magnetic field is shown in \cref{fig:field-loop-initial}.
Under the periodic boundary
conditions, $t=2$ corresponds to one period, and the magnetic field
should return to its initial profile.
The solution is evolved until $t=10$. 
 As shown in \cref{fig:field-loop-comparison}, the magnetic field lines simulated by the
\hwenobase{} scheme exhibit slight distortion at $t = 2$ and severe distortion at
$t = 10$. Moreover, the magnitude of the magnetic field suffers significant
dissipation, indicating that the \hwenobase{} scheme is not stable for this
simulation. In contrast, the structures of magnetic field lines are well preserved by the
\hwenodf{} scheme at both $t = 2$ and $t = 10$, and the magnitude of the magnetic
field exhibits only slight dissipation, which agrees with \cite{Li2011CentralDiscontinuous, Fu2018GloballyDivergenceFree}.
Furthermore, the results of the \hwenodf{} scheme are slightly better than those in \cite{Li2011CentralDiscontinuous, Fu2018GloballyDivergenceFree}.
These results further demonstrate the effectiveness of the proposed approach for 
handling the divergence-free constraint.
\begin{figure}[htbp]
	\begin{subfigure}{.45\textwidth}
		\centering
		\includegraphics[width=\linewidth]{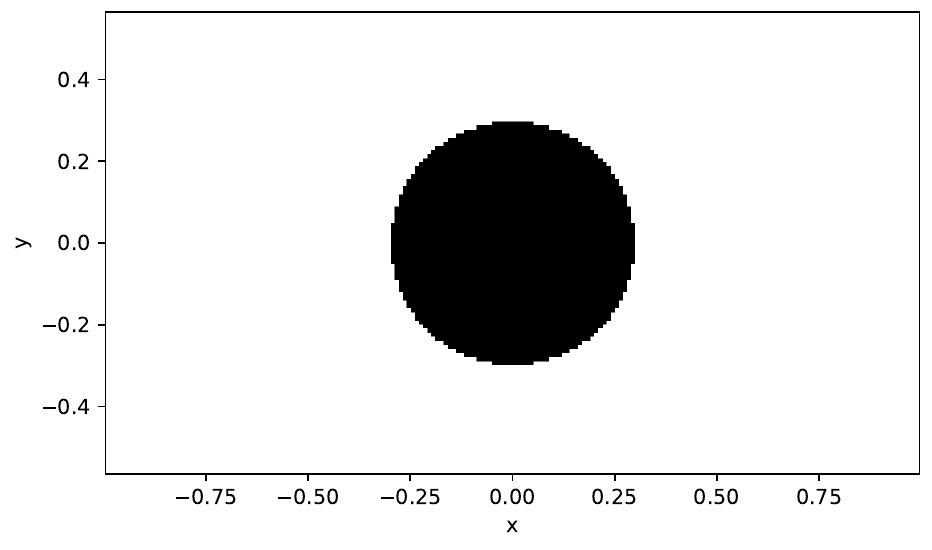}
		\caption{Gray scale plot of $\norm{\theB}^2$}
	\end{subfigure}
	\hfill
	\begin{subfigure}{.45\textwidth}
		\centering
		\includegraphics[width=\linewidth]{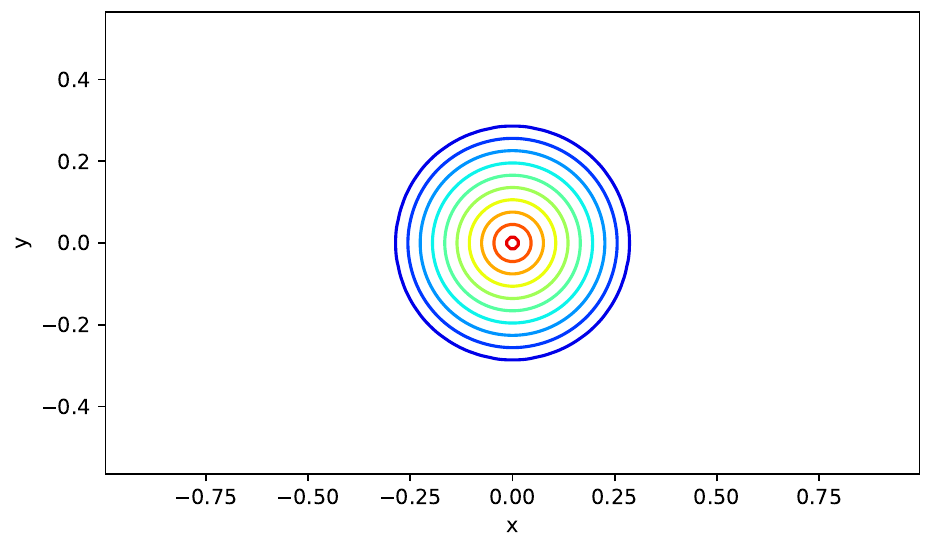}
		\caption{Magnetic field lines}
	\end{subfigure}
	\caption{\cref{ex:field-loop}. Initial magnetic field for the field loop advection problem.}
	\label{fig:field-loop-initial}
\end{figure}
\begin{figure}[htbp]
	\centering
	\begin{subfigure}{.24\textwidth}
		\centering
		\includegraphics[width=\linewidth]{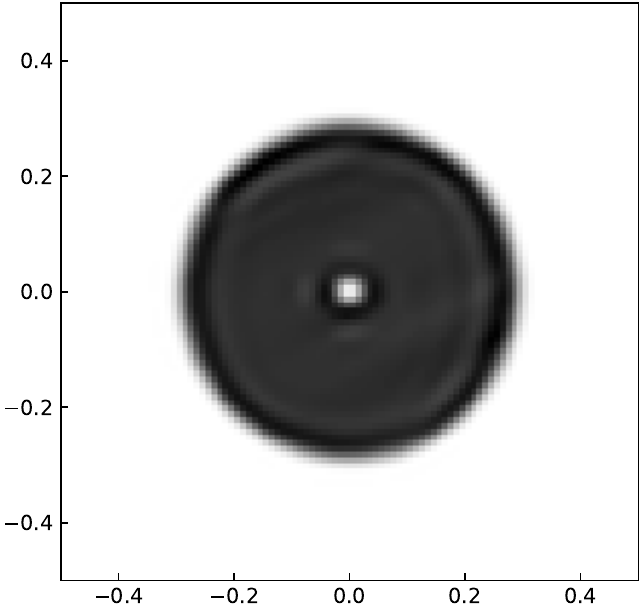}
		\caption{\hwenodf{}, $t=2$}
	\end{subfigure}
	\begin{subfigure}{.24\textwidth}
		\centering
		\includegraphics[width=\linewidth]{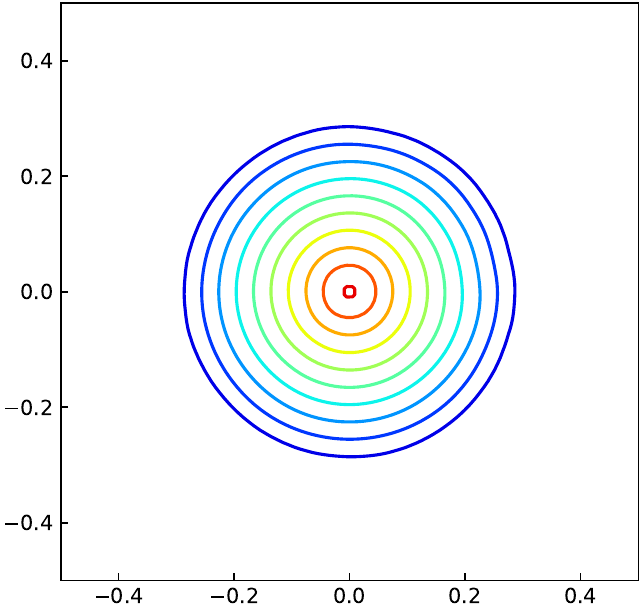}
		\caption{\hwenodf{}, $t=2$}
	\end{subfigure}
	\begin{subfigure}{.24\textwidth}
		\centering
		\includegraphics[width=\linewidth]{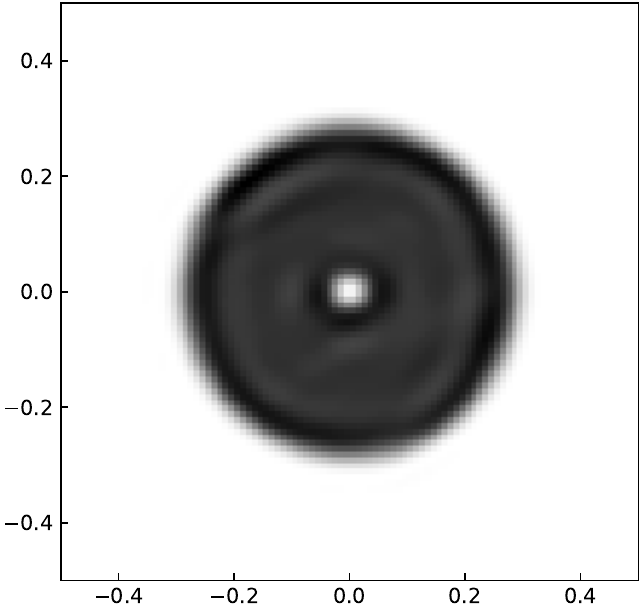}
		\caption{\hwenodf{}, $t=10$}
	\end{subfigure}
	\begin{subfigure}{.24\textwidth}
		\centering
		\includegraphics[width=\linewidth]{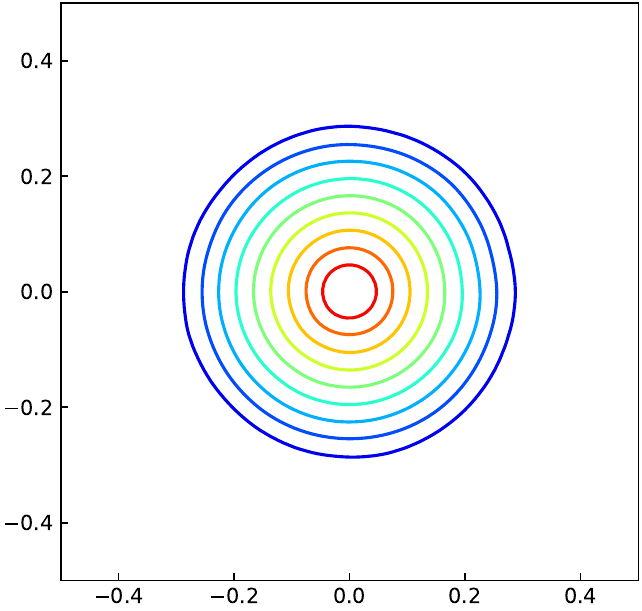}
		\caption{\hwenodf{}, $t=10$}
	\end{subfigure}

	\begin{subfigure}{.24\textwidth}
		\centering
		\includegraphics[width=\linewidth]{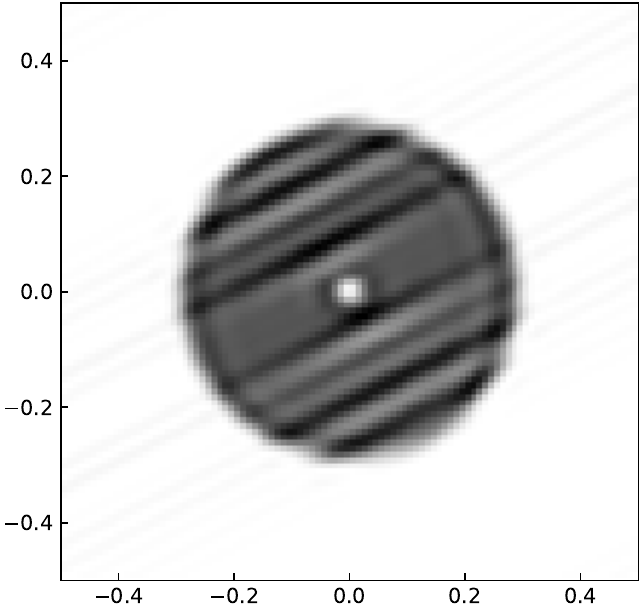}
		\caption{\hwenobase{}, $t=2$}
	\end{subfigure}
	\begin{subfigure}{.24\textwidth}
		\centering
		\includegraphics[width=\linewidth]{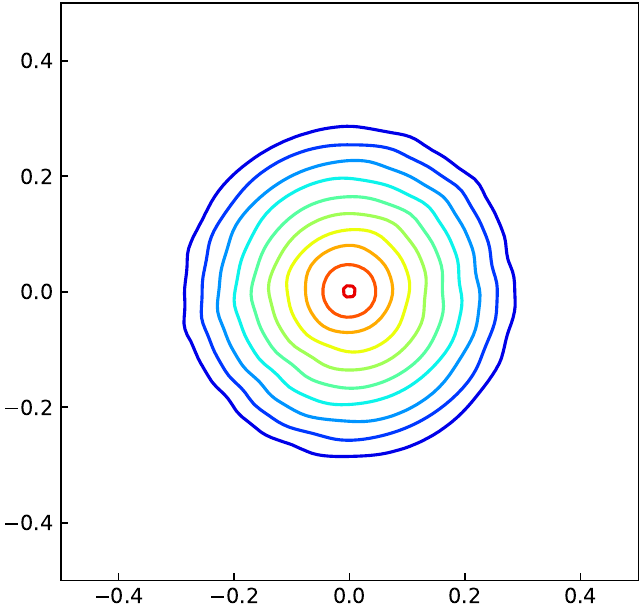}
		\caption{\hwenobase{}, $t=2$}
	\end{subfigure}
	\begin{subfigure}{.24\textwidth}
		\centering
		\includegraphics[width=\linewidth]{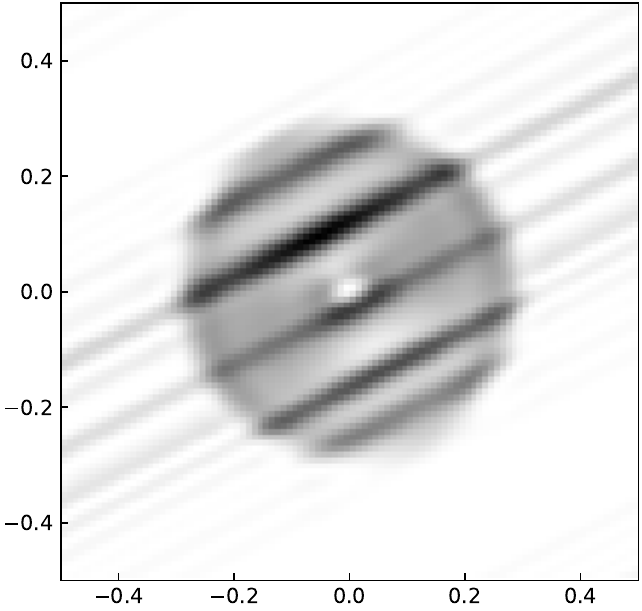}
		\caption{\hwenobase{}, $t=10$}
	\end{subfigure}
	\begin{subfigure}{.24\textwidth}
		\centering
		\includegraphics[width=\linewidth]{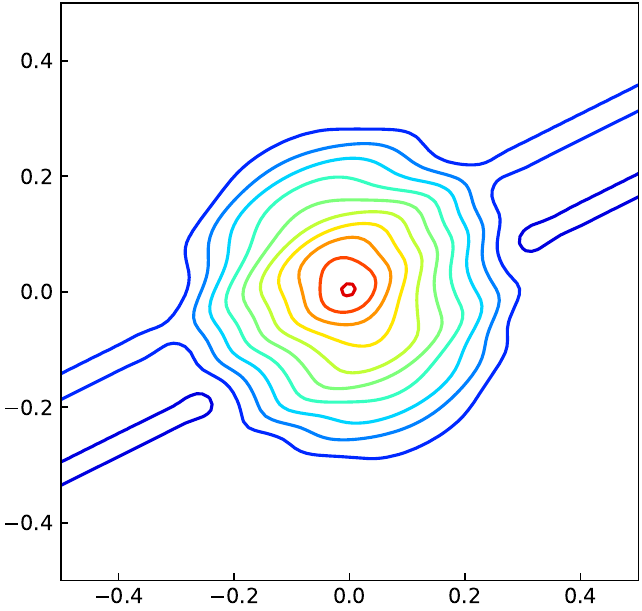}
		\caption{\hwenobase{}, $t=10$}
	\end{subfigure}
	\caption{\cref{ex:field-loop}. Field loop advection, simulated by \hwenodf{} (top row) and \hwenobase{} (bottom row) with $200\times 100$ meshes. From left to right: $t=2$ and $t=10$. In each row, the first two panels show the gray scale plot of $\norm{\theB}^2$ and the magnetic field lines at $t=2$, and the last two panels show the corresponding plots at $t=10$. 9 equally spaced contour lines are plotted for the magnetic field line panels.}
	\label{fig:field-loop-comparison}
\end{figure}
\end{example}

\begin{example}[2D rotated shock tube]
\label{ex:brio-wu}
The Brio--Wu shock tube problem \cite{Brio1988UpwindDifferencing} is a one-dimensional Riemann problem. We rotate the shock tube by an angle $\alpha=\arctan(0.5)$ to obtain a two-dimensional test problem \cite{Toth2000B0Constraint}. The computational domain is $\Omega=[0,1]\times\left[0,\frac{2}{N_x}\cot\alpha\right]$, and a uniform mesh of size $N_x\times N_y=800\times2$ is employed. The left and right boundary values are fixed according to the initial data, while shifted periodic boundary conditions are imposed on the top and bottom boundaries \cite{Toth2000B0Constraint}. Introduce the rotated coordinate $\xi=x+y\tan\alpha=x+\frac{y}{2}$, so that the initial discontinuity is located at $\xi=\frac12$, or equivalently $2x+y=1$. In the local normal--tangential coordinate system, the initial primitive variables are given by
\begin{equation*}
(
\rho,
u_{\parallel},u_{\perp},u_z,B_{\parallel},B_{\perp},B_z,p
)
=
\left\{
\begin{aligned}
&
(
1,0,0,0,0.75,1,0,1
),
\quad && \xi<\frac12,\\[1ex]
&
(
0.125,0,0,0,0.75,-1,0,0.1
),
\quad && \xi\geq\frac12,
\end{aligned}
\right.
\end{equation*}
with an adiabatic index of $\gamma=5/3$. The corresponding Cartesian components are obtained in the same way as in \cref{ex:alfven-wave}. The solution is evolved until $t=0.1\cos\alpha=\frac{0.2}{\sqrt{5}}$. The numerical solutions are transformed back to the local normal--tangential coordinate system and plotted against the rotated coordinate $\xi$. As shown in \cref{fig:brio-wu}, despite reasonable oscillations, the \hwenodf{} scheme preserves the constant parallel magnetic-field component $B_{\parallel}=0.75$ better. In contrast, the \hwenobase{} scheme produces noticeable oscillations near the discontinuities. The results of the \hwenodf{} scheme agree well with those in \cite{Liu2025GloballyDivergencefreea}. This example further demonstrates the effectiveness of the proposed divergence-free correction of magnetic field derivatives in suppressing spurious oscillations and improving the robustness of the scheme.
\begin{figure}[htbp]
\centering
\begin{subfigure}{.32\textwidth}
\centering
\includegraphics[width=\linewidth]
{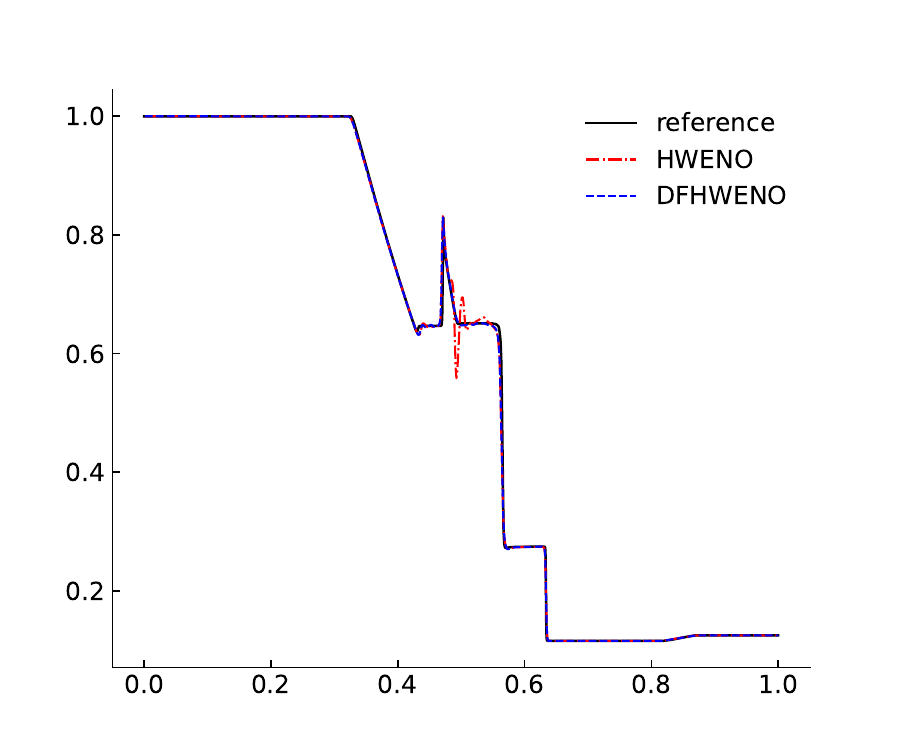}
\caption{$\rho$}
\end{subfigure}
\hfill
\begin{subfigure}{.32\textwidth}
\centering
\includegraphics[width=\linewidth]
{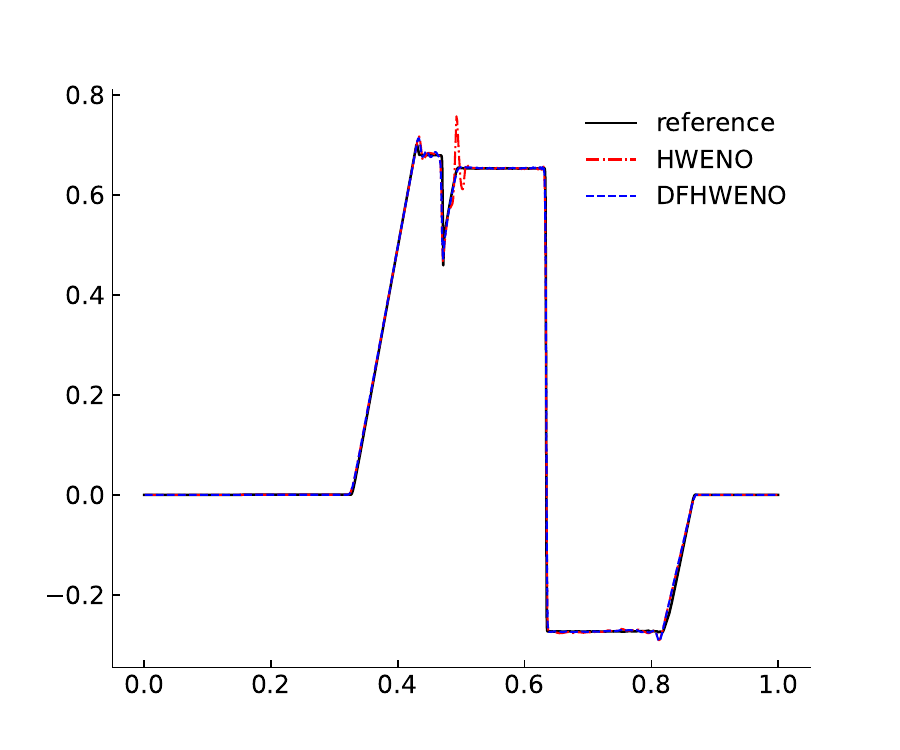}
\caption{$u_{\parallel}$}
\end{subfigure}
\hfill
\begin{subfigure}{.32\textwidth}
\centering
\includegraphics[width=\linewidth]
{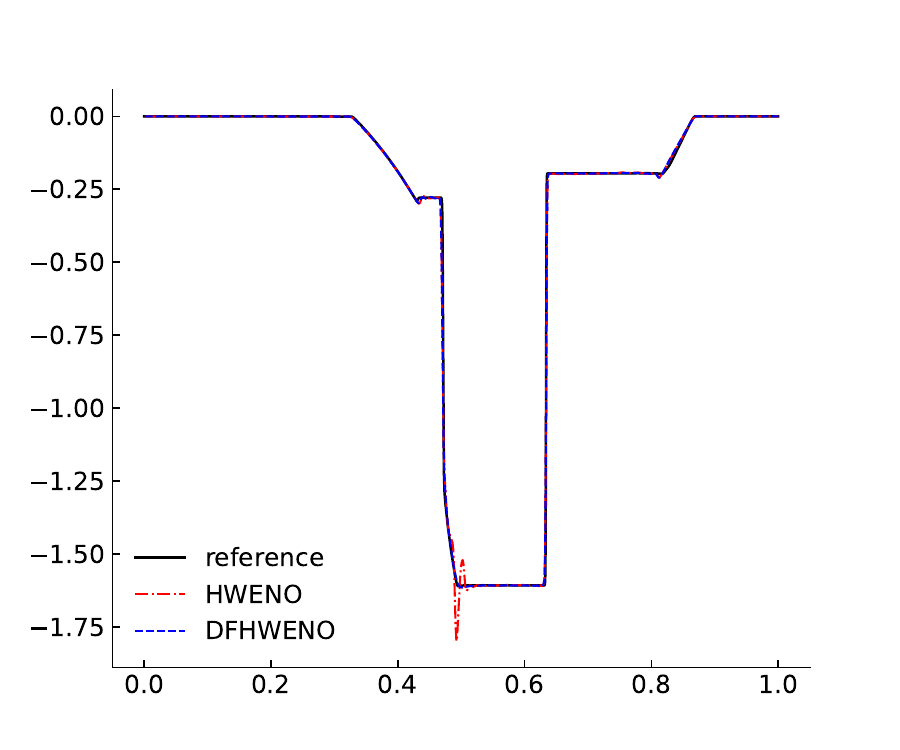}
\caption{$u_{\perp}$}
\end{subfigure}

\begin{subfigure}{.32\textwidth}
\centering
\includegraphics[width=\linewidth]
{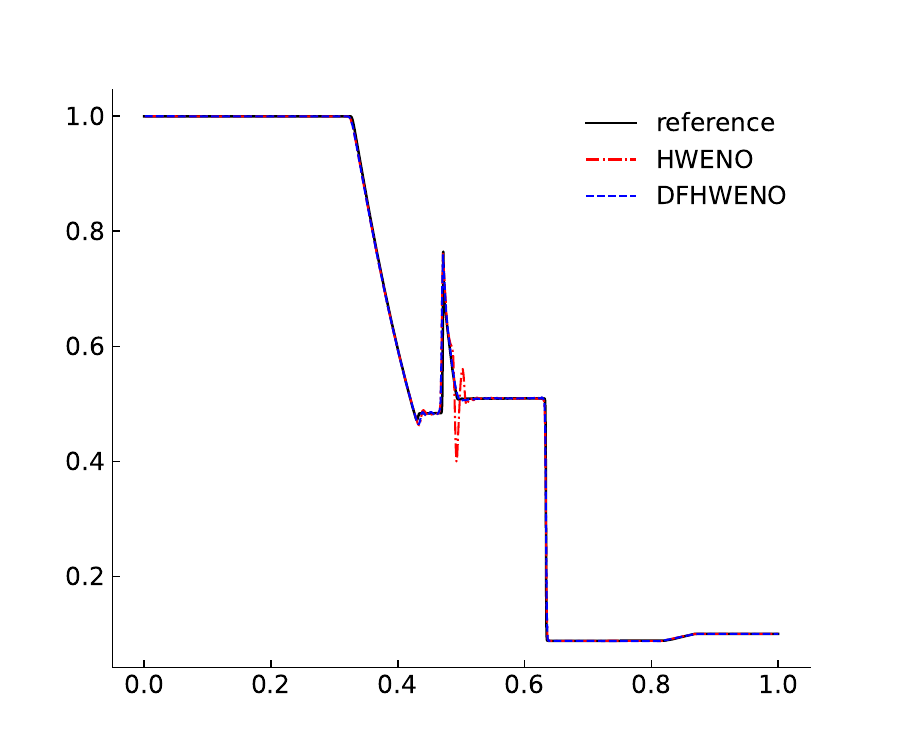}
\caption{$p$}
\end{subfigure}
\hfill
\begin{subfigure}{.32\textwidth}
\centering
\includegraphics[width=\linewidth]
{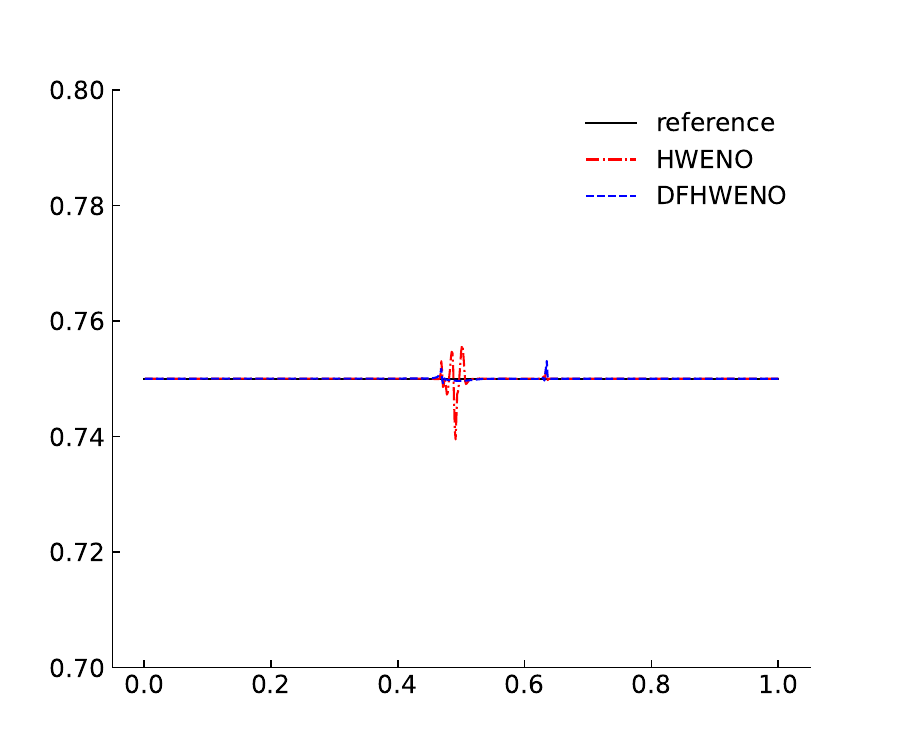}
\caption{$B_{\parallel}$}
\end{subfigure}
\hfill
\begin{subfigure}{.32\textwidth}
\centering
\includegraphics[width=\linewidth]
{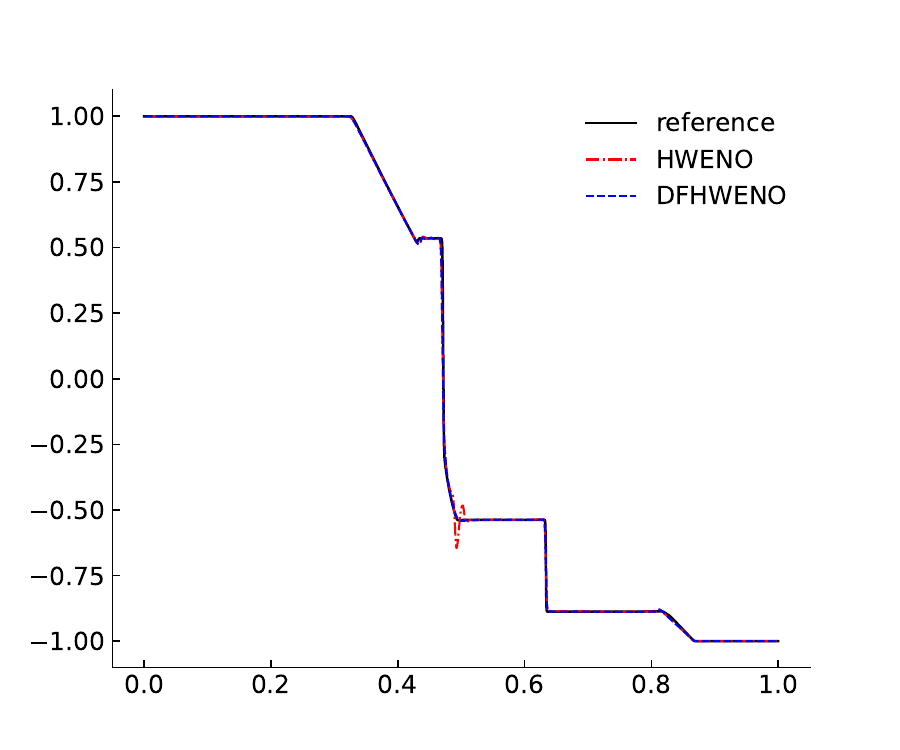}
\caption{$B_{\perp}$}
\end{subfigure}

\caption{\cref{ex:brio-wu}. Rotated Brio--Wu shock tube at $t=0.2/\sqrt{5}$, simulated by \hwenobase{} and \hwenodf{} with $800\times 2$ meshes. The solid line denotes the reference solution simulated by the one-dimensional finite difference WENO scheme \cite{Jiang1999HighOrderWENO}.}
\label{fig:brio-wu}
\end{figure}
\end{example}

\section{Conclusions}
\label{sec:conclusions}
In this paper, we proposed a fifth-order finite difference divergence-free Hermite WENO scheme for the ideal MHD equations.
We constructed a novel discretization of the magnetic field divergence by exploiting the natural incorporation of first-order derivative information in the HWENO framework.
Based on this discrete divergence, we developed a divergence-free correction 
method for handling the divergence-free constraint of the ideal MHD equations. 
The solution and its spatial derivatives are evolved in time by the \hwenobase{} scheme first,
and then the divergence-free correction is applied to the spatial derivatives of the magnetic field. Specifically,
the discrete divergence is computed by $ \divh{\theB} = \theV_{i,j}^{5} + \theW_{i,j}^{6} $ and evenly distributed to the spatial derivatives of the magnetic field $ \theV_{i,j}^{5} $ and $ \theW_{i,j}^{6} $, i.e., $ \hat{\theV}_{i,j}^{5} = \theV_{i,j}^{5} - \frac{1}{2}\divh{\theB} $, $ \hat{\theW}_{i,j}^{6} = \theW_{i,j}^{6} - \frac{1}{2} \divh{\theB} $,  which ensures that the magnetic field is discretely divergence-free at the new time level.

The discrete divergence is
directly eliminated through a simple linear correction procedure applied at every Runge--Kutta stage, which is easy to implement and introduces almost negligible additional computational cost. Moreover, both the conservation property and accuracy  are preserved, since only the spatial derivatives of the magnetic field are corrected by a high-order perturbation, while the conserved quantities remain unchanged.
With this divergence-free correction approach, we apply the fifth-order finite difference HWENO scheme 
to the simulation of ideal MHD for the first time. 
Extensive numerical experiments are conducted. The results show that the proposed scheme is fifth-order accurate and robust, and the proposed divergence-free correction method for handling the divergence-free constraint is crucial and effective. This correction technique can also be extended to the finite volume framework, and the research is  ongoing.

%

%

%

\appendix
\section{Review of the HWENO procedure}
\label{sec:hweno-detail}

This section reviews the one-dimensional HWENO reconstruction of the numerical
fluxes and the limiter for the auxiliary derivative variables.
In multidimensional cases, the HWENO procedure is applied dimension by dimension \cite{Jiang1996EfficientImplementation}.
Without loss of generality, we describe the HWENO procedure in a scalar manner.
We describe in detail the reconstruction of $\hat{f}_{i+\half,j}$ and
$\hat{f}^{v}_{i+\half,j}$. For simplicity, the index $j$ is omitted in the
following description.

\subsection{Reconstruction of the numerical fluxes}
\begin{step}
	\label{step:flux-splitting}
	Apply flux splitting.
	The global Lax--Friedrichs flux splitting is employed
	\begin{equation*}
		f^{\pm}_{i} = \frac{1}{2}(f(u_{i}) \pm \alpha u_{i}),\quad {f}^{v,\pm}_{i} = \frac{1}{2}(f^v(u_{i},v_{i}) \pm \alpha v_{i}),
	\end{equation*}
	where $ \alpha = \max \abs{f'(u)} $.
\end{step}

\begin{step}
	Reconstruct the fluxes.
	Using $f^+_{i-1}$, $f^+_{i}$, $f^+_{i+1}$, and
	$\hat{f}^{v,+}_{i-1}$, $\hat{f}^{v,+}_{i}$,
	$\hat{f}^{v,+}_{i+1}$, we construct three lower-degree
	polynomials $p_{1}(x)$, $p_{2}(x)$, and $p_{3}(x)$, together
	with a higher-degree polynomial $p_{0}(x)$, by requiring them
	to match the corresponding stencil data exactly as given in
	\cref{eq:app-flux-polynomials},
	\begin{equation}
\label{eq:app-flux-polynomials}
\begin{aligned}
\frac{1}{\dx}\int_{I_{i+\ell}} p_k(x) \d{x} &= f_{i+\ell}^{+},
&& \ell \in S_k, \\
\frac{1}{\dx}\int_{I_{i+\ell}} p_k'(x) \d{x} &= f^{v,+}_{i+\ell},
&& \ell \in S_k',
\end{aligned}
\end{equation}
where the stencils $S_k$ and $S'_k$ are given in \cref{tab:app-flux-stencils}.
\begin{table}[htbp]
\centering
\caption{Stencils for the flux reconstruction polynomials.}
\label{tab:app-flux-stencils}
\begin{tabular}{c|cccc}
\toprule
$k$ & 1 & 2 & 3 & 0 \\
\midrule
$S_k$ & $\{-1,0\}$ & $\{0,1\}$ & $\{-1,0,1\}$ & $\{-1,0,1\}$ \\
$S'_k$ & $\{-1,0\}$ & $\{0,1\}$ & $\{0\}$ & $\{-1,0,1\}$ \\
\bottomrule
\end{tabular}
\end{table}
	The linear weights $\gamma_{1,2,3}$ are determined by requiring
	the linear combination of $p_{1,2,3}(x_{i+\half})$ to reproduce
	$p_0(x_{i+\half})$ \cite{Jiang1996EfficientImplementation}.
	Similarly, the linear weights for the derivative reconstruction,
	$\gamma'_{1,2,3}$, are determined by requiring the linear combination
	of $p_{1,2,3}'(x_{i+\half})$ to reproduce
	$p_0'(x_{i+\half})$ \cite{Liu2014FiniteDifference}.

	\begin{step}
		\label{step:flux-reconstruct}
		Reconstruct $\hat{f}_{i+\half}^{+}$
		\cite{Zhao2023WellbalancedFifthorder}
		\begin{equation*}
			f_{i+\half}^{+} = \omega_1 p_1(x_{i+\half}) + \omega_2 p_2(x_{i+\half}) + \omega_3 p_3(x_{i+\half}),
		\end{equation*}
		where the nonlinear weights $\omega_k$ are given by
		\begin{equation*}
			\omega_k = \frac{\alpha_k}{\sum_{l=1}^3 \alpha_l},\quad \alpha_k = \frac{\gamma_k}{(\epsilon + \beta_k)^2}, \quad \beta_k = \sum_{\ell=1}^{\deg p_{k}} \frac{1}{\dx} \int_{I_{i}}^{} (\dx^{\ell} \odv*[\ell]{p_{k}(x)}{x})^{2} \d{x}.
		\end{equation*}
	\end{step}

	\begin{step}
		\label{step:flux-derivative-reconstruct}
		Reconstruct $\hat{f}^{v,+}_{i+\half}$
		\cite{Liu2014FiniteDifference}
		\begin{equation*}
			\hat{f}^{v,+}_{i+\half} = \omega_1 p_1'(x_{i+\half}) + \omega_2 p_2'(x_{i+\half}) + \omega_3 p_3'(x_{i+\half}),
		\end{equation*}
		where the nonlinear weights $\omega_k$ are given by
		\begin{equation*}
			\omega_k = \frac{\alpha_k}{\sum_{l=1}^3 \alpha_l},\quad \alpha_k = \frac{\gamma_k'}{(\epsilon + \beta_k)^2}, \quad \beta_k = \sum_{\ell=2}^{\deg p_{k}} \frac{1}{\dx} \int_{I_{i}}^{} (\dx^{\ell} \odv*[\ell]{p_{k}(x)}{x})^{2} \d{x}.
		\end{equation*}
	\end{step}
\end{step}

\begin{step}
	Finally, combine the positive and negative flux components
	\begin{equation*}
		\hat{f}_{i+\half} = \hat{f}_{i+\half}^{+} + \hat{f}_{i+\half}^{-},\quad \hat{f}^{v}_{i+\half} = \hat{f}^{v,+}_{i+\half} + \hat{f}^{v,-}_{i+\half}.
	\end{equation*}
\end{step}

The negative fluxes $\hat{f}^{-}_{i+\half,j}$ and
$\hat{f}^{v,-}_{i+\half,j}$ can be reconstructed analogously using
right-biased stencils.
The fluxes $\hat{g}_{i,j+\half}$ and
$\hat{g}^{w}_{i,j+\half}$ are reconstructed in the same manner
using stencils in the $y$-direction.
The fluxes involving mixed derivatives are obtained using the following
linear approximations
\begin{equation*}
\begin{aligned}
    \hat{f}^{w}_{i+\half,j} &= -\frac{1}{12} f^w_{i-1,j} + \frac{7}{12} f^w_{i,j} + \frac{7}{12} f^w_{i+1,j} - \frac{1}{12} f^w_{i+2,j},\\
    \hat{g}^{v}_{i,j+\half} &= -\frac{1}{12} g^v_{i,j-1} + \frac{7}{12} g^v_{i,j} + \frac{7}{12} g^v_{i,j+1} - \frac{1}{12} g^v_{i,j+2}.
\end{aligned}
\end{equation*}

\subsection{The HWENO limiter for derivatives}
Then we describe the HWENO limiter for the derivative variables $v$.
Using $ u_{i-1}, u_{i}, u_{i+1} $ and $ v_{i-1}, v_{i}, v_{i+1} $,
we construct three lower-degree
	polynomials $p_{1}(x)$, $p_{2}(x)$, and $p_{3}(x)$, together
	with a higher-degree polynomial $p_{0}(x)$, by requiring them
	to match the corresponding values and derivative values pointwise exactly as given in
	\cref{eq:limiter-polynomials},
	\begin{equation}
\label{eq:limiter-polynomials}
\begin{aligned}
p_k(x_{i+\ell}) &= u_{i+\ell},
&& \ell \in T_k, \\
p_k'(x_{i+\ell}) &= v_{i+\ell},
&& \ell \in T_k',
\end{aligned}
\end{equation}
where the stencils $T_k$ and $T'_k$ are given in \cref{tab:limiter-stencils}.
\begin{table}[htbp]
\centering
\caption{Stencils for the HWENO limiter polynomials.}
\label{tab:limiter-stencils}
\begin{tabular}{c|cccc}
\toprule
$k$ & 1 & 2 & 3 & 0 \\
\midrule
$T_k$ & $\{-1,0\}$ & $\{-1,0,1\}$ & $\{0,1\}$ & $\{-1,0,1\}$ \\
$T'_k$ & $\{-1\}$ & $\emptyset$ & $\{1\}$ & $\{-1,1\}$ \\
\bottomrule
\end{tabular}
\end{table}
	The linear weights $\gamma_{1,2,3}$ are determined by requiring
	the linear combination of $p'_{1,2,3}(x_{i})$ to reproduce
	$p'_0(x_{i})$ \cite{Zhao2020ModifiedFifth}. And the limited derivative $\tilde{v}_{i}$ is reconstructed by
	\begin{equation*}
		\tilde{v}_{i} = \omega_1 p_1'(x_{i}) + \omega_2 p_2'(x_{i}) + \omega_3 p_3'(x_{i}),
	\end{equation*}
where the nonlinear weights $\omega_k$ are given by
	\begin{equation*}
		\omega_k = \frac{\alpha_k}{\sum_{l=1}^3 \alpha_l},\quad \alpha_k = \frac{\gamma_k}{(\epsilon + \beta_k)^2}, \quad \beta_k = \sum_{\ell=1}^{\deg p_{k}} \frac{1}{\dx} \int_{I_{i}}^{} (\dx^{\ell} \odv*[\ell]{p_{k}(x)}{x})^{2} \d{x}.
	\end{equation*}
In two-dimensional cases, the HWENO limiter is applied dimension by dimension \cite{Zhao2020ModifiedFifth}.
\begin{remark}
	For MHD systems, the HWENO reconstruction and limiter described above are
	performed in the characteristic fields \cite{Jiang1999HighOrderWENO}.
\end{remark}

\bibliographystyle{siamplain}
\bibliography{references}
\nocite{*}
\end{document}